\documentclass[draft,11pt]{jsg}
 \usepackage{pdfsync}
 \usepackage{amssymb}
 \usepackage{latexsym}
 \usepackage{amsmath,amsfonts}
 \usepackage{amscd}
 \usepackage[mathscr]{eucal}
\DeclareFontFamily{OT1}{pzc}{}
 \DeclareFontShape{OT1}{pzc}{m}{it}{<-> s * [1.050] pzcmi7t}{}
\DeclareMathAlphabet{\mathpzc}{OT1}{pzc}{m}{it}

 \usepackage{xy}
 \xyoption{all}
 \CompileMatrices
\usepackage{graphicx}
\usepackage{calc}

\newlength{\wcwidth}
\newlength{\wcheight}
\newcommand{\widecheck}[1]{\ensuremath{
\settowidth{\wcwidth}{#1}
\settoheight{\wcheight}{#1}
\addtolength{\wcheight}{1pt}
\makebox[0cm][l]{%
\raisebox{\depth+\wcheight}[0cm][0cm]{%
\scalebox{-1}{$\widehat{\hphantom{#1}}$}}}#1
\rule{0pt}{\wcheight+2.5pt}}}

 \renewcommand{\theequation}{{\thesection.\arabic{equation}}}
 \newcommand{\p}{\partial}
 
\newcommand{\pt}{\textup{p}}

\theoremstyle{plain}
\newtheorem*{thm*}{Theorem}
\newtheorem{thm}{Theorem}[section]\newtheorem{prop}[thm]{Proposition}
\newtheorem{lemma}[thm]{Lemma}\newtheorem{cor}[thm]{Corollary}

\theoremstyle{definition}
\newtheorem{defn}[thm]{Definition}

\theoremstyle{remark}

 \newcommand{\bth}{\begin{*thm}}
 \newcommand{\ethm}{\end{*thm}}
 \newcommand{\bco}{\begin{*cor}}
 \newcommand{\eco}{\end{*cor}}
 \newcommand{\bcj}{\begin{*conj}}
 \newcommand{\ecj}{\end{*conj}}
 \newcommand{\bpr}{\begin{*prop}}
 \newcommand{\epr}{\end{*prop}}
 \newcommand{\bprs}{\begin{**prop}}
 \newcommand{\eprs}{\end{**prop}}
 \newcommand{\ble}{\begin{*lemma}}
 \newcommand{\ele}{\end{*lemma}}
 \newcommand{\bles}{\begin{**lemma}}
 \newcommand{\eles}{\end{**lemma}}
 \newcommand{\bsl}{\begin{*sublemma}}
 \newcommand{\esl}{\end{*sublemma}}
 \newcommand{\bre}{\begin{*rem}}
 \newcommand{\ere}{\end{*rem}}
 \newcommand{\bres}{\begin{**rem}}
 \newcommand{\eres}{\end{**rem}}
 \newcommand{\bnt}{\begin{*nota}}
 \newcommand{\ent}{\end{*nota}}
 \newcommand{\bnts}{\begin{**nota}}
 \newcommand{\ents}{\end{**nota}}
 \newcommand{\bde}{\begin{*defn}}
 \newcommand{\ede}{\end{*defn}}
 \newcommand{\bdes}{\begin{**defn}}
 \newcommand{\edes}{\end{**defn}}

\DeclareMathAlphabet\eusm{U}{eus}{m}{n}
\def\makebb#1{\expandafter\def\csname bb#1\endcsname{{\mathbb{#1}}}\ignorespaces}
\def\makerm#1{\expandafter\def\csname rm#1\endcsname{{\rm #1}}\ignorespaces}
\def\makebf#1{\expandafter\def\csname bf#1\endcsname{{\bf #1}}\ignorespaces}
\def\makegr#1{\expandafter\def\csname gr#1\endcsname{{\mathfrak{#1}}}\ignorespaces}
\def\makescr#1{\expandafter\def\csname scr#1\endcsname{{\mathscr{#1}}}\ignorespaces}
\def\makecal#1{\expandafter\def\csname cal#1\endcsname{{\cal #1}}\ignorespaces}
\def\makeudl#1{\expandafter\def\csname udl#1\endcsname{{\underline{#1}}}\ignorespaces}
\def\doLetters#1{%
  #1A #1B #1C #1D #1E #1F #1G #1H #1I #1J #1K #1L #1M
  #1N #1O #1P #1Q #1R #1S #1T #1U #1V #1W #1X #1Y #1Z}
\def\doletters#1{%
  #1a #1b #1c #1d #1e #1f #1g #1h #1i #1j #1k #1l #1m
  #1n #1o #1p #1q #1r #1s #1t #1u #1v #1w #1x #1y #1z}
\doLetters\makebb\doLetters\makecal\doLetters\makerm\doLetters\makebf\doLetters\makescr
\doletters\makerm \doletters\makebf \doLetters\makegr
\doletters\makegr \doLetters\makeudl \doletters\makeudl

 \newcommand{\op}{\operatorname}
\def\co{\colon\thinspace}
 \def\pf{\noindent{\it Proof.}\enspace}

 \def\Hom{\mathop{\mathrm{Hom}}\nolimits}

 \def\rrr{\textsc{r}}
\def\uuu{\op{\textsf{u}}}
\def\hhh{\op{\mathpzc{h}}}
\def\zzz{\textsl{z}}
\def\smE{\textsc{e}}
\def\smm{\textsc{m}}
\def\ff{\mathpzc{f}}
\def\upN{\textup{N}}
\def\hata{\op{\hat{a}}}
 
 \def\End{\mathop{\mathrm{End}}\nolimits}
 
 \def\ker{\mathop{\mathrm{Ker}}\nolimits}
 \def\cok{\mathop{\mathrm{Coker}}\nolimits}

 \def\C{{\mathbb C}}
 
 \def\U{\mathop{\mathrm{U}}\nolimits}
 
 \def\SO{\mathop{\mathrm{SO}}\nolimits}

 \def\Z{{\mathbb Z}}
 \def\dist{\mathop{\mathrm{dist}}\nolimits}
 \def\cl{\mathop{\mathrm{cl}}}
 
 \def\hatc{\hat{\mathop{\mathrm{c}}}}
 
 \def\R{{\mathbb R}}
 \def\Spin{\mathop{\mathrm{Spin}}\nolimits}

\def\BTitem#1\ETitem{\begin{equation}\hss\left\{\mkern-20mu\parbox{0.9\hsize}{%
\begin{itemize}#1\end{itemize}}\right.\end{equation}}

\begin{document}

 \title[Periodic Floer homology\(=\) Seiberg-Witten Floer cohomology]{Periodic Floer homology and \\Seiberg-Witten Floer cohomology}

 \author{Yi-Jen Lee} 
 \address{Department of Mathematics\\Purdue University\\
          W. Lafayette, IN47907\\USA}
 \email{yjlee@math.purdue.edu}
 
 \author{Clifford Henry Taubes}
 \address{Department of Mathematics\\Harvard University\\
          Cambridge, MA 02138\\USA}
 \email{chtaubes@math.harvard.edu}

 \begin{abstract}
Various Seiberg-Witten Floer
 cohomologies are defined for  a closed, oriented 3-manifold; and if
 it is the mapping torus of an area-preserving surface automorphism, it
 has an  associated periodic Floer homology as defined by Michael
 Hutchings.  We construct
 an isomorphism between a certain
 version of  
 Seiberg-Witten Floer cohomology and the corresponding 
periodic Floer homology, and describe some immediate consequences.

 \end{abstract}

 \maketitle

 \section{Introduction}\label{sec:1}

Suppose that \(F\) is a closed, oriented 2-manifold with an area form and
 a given volume preserving diffeomorphism \(\ff\co F\to F\).  Let \(M\)
 denote the 3-manifold \((\R \times F)/\Z\) where the \(\Z\)-action has
 1 sending any given \((t, x)\in \R \times F\) to \((t + 2\pi, \, \ff(x))\).  Michael
Hutchings \cite{H1}, \cite{HS} defined a version of Floer homology in this
context which he called \textit{periodic Floer homology}.  To a first
approximation, the chain complex is generated by sets of pairs where
each pair consists of an irreducible periodic orbit of \(\ff\) and a
positive integer.  The differential is defined using pseudoholomorphic
curves in \(\bbR \times M\) as defined by an appropriately chosen,
\(\bbR\)-invariant almost complex structure.  Hutchings conjectured
that this periodic Floer homology is isomorphic to a version of the
Seiberg-Witten Floer cohomology.  The purpose of this article is to
explain how the analysis used in \cite{T1}-\cite{T4} to
establish the equivalence between Seiberg-Witten Floer cohomology and embedded
contact homology can be used to prove Hutchings's conjectured
equivalence between periodic Floer homology to Seiberg-Witten Floer
cohomology.  A precise statement of the equivalence is given in Section
\ref{sec:1(c)} for the monotone case. The general version of the
isomorphism theorem requires additional preparation to state, and is
therefore postponed until Section \ref{sec:6(c)}.  What
follows directly sets the stage with a more detailed description of
periodic Floer homology and of the relevant version of Seiberg-Witten
Floer cohomology.

 \subsection{Periodic Floer homology: the monotone case}\label{sec:1(a)}

The description that follows of periodic Floer homology for the most
part paraphrases what is presented in \cite{H1} and \cite{HS}.

\subsubsection*{The Geometry of \(M\):}
The projection map from \(\bbR  \times  F\) to \(\bbR  \) descends to \(M\) so as
to give a fibration \(\pi : M \to  S^1\).  Here, \(S^1\) is identified with
\(\bbR /(2\pi \bbZ )\). 
Conversely, a closed, oriented 3-manifold that
fibers over \(S^1\) admits a fiber preserving diffeomorphism to the
manifold \((\bbR  \times  F)/\bbZ  \) where \(F\) here is the fiber of the
fibration, and where the \(\bbZ  \)-action is defined by the
area-preserving diffeomorphism, \(\ff\), of \(F\).  

The push forward of the Euclidean vector field on \(\bbR  \) defines a
nowhere zero vector field, \(\partial_t\), which is transverse to the
fibers of \(\pi \).  Closed integral curves of \(\partial_t\) are determined
by the fixed points of the iterates of \(\ff\).  To elaborate, introduce the
notion of a \textit{periodic orbit} of \(\ff\).  This is a finite,
cyclically ordered set of distinct points in \(F\) that are
cyclically permuted by \(\ff\).  Let \(\gamma  \subset  F\) denote such a set.
The number of points in \(\gamma  \) is said to be the {\em period}.  The
periodic orbits of \(\ff\) are in 1-1 correspondence with the compact,
embedded 1-dimensional submanifolds that are integral curves of
\(\partial_t\).  The symbol \(\gamma  \) is used in what follows to denote both
a periodic orbit and also the associated embedded integral curve of
\(\partial_t\). 

A periodic orbit of period \(q\) is said to be {\em non-degenerate} when the
linear map \(1 - d(\ff^{kq})\) on \(TF|_{\gamma }\) is invertible for all
integers \(k > 0\).  Here, \(\ff^{m}\) is shorthand for the \(m\)'th iterate of the
map \(\ff\).  We assume in what follows that 
\begin{equation}\label{f-assumption}
\textit{\(\ff\) is chosen so that all orbits are
non-degenerate.}  
\end{equation}
An orbit is said to be \textit{hyperbolic} if \(d\ff^{q}\)
has real eigenvalues.  Otherwise the orbit is said to be
\textit{elliptic}.  The hyperbolic orbit is \textit{positive} or
\textit{negative} when the eigenvalues of \(d\ff^{q}\) are respectively
positive or negative.

Any given closed integral curve of \(\partial_t\) has a canonical
orientation, namely the orientation along \(\partial_t\).  This understood,
such an integral curve defines a class in \(H_{1}(M; \bbZ )\).  If \(\gamma\)
denotes the integral curve, then \([\gamma ]\) will denote the
corresponding homology class.  Meanwhile, the fiber \(F\) of \(\pi  \)
is oriented by its area form and so defines a class, \([F]\), in \(H_{2}(M;
\bbZ )\).  The Poincar\'e dual of this class pairs with \([\gamma ]\) to
give the period of the corresponding periodic orbit of \(\ff\).  

Two other cohomology classes play a distinguished role in what
follows.  To describe the first of these, introduce \(w _{F}\) to
denote the area form of the surface \(F\).  Since \(\ff\) preserves \(w _{F}\),
this form descends to \(M\) as a closed, nowhere zero 2-form on \(M\) that
annihilates the vector field \(\partial_t\).  We denote this 2-form on
\(M\) by \(w_\ff\); and use \([w _{\ff}]\)
to denote its cohomology class in \(H^{2}(M; \bbR )\).  The
second is the Euler class of the kernel of \(\pi \)'s differential.  This
class is denoted as \(c_{1}(K^{-1})\).  
Note that while the class \(c_1(K^{-1})\) remains fixed
under symplectic isotopies of \(\ff\), the class \([w_\ff]\) varies in
the following manner: By the Mayer-Vietoris sequence, we have the
following exact sequence: 
\[ 0\to \big(\ker (1-\ff_*)|_{H_1(F)}\big)^*\to H^2(M;\bbR)\to \big(
\bbR[F] \big)^*\to 0.\] 
By the definition of \(w_f\), the cohomology class \([w_f]\)  belongs
to the coset of elements in \(H^2(M)\) which maps to \([w_F]\in H^2(F)\).
As \(\ff\) varies through a 
 symplectic isotopy, the class \([w_\ff]\) varies by the image of the
 flux of this symplectic isotopy under the composite map: \(H^1(F)\to \big(\ker (1-\ff_*)|_{H_1(F)}\big)^*\to H^2(M;\bbR)\). Thus, \([w_\ff]\) is invariant under Hamiltonian
 isotopies of \(\ff\), and, because of the surjectivity of the flux
 homomorphism, it can take values in any element in the half space
\(\{\, e\, |\, \langle e, [F]\rangle>0\}\subset H^2(M;\bbR)\) under
symplectic isotopies and rescaling.  

\begin{defn}
Given a class \(\Gamma  \in
H_{1}(M;\bbZ)\), Let \(e_\Gamma\) denote the Poincar\'e dual of
\(\Gamma\), and let 
\[
c_\Gamma:=2e_\Gamma +c_1(K^{-1}).
\]
The class \(\Gamma\) is said to be \textit{monotone} (with respect to
\([w_\ff]\)) if 
\begin{equation}\label{def:monotone}
[w _{\ff}]=-\lambda \, c_\Gamma \quad \text{for some \(\lambda\in \bbR\)}. 
\end{equation}
Note that since \([w_\ff]\) is nontrivial, by this definition
\(c_\Gamma\) can not be torsion if \(\Gamma\) is
monotone, and \(\lambda\neq 0\).  We call \(\Gamma\) {\em positive
  monotone} when \(\lambda>0\); conversely, {\em negative monotone} if \(\lambda<0\).
\end{defn}
Fix a monotone \(\Gamma\) throughout this subsection. 

Let
\[
d_\Gamma:=e_\Gamma  \, ([F]).
\]
Note that
since \(c_\Gamma\) is non-torsion,
\(d_\Gamma\neq \textit{genus}\, (F)-1\). This condition guarantees that any given map \(\ff\) can be isotoped via a symplectic isotopy to obtain monotonicity. (A generic map \(\ff\) will not have any monotone classes).

\subsubsection*{The chain complex:}  
Introduce the set \(\mathcal{A}\) whose elements are described next. A given element, \(\Theta \), is a finite set of pairs, each of the form \((\gamma ,
m)\) where \(\gamma  \) is a periodic orbit of \(\ff\) and where \(m\) is a positive
integer.  These are constrained as follows:
\BTitem\label{(1.1)}
 \item  \(\sum_{(\gamma ,m) \in   \Theta } m [\gamma] = \Gamma \).
\item  Distinct pairs have distinct periodic orbit
components.
\item  The integer \(m = 1\) when \(\gamma\) is hyperbolic.
\ETitem

Note that \(\Theta\) can be the empty set, and that \(\mathcal{A}\) is
a finite set.  Indeed, such is the case because
the integer \(d_\Gamma \) is the sum over each \((\gamma,m)\) in  
\(\Theta\) of \(m\) times the period of \(\gamma\).
Meanwhile, there are only
finitely many periodic points with period bounded by any given
integer.  Let \(\mathcal{A} _{+}\) denote the set of pairs of the form \((\Theta ,
\gro)\) where \(\Theta  \in  \mathcal{A}\)  and where \(\gro \) denotes an ordering of
the elements in \(\Theta  \) whose periodic orbit component is positive and
hyperbolic.  The chain complex for periodic Floer homology is defined
to be
\begin{equation}\label{CP}
CP_*(\ff, \Gamma) = \bbZ \mathcal{A}_+/\sim,
\end{equation}
where \(\sim \) denotes the equivalence relation that that has
\((\Theta , \gro) \sim  -(\Theta , \gro')\) when \(\gro'\) differs
from \(\gro \) by an odd permutation.  This is a free \(\bbZ
\)-module.  One can assign to each
\(\Theta  \in  \mathcal{A}\)  a fixed ordering, \(\gro\), of the elements in \(\Theta\)
whose periodic orbit component is positive and hyperbolic.  Doing so
identifies the chain complex \(CP_*\) with \(\bbZ \mathcal{A}\).  Such an assignment is
assumed implicitly in what follows unless stated to the contrary.
This the case, each \(\Theta  \in  \mathcal{A}\) can be viewed as a generator of
the chain complex.  
\bigbreak

\subsubsection*{A relative grading:}  The chain
complex has a relative \(\bbZ /p\bbZ  \)-grading where \(p\) is the
divisibility in \(H^{2}(M; \bbZ )/\op{Tors}\) of the class \(c_\Gamma\).
 To define this grading, suppose that
\(\Theta  _+\) and \(\Theta _-\) are two elements from \(\mathcal{A}\).  Use \(H_{2}(M; \Theta_+
, \Theta _-)\) to denote the set of relative homology classes of
integer-valued 2-chains \(Z \subset  M\) with 
\[
\partial Z = \sum_{(\gamma,m) \in   \Theta   _+} m \gamma-\sum_{(\gamma
  ,m)\in  \Theta_-} m \gamma .
\]  
This is a torsor modeled on \(H_{2}(M; \bbZ )\).  Any
given \(Z \in  H_{2}(M; \Theta _+, \Theta _-)\) can be used to associate an
integer, \(I(\Theta _+, \Theta _-; Z)\), to the pair \((\Theta_+ , \Theta _-)\).
This integer is defined in Section 2 of \cite{H1}.  As explained in the
latter reference, this integer has two nice features.  
\begin{itemize}
\item[(1)] \(I(\Theta_+
, \Theta _-; Z)\) changes by a multiple of \(p\) when the class \(Z\)
is changed. More precisely,
\begin{equation}\label{spec-flow}
I(\Theta_+, \Theta _-; Z)-I(\Theta_+, \Theta _-; Z')=\langle c_\Gamma, [Z-Z']\rangle.
\end{equation}
\item [(2)] \(
I(\Theta , \Theta '; Z) +I(\Theta ', \Theta ''; Z') = I(\Theta ,\Theta
''; Z + Z').\)  
\end{itemize}
These two features imply that \[
I(\Theta_+, \Theta_-):=I(\Theta_+ , \Theta _-; \cdot ) \mod  p \bbZ \]
gives a relative \(\Z/p\bbZ \)-grading to the chain complex.

Some digression is required before introducing the differential.

\bigbreak

\subsubsection*{The set of complex structures \(\mathcal{J}_f\):} Let
\(\mathcal{J}_f\) be the set of almost complex structures \(J\) on
\(\R\times M\) satisfying the following properties:
\BTitem\label{def:J}
\item[\it 1)]  \(J\) is invariant with respect to the
action of \(\bbR  \) on \(\bbR  \times  M\) that translates a constant amount
along the \(\bbR  \) factor. 
\item[\it 2)] Let \(s\)  denote the Euclidean coordinate along the \(\bbR  \) factor in \(\bbR  \times
M\). Then \(J\cdot \partial _{s} = \partial_t\).
\item[\it 3)] Let \(\gamma  \) denote any periodic orbit from a set in \(\mathcal{A}\) .  View \(\gamma  \) as an embedded
circle in \(M\).  Then \(J|_{\gamma}\) maps \(\ker\,(\pi _*)\) to
itself.
\item[\it 4)]  \(J\) defines a \textit{tame} almost
complex structure with respect to the symplectic form 
\[
\Omega _{F} = ds \wedge  dt + w _{\ff}
\] 
on \(\bbR  \times  M\).  That is, the quadratic form \(\Omega _{F}(\cdot , J(\cdot ))\) on \(T(\bbR  \times  M)\) is
positive. 
\ETitem
Endow \(\mathcal{J}_f\) with the \(C^{\infty }\)-Frechet topology. 
We say that \(J\) is {\em admissible} if condition {\it 3)} above is
replaced by the stronger condition:
\begin{itemize}
\item[\(\textit{3}'\)\it )] \(J\) maps \(\ker \pi_*\) to itself at any point on \(M\).
\end{itemize}

\subsubsection*{Pseudoholomorphic subvarieties in \(\R\times M\):}
For the purposes of this article, a {\em pseudoholomorphic subvariety}
in \(\bbR  \times  M\) is a closed subset \(C  \subset  \bbR  \times  M\) with the following properties:

\begin{itemize}
\item \(C\) has no point components. 
\item  The complement in \(C  \) of a
finite set of points is a smooth, non-empty, 2-dimensional submanifold
of \(\bbR  \times  M\) whose tangent space is \(J\)-invariant.  Endow
\(C\) with the complex structure and orientation induced from \(J\)
away from these points. 
\item  The integral over \(C \) of \(w _{\ff}\) is finite.
\end{itemize}
\subsubsection*{\it Remark.}
The third item of this definition is equivalent to the
requirement that \(C\) has finite engery with respect to the
symplectic form \(\Omega_F\) in the sense of Hofer {\it et al.}
(Cf. e.g. Section 5.3 in \cite{BEHWZ}). Indeed, the integral \(\int_C
w_\ff\) is what is termed ``\(\omega\)-energy'' in the aforementioned
reference. On the other hand, the integral \(\int_{C\cap \{s\}\times M} dt \) is 
independent of \(s\) whenever it is well-defined, since slices of
\(C\) at any \(s\in \bbR\) are all homologous in \(M\). Denote this
constant by \(l\). (In our context, which will be described
momentarily, \(l=d_\Gamma\)). Then 
\[
\int_{C\cap [s-L, s+L]\times M} ds\wedge dt= 2L\, l \quad \forall s\in \bbR.
\]  
This implies that the ``\(\lambda\)-energy'' in \cite{BEHWZ} is bounded automatically.

To continue the digression, 
if \(\gamma  \subset  M\) is a closed integral curve of \(\partial_t\), then
the cylinder \(\bbR  \times  \gamma  \subset  \bbR  \times  M\) is
pseudoholomorphic.  Cylinders of this sort are the only
\(\bbR\)-invariant, irreducible pseudoholomorphic subvarieties.  (A
pseudoholomorphic subvariety is said to be irreducible when no finite
set has disconnected complement).

Suppose that \(C\) is a pseudoholomorphic
subvariety.  The condition that \(w _{\ff}\) have finite integral on \(C\)
can be used to show the following:  There exists \(s_{0} \in  \bbR  \)
such that the \(|s| \geq  s_{0}\) portion of \(C\) is a disjoint union of embedded
submanifolds on which \(s\) restricts with no critical points.  In
particular, each such submanifold is a cylinder that lies in a small
radius tubular neighborhood of either the \(s\geq s_0\) or \(s\leq -s_0\) part of one of the \(\bbR\)-invariant, pseudoholomorphic cylinders, thus a in a neighborhood
of \(\bbR\times \gamma\) with \(\gamma\subset M\) being a closed, integral curve of \(\partial_t\).  These
cylinders are called the \textit{ends }of \(C\).  An end where \(s\) is
unbounded from below is said to be a \textit{negative} end, and one
where s is unbounded from above is said to be a \textit{positive} end.  

The constant-\(s\) slices of any given end of \(C\) limit as \(|s| \to  \infty
\) to a closed integral curve of \(\partial_t\).  To elaborate, let \(\gamma  \in  M\)
denote such a curve.  The non-degeneracy of \(\gamma\) implies the
following: There exists a disk \(D \subset  \bbC  \) about the
origin and an embedding \(\varphi_\gamma : S^1 \times  D \to  M\)
such that \(\varphi_\gamma (\cdot , 0) = \gamma \), with the 
additional properties described below:  
Suppose first that \(\mathpzc{E} \subset C\) is a positive end
whose constant \(s\) slices limit to \(\gamma  \) as \(s \to  \infty
\).   Let \(\hat{\varphi}_\gamma\) denote the map \(\op{id}\times
\varphi_\gamma: \R\times (S^1\times D) \to \R \times M\). Then
there exists \(s_{\mathpzc{E}} > 0\) and an integer \(q_{\mathpzc{E}} \geq  1\), such that
\(\hat{\varphi}_\gamma^{-1}(\mathpzc{E} \, \cap \, ([s_{\mathpzc{E}}, \infty )
\times M))\) 
is the image of a map from \([s_{\mathpzc{E}}, \infty ) \times  \bbR /2\pi
q_{\mathpzc{E} }\bbZ  \) to \(\bbR  \times  (S^1\times  D)\)  that has
the form \((s, \tau) \mapsto (s, (\tau, \varsigma(\tau, s))\), 
where \(\varsigma\) can be written as
\begin{equation}\label{(1.3)}
 \varsigma (\tau, s) =
 e^{-\lambda_{q_{\mathpzc{E}}}s}\varsigma_{q_{\mathpzc{E}}}(\tau) (1+ \grr (\tau, s))
\end{equation}
with \(|\grr (\tau, s)| \leq     e^{-\varepsilon |s|}\) for some \(\varepsilon  = \varepsilon_{\mathpzc{E}} \geq  0\).  Here,
\(\lambda _{q_{\mathpzc{E}}}\) is a positive constant.  If \(\mathpzc{E} \subset  C\) is a negative end, then
(\ref{(1.3)}) holds when \(s\leq  -s_{\mathpzc{E}} \) with
\(\lambda _{q_{\mathpzc{E}}}\) a negative constant.  The integer
\(q_{\mathpzc{E}}\) is said to be the {\em multiplicity}
of the end \(\mathpzc{E}\).  The function \(\varsigma_{q_{\mathpzc{E}}}(\tau)\)
can be assumed to be a \((2\pi q_{\mathpzc{E}})\)-periodic eigenfunction of a
symmetric, \(\bbR \)-linear operator  
\begin{equation}\label{(1.4)}
P_\gamma \co C^{\infty}(\bbR ; \bbC )\to C^{\infty}(\bbR ; \bbC ),
\quad \eta  \mapsto\frac{i}{2} \frac{d}{d\tau} \eta+ \nu\eta + \mu \bar{\eta}.
\end{equation}
Here, \(\nu  \) is a real-valued function with period \((2\pi  )\), and
\(\mu  \) is a \(\bbC \)-valued periodic function with period \((2\pi  )\).  They are determined by \(\ff\) and
\(J\).  The constant \(\lambda _{\mathpzc{E}}\) in (\ref{(1.3)}) is the corresponding
eigenvalue.  With regard to terminology:  A function, \(\eta \), on \(\bbR\)
is said to be \((2\pi q)\)-periodic for some positive integer \(q\) if 
\(\eta \, (\cdot  + 2\pi q) = \eta \, (\cdot)\). The {\em period}
of \(\eta\) is the least positive integer for which the above equality
holds.  
The assumption that \(\gamma  \) is non-degenerate
guarantees that the kernel of the operator \(P_\gamma\) 
has no non-trivial elements with period \((2\pi q)\) for any \(q\in \Z^+\).

\subsubsection*{The differential:}  
It follows from (\ref{(1.3)}) that the image via the projection
from \(\bbR  \times  M\) to \(M\) of a pseudoholomorphic
subvariety \(C\) defines a 2-dimensional integer cycle, which we
denote by \([C]\).

Fix elements \(\Theta _+\) and \(\Theta _-\) from the set \(\mathcal{A}\)
that gives the generators for periodic Floer homology.
We define the set \(\mathcal{M}_{1}(\Theta _+, \Theta _-)\) as
follows: If the relative \(\bbZ/p \bbZ\)-grading \(I(\Theta _+,
\Theta _-)\neq 1\), let \(\mathcal{M}_{1}(\Theta _+, \Theta _-)=\emptyset\). Otherwise, if \(I(\Theta _+,
\Theta _-)= 1\mod p\), let \(\mathcal{M}_{1}(\Theta _+, \Theta
_-)\) be the set consisting of elements of the following form:
\begin{itemize}
\item[(1)] An element  \(\Sigma\in \mathcal{M}_{1}(\Theta _+, \Theta
_-)\) is itself a finite set consisting of  pairs
the form \((C, m)\), where \(C\) is an embedded, connected pseudoholomorphic
submanifold in \(\bbR\times M\), and \(m\) is a positive integer.  
The integer \(m=1\) unless \(C=\bbR\times \gamma\), where \(\gamma\)
is an elliptic periodic orbit. 
 \item[(2)] If \((C,m)\), \((C',m')\) are distinct pairs in \(\Sigma\),
   then \(C, C'\) are disjoint.

\item[(3)] 
The weighted pseudoholomorphic curves constituting \(\Sigma\) have the 
appropriate asymptotic behaviors and weights such that 
the cycle \[
Z_{\Sigma } = \sum_{(C,m)\in \Sigma } m \,C
\] 
defines an element in \(H_{2}(M\, ; \Theta _+, \Theta _-)\).  We call
this element the {\em relative homology class} of \(\Sigma\).
\item[(4)] Denote the above relative homology class by 
\(
[Z_\Sigma] \in  H_{2}\, (M\, ; \Theta _+,\Theta _-).
\) 
Then \([Z_\Sigma]\) must be among the
relative homology classes satisfying 
\[ I\, (\Theta _+, \Theta _-, [Z_{\Sigma }])=1\in \bbZ.\]
\end{itemize}

The assumption that \(\Gamma\) is monotone has the following important
consequence: First, note that because of the property (\ref{spec-flow}) of the relative index \(I\),
the set of relative homology classes meeting the constraint (4) above
is a torsor over \(\ker c_\Gamma\), where \(c_\Gamma\) is
viewed as a homomorphism from \(H_2\, (M; \bbZ)\) to \(\bbZ\).
The monotonicity of \(\Gamma\) implies that \([w_\ff]\), also viewed as a
homomorphism from \(H_2\, (M; \bbZ)\) to \(\bbZ\), restricts to a
trivial map from \(\ker c_\Gamma\). That is to say, for every
\(\Sigma\in  \mathcal{M}_{1}(\Theta _+, \Theta _-)\), the integral
\[
\int_\Sigma w_f=\sum_{(C_i,m_i)\in \Sigma}m_i \int_{C_i}w_f, 
\]
has the same value. 
This provides the type of ``energy bound'' needed for
the typical Gromov compactness argument to establish the 
compactness of \(\mathcal{M}_{1}(\Theta _+, \Theta _-)\).
(Cf. e.g. Lemma 9.8 in \cite{H1}). 

With this understood, we shall define in Section \ref{sec:2(b)}
a residual subset \(\mathcal{J} _{1f} \subset  \mathcal{J} _{f}\),
with the following property: 
If \(J \in  \mathcal{J} _{1f}\), then the set \(\mathcal{M}
_{1}(\Theta _+, \Theta _-)\) has the structure of a
smooth, 1-dimensional manifold with a finite set of  components. 
Moreover, each component is a free orbit of
the \(\bbR  \)-action that is induced by the action of \(\bbR  \) on
\( \bbR  \times  M\) by constant translations along the \(\bbR\)
factor.  
As a parenthetical remark, note that
an admissible almost complex structure need not be in \(\mathcal{J} _{1f}\). 

Assume in what follows that \(J\) is from the set \(\mathcal{J} _{1f}\).

Let \(\gro_+\), \(\gro_-\) respectively denote a choice of the ordering of the
 elements in \(\Theta_+\) and \(\Theta_-\) whose periodic orbit component is positive hyperbolic.  
Given a \(J\in \mathcal{J} _{1f}\), each component of \(\mathcal{M}
_{1}(\Theta _+, \Theta _-)\) can be assigned a weight
 from the set \(\{\pm 1\}\) once a choice of the orderings, \(\gro_+\),
 \(\gro_-\) is made. 
 This is described in Section 9 of \cite{HT1}, where similar weights
 are assigned to the pseudoholomorphic subvarieties that
 are used to construct the differential for embedded
 contact homology.  Note that we are implicitly assuming here and in what follows that a product structure has been chosen for a certain natural \(\bbZ/2\)-bundle over the set of all hyperbolic orbits of \(\partial_t\). There is
canonical isomorphism between the respective chain complexes that are defined by
different product structures.

If \(\mathcal{M} _{1}(\Theta _+, \Theta _-)\)
 is not empty, use \(\sigma (\Theta _+, \Theta _-) \in  \bbZ    \) to denote the sum of the \(\pm 1\) weights 
associated to its components with respect to the choice of orderings
\(\gro\) made following the definition of \(CP\) in (\ref{CP}).  Set \(\sigma (\Theta _+,
\Theta _-)\) equal to 0 when \(\mathcal{M} _{1}(\Theta _+, \Theta_-)\) is empty.    

With the digression now over, we note that the differential
for periodic Floer homology is defined by taking the
linear extension of the following action on the
generators:
\begin{equation}
 \partial \Theta_+= \sum_{\Theta _-\in \mathcal{A}}
\sigma (\Theta _+, \Theta _-) \, \Theta  _-.
\label{(1.5)}\end{equation}

This endomorphism of \(\bbZ \mathcal{A}\)  has square zero; a proof can be had by
taking almost verbatim the proof in \cite{HT2} that the corresponding
embedded contact homology differential has square zero.  The
endomorphism \(\partial\) also decreases the \(\bbZ /p\bbZ  \) degree
by 1.  

In summary, we have defined a Floer chain complex with coefficient
\(\bbZ\) to each parameter set:
\begin{equation}\label{PFH-par}
\boldsymbol{\mu}_P =\{(F, w_F), \ff, \Gamma, J\},
\end{equation}
where \((F, w_F)\) is an oriented closed surface with volume form \(w_F\),
\(\ff\) is a nondegenerate volume-preserving automorphism of \(F\), \(\Gamma\) is a
monotone class in \(H_1(M; \bbZ)\) with respective to \(w_\ff\), and
\(J\in \mathcal{J}_{1\ff}\) is an almost complex structure on \(\bbR\times M\).
The
resulting \(\bbZ /p\bbZ  \)-graded homology is by definition the periodic
Floer homology associated to \(\boldsymbol{\mu}_P\), which we denote by
\[
HP_* \big(\ff\co(F, w_F) \circlearrowleft   , \Gamma\big)_J=H (CP_*, \partial),
\] 
or simply \(HP_* (\ff, \Gamma)\) when there is no danger of confusion.

\subsection{Seiberg-Witten Floer cohomology: the monotone case}\label{sec:1(b)}

Some stage setting is needed before describing the
relevant version of the Seiberg-Witten Floer cohomology.
What follows is a brief description of how this cohomology
is defined.  The reader should look at Chapter 29 of \cite{KM} for
the detailed story.

\subsubsection*{\(\Spin^{\bbC  }\) structures:}  Given an orientated 3-manifold
\(M\),  a \(\Spin^{\bbC }\) structure \(\grs\) on \(M\) is an equivalence class of a
pair consisting of a Riemannian metric on \(M\) and a lift, \(\pi
\co \scrF \to   M\), of the principle \(\SO (3)\)-bundle of oriented, orthonormal
frames for \(TM\) to a \(\U (2)\)-principle bundle.  Two such pairs are equivalent if there is a path of
metrics and a corresponding lift along the path that interpolates
between one pair and the other.  The set of equivalence classes is a
torsor modeled on \(H^{2}(M; \bbZ )\).   

Let \(\bbS\) denote the associated
bundle \(\scrF \times _{\U (2)} \bbC ^{2}\).   The standard Hermitian metric
on \(\bbC ^{2}\) gives \(\bbS \) a Hermitian fiber metric.  The bundle \(\bbS\)
is a Clifford module for \(T^*M\).  This is to say that there is an
endomorphism, \(\cl \co 
T^*M \to   \End\,(\bbS)\), such that
\[\cl(a_{1})\cl(a_{2}) = -\langle a_{1}, a_{2}\rangle  - \cl\, (*(a_{1} \wedge
a_{2})),\] 
where \(\langle  \cdot,\cdot \rangle \) denotes the metric inner product.
Endomorphisms in the image of \(\cl\) are anti-Hermitian.  The endomorphism
\(\cl\) induces two auxilliary homomorphisms. The
first, \(\hatc\), is defined as:
\[\hatc\co \bbS \otimes  T^*M  \to  \bbS, \quad \eta  \otimes
b\mapsto \cl(b)\,\eta .\]  
The second map is a quadratic, bundle-preserving map from \(\bbS  \) to
\(iT^*M\). Given \(\eta  \in \bbS \),  we will write in what follows its image under
this map as \(\eta ^{\dag }\tau\eta \); the latter is in turn defined
by the rule: 
\[ \langle b, \eta ^{\dag }\tau\eta \rangle  = \eta ^{\dag }\cl(b)\,\eta .\] 

Let \(\det \,(\bbS)\) denote the complex, Hermitian
line bundle \(\bigwedge^{2}\bbS \).  Its first Chern class,
\(c_1(\det\, \bbS)\in H^{2}(M; \bbZ )\), is called the \textit{canonical} cohomology class of the \(\Spin^{\bbC }\)
structure \(\grs\), and is often denoted as \(c_1(\grs)\). 
 A Hermitian connection, \(\bbA\), on \(\det\, (\bbS)\) and the
Levi-Civita connection on \(TM\) together define a metric compatible
covariant derivative, \(\nabla_{\bbA}\), on the space of sections of
\(\bbS\).  This covariant derivative is used to construct the
corresponding Dirac operator, 
\[
D^{\bbA}=\hatc (\nabla_{\bbA}):C^{\infty }(M; \bbS) \to
C^{\infty}(M; \bbS).
\] 
In what follows,
\(B_{\bbA}\) denotes the Hodge star of the curvature 2-form of the
connection \(\bbA\). This is an \(i\bbR\)-valued 1-form on \(M\).  
\medbreak
 
\subsubsection*{The Seiberg-Witten equations:}  Fix an orientation for
\(M\), a \(\Spin^{\bbC }\) structure.
Choose a metric on \(M\) and a lift, \(\pi   \co \scrF \to   M\), of its oriented orthonormal
frame bundle to a principle \(\U (2)\)-bundle, so that the pair is in
the equivalence class specified by the given \(\Spin^\C\) structure.  
Let \(\varpi\) be a closed 2-form on \(M\). 
We call a pair, \((\bbA, \Psi )\), consisting of a Hermitian connection on \(\det\,
(\bbS)\) and a section of \(\bbS\) a {\em configuration}. 
The group \(C^{\infty }(M; \U (1))\) acts on the space of configurations
in the following fashion:  Let \(u: M \to  \U (1)\).  Then \(u\) sends
a configuration, \((\bbA, \Psi )\), to \((\bbA - 2u^{-1}du, u\Psi )\).   Two solutions
obtained one from the other in this manner are said to be
\textit{gauge equivalent}.  The group \(C^{\infty }(M; \U (1))\) is called
the \textit{gauge group}.

In the most general form, the Seiberg-Witten equations ask that  a configuration
\((\bbA, \Psi )\) obey
\begin{equation}\begin{cases}
B_{\bbA } - \Psi^{\dag }\tau\Psi + i*\varpi  -\grT= 0  & \textit{  and    }\\
D^{\bbA}\Psi  -\grS= 0 , & \end{cases}
\label{(1.7)}
\end{equation}
where the pair \((\grT, \grS)\) is a small perturbation arising as
the formal gradient of a  gauge-invariant function of $(\bbA, \Psi
)$. It is in general needed
to guarantee the transversality properties necessary for the
definition of Seiberg-Witten-Floer cohomology. See Chapters 10 and 11
in \cite{KM}. 

Since the Seiberg-Witten equations are gauge invariant,
the gauge group acts on the space of Seiberg-Witten solutions as well.
Use \(\scrC\) in what follows to denote the set of gauge equivalence
classes of solutions to (\ref{(1.7)}).

 \subsubsection*{The cochain complex:}  Assume from this point on that:
\begin{equation}\label{irreducible}
\textrm{The gauge group acts freely on the space of solutions to (\ref{(1.7)}).}
\end{equation}
Then, for a suitably generic choice of \(\varpi \), \(\grT\), and \(\grS\),
the set \(\scrC\) of gauge equivalence classes of solutions to (\ref{(1.7)})
is finite.  (In fact \((\grT, \grS)\) can be set to be trivial for
this purpose). 

Paraphrasing Definition 29.1.1 in \cite{KM}, the case when 
\begin{equation}\label{SW-monotone}
2\pi c_1(\grs)-[\varpi]=\textsf{t} \, 2\pi c_1(\grs)
\end{equation}
when \(\varpi\neq 0\) for some real \(\textsf{t} > 0\) is said to be {\it positive monotone}, and it
is said to be {\em negative monotone} when \(\textsf{t}<0\). The condition
(\ref{irreducible}) holds in the monotone case. (Cf. e.g. Lemma 29.1.2
in \cite{KM}).

Assume \(\varpi\neq 0\) and monotonicity for the rest of this subsection.
In this case, the cochain complex used to
define the Seiberg-Witten-Floer cohomology is
\[
CM^*= \bbZ\scrC, \]
the free \(\bbZ  \)-module that is generated by \(\scrC\).

\subsubsection*{A relative grading:}  The complex
\(\bbZ \scrC\)  has a relative \(\bbZ /p\bbZ  \)-grading where \(p\) is the
divisibility in \(H^{2}(M; \bbZ )/\op{Tors}\) of the class \(c_{1}(\grs)\).  This
grading is defined using  the spectral flow of a 1-parameter family of unbounded,
self-adjoint, Fredholm operators \(\grL_\grc\) on \(L^{2}(M; iT^*M \oplus  \bbS \oplus
i\, \bbR )\), whose end members have trivial cokernel. The parameter \(\grc\)
of the family is from the space of configurations.
The precise definition of \(\grL_\grc\) will be given in
(\ref{(3.16)}); for now it
suffices to say that the operator \(\grL_\grc\) is obtained from the
linearization of (\ref{(1.7)}) at a given configuration \(\grc=(\bbA,
\Psi )\), and the operators associated to gauge equivalent
configurations are conjugate to each other. 
A configuration \(\grc\) 
is said to be {\em nondegenerate} when \(\grL_\grc\) has 
trivial cokernel. As just noted, this notion only depends on the gauge
equivalence class of \(\grc\).

Let \(\grc_-, \grc_+\) be two nondegenerate gauge equivalence classes
of configurations and use
\(\grP = \grP (\grc _{-}, \grc _{+})\) to denote the space of
piecewise differentiable maps from \(\bbR\) to the configuration space
\(\op{Conn}(\det \bbS)) \times C^{\infty}(M; \bbS )\) which have \(s \to
-\infty\) limit that is a configuration in the gauge equivalence class of
\(\grc _{-}\) and \(s \to \infty \) limit that is a configuration in the gauge equivalence class of \(\grc_{+}\).
The group \(C^\infty(M, \U \, (1))\) acts on \(\grP\) by the following
rule: an element \(u\in C^\infty(M, \U \, (1))\) sends a path
\((\bbA(s),\Psi (s))\) to the path \((\bbA(s)-2u^{-1} du, u\Psi (s))\). Two
elements in \(\grP\) in the same orbit under this action are said to
be gauge equivalent. 

The family of operators \(\grL_{\grc (s)}\) associated to every
element \(\grc(s)\) in \(\grP (\grc_-, \grc_+)\) has a 
well defined spectral flow.  
It turns out that this spectral flow depends only on \(\grc_-,
\grc_+\), and the path component of
the gauge equivalence class of \(\grc (s)\) in the orbit space \[
\scrB(\grc_-, \grc_+)=\grP (\grc_-, \grc_+)/C^\infty(M, \U \, (1)).
\]
With this understood, we call this class in \(\pi_0\big(\scrB(\grc_-, \grc_+)\big)\)
the {\em relative homotopy class} of \(\grc (s)\), and 
denote by \(\scrI (\grc_-, \grc_+; \grh)\in \bbZ\) the spectral
flow of elements in \(\grP (\grc_-, \grc_+)\) with relative
homotopy class \(\grh\). Meanwhile, note that
 the set of relative homotopy classes, \(\pi_0\big(\scrB(\grc_-,
 \grc_+)\big)\), is a torsor over the groups of components of the
 space of
 gauge transformations, i.e. \(H^1(M; \bbZ)\simeq H_2(M;
\bbZ)\). Moreover, given two relative homotopy classes \(\grh\),
\(\grh'\) with their difference \(\grh-\grh'\) viewed as an element in
\(H_2(M;\bbZ)\), their respective spectral flows differ by
\begin{equation}\label{SW-spec}
\scrI (\grc_-, \grc_+; \grh)-\scrI (\grc_-, \grc_+; \grh')=\langle c_1(\grs),\grh-\grh'\rangle.
\end{equation}

Finally, note that the genericity condition on \((\varpi  , \grT, \grS)\) is such that 
\(\grc=(\bbA,\Psi )\) is nondegenerate when it solves (\ref{(1.7)}).   
Assume in what follows that the
triple \((\varpi, \grT, \grS)\) is chosen so that this condition hold.
Applying the above discussion to the cases when \(\grc_-\),
\(\grc_+\in \scrC\), it follows from (\ref{SW-spec}) that  the mod-\(p\) reduction of the spectral flow, 
\[
\scrI (\grc_-, \grc_+):=\scrI (\grc_-, \grc_+; \grh) \mod p,
\] 
provides a relative \(\bbZ/p \bbZ\)-grading between any two
generators \(\grc_-, \grc_+\) of the Seiberg-Witten
Floer cochain complex.

\subsubsection*{The differential:}  The
differential for the Seiberg-Witten Floer cohomology is
defined using solutions to the Seiberg-Witten
equations on \(\bbR  \times  M\).  These equations are viewed
as a system of equations for a smooth map \(s\mapsto
(\bbA(s),\Psi(s))\) from \(\bbR\) to \(\op{Conn}(\det (\bbS))\times
C^\infty (M; \bbS)\).
The most general form of the equations read:
\begin{equation}\label{(1.8)}
 \begin{cases}
\frac{d}{ds}\bbA+B_{\bbA}- \Psi^{\dag}\tau\Psi+ i*\varpi  - \grT  = 0\\
\frac{d}{ds}\Psi+ D^{\bbA}\Psi  - \grS = 0,
\end{cases}
\end{equation}
where \(\varpi\) and \((\grT, \grS)\) are as in
(\ref{(1.7)}). 
Given \(\grc_-, \grc_+\in \scrC\), an \textit{instanton from \(\grc_-\) to
\(\grc_+\)} is an element \((\bbA(s),\Psi(s))\) in \(\grP (\grc_-,
\grc_+)\) solving the above equations. 
It is assigned a relative homotopy class and a spectral flow in the
manner described in
the previous paragraph. Moreover, a gauge transformation 
of an instanton is another instanton.

Let \(\EuScript{M}_1(\grc _{-}, \grc _{+})\) denote the following:
It is the empty set when  \(\scrI(\grc _{-}, \grc_{+})\neq 1\). When
\(\scrI(\grc _{-}, \grc_{+})= 1 \mod p\), let it be the space consisting
of gauge equivalence classes of instantons 
from \(\grc _{-}\) to 
\(\grc _{+}\) with spectral flow 1. By (\ref{SW-spec}), 
the relative homotopy classes of  such instantons lie in 
a torsor over \(\ker c_1(\grs)\).   As explained in Section
30.1 of \cite{KM}, this fact, together with the monotonicity
condition, ensures the compactness of \(\scrM_1 (\grc_-, \grc_+)\).
With suitably generic choice of \((\varpi, (\grT, \grS))\), the
space \(\scrM_1(\grc _{-},\grc_{+})\) will possess the
following desirable properties: 
It is a 1-dimensional smooth manifold diffeomorphic to a
disjoint union of finitely many copies of \(\bbR\),  
Moreover,  the
1-parameter group \(\bbR\) acts freely on \(\scrM_1(\grc _{-},
\grc_{+})\) by translating any instanton a constant amount with
respect to the parameter \(s\in \bbR \), and the aforementioned
diffeomorphism from \(\bbR\) to each component 
induces this \(\bbR\)-action. 
In addition, each 1-dimensional
component has an assigned weight, either \(+1\) or \(-1\).  Suffice it to say
for now that this weight is obtained by comparing two canonical
orientations.  The first orientation is that defined by the generator
of the \(\bbR \)-action; and the second is defined using the determinant
line of a certain family of Fredholm operators that is obtained from
the linearization of the expression on the left hand side of
(\ref{(1.8)}).   Use \(\sigma(\grc _{-}, \grc _{+})\) to denote the sum of
the \(\pm 1\) weights that are assigned to the 1-dimensional components of
\(\EuScript{M}_1(\grc _{-}, \grc _{+})\) when the set of such components is
non-empty.  Set \(\sigma (\grc _{-}, \grc _{+})\) to equal zero otherwise.

Given \(\grc  \in \scrC\), 
the action of the differential on \([\grc ]\) is given by the following
formula:
\begin{equation}\label{(1.9)}
\delta\grc _+= \sum_{\grc _-\in \scrC} \sigma (\grc _-, \grc _+) \, \grc _-.
\end{equation}
This differential has square zero and {\em increases} the relative \(\bbZ /p\bbZ
\)-grading by 1.

It turns out that in our setting, the perturbation terms \(\grT \) and \(\grS
\) may be set as zero; and the assumption (\ref{irreducible}) and 
and the various properties of \(\scrC\) and \(\scrM_1(\grc_-,
\grc_+)\) required for the definition of the Seiberg-Witten
cohomology, as sketched above,  may be obtained by a suitable generic
choice of \(\varpi\) alone. 

\subsubsection*{The Seiberg-Witten cohomology:} 
To summarize, in this subsection we defined a Floer cochain complex
with coefficient \(\bbZ\) associated
to each parameter set:
\[
\boldsymbol{\mu}_S=\{M, \grs, \varpi, g, \grq\},
\]
where \(M\) is an oriented closed 3-manifold; \(\grs\) and \(\varpi\) are
a \(\Spin^\bbC\) structure and a closed 2-form satisfying the monotonicity assumption (\ref{SW-monotone});
\(g\) is a Riemannian metric on \(M\); and \(\grq=(\grT, \grS)\) is as
in (\ref{(1.7)}). 
The Seiberg-Witten Floer \textit{co}homology is defined to be the homology of
the above Seiberg-Witten cochain complex, 
and it has a relative \(\bbZ /p\bbZ \)-grading. 
In the notation of Definition~30.1.1 in \cite{KM}, it is written as: 
\[
HM^*(M, \grs, -\pi [\varpi]): =H\, (CM^*, \delta).
\]
It is explained in 
\cite{KM} that the Seiberg-Witten Floer cohomology depends only on 
the following triple:
\BTitem\label{parameter-SW}
\item an oriented closed 3-manifold \(M\),
\item a \(\Spin^\bbC\)-structure \(\grs \) on \(M\), and
\item the nontrivial cohomology class \[
\mathpzc{p}_{\grs, \varpi}:=2\pi  ^2c_1(\grs)-\pi[\varpi]\in H^2(M;\bbR),
\]
modulo multiplication by positive numbers. 
\ETitem
Thus, the positive and negative monotone cases correspond to two 
different versions of Seiberg-Witten Floer cohomologies. In
particular, the positive monotone version agrees with the ``ordinary
Seiberg-Witten cohomology'':
\[
HM^* (M, \grs)=HM^* (M, \grs, 0).
\]

\subsubsection*{\it Remarks.} (1) The condition (\ref{irreducible}) implies that
the \(\widehat{HM}\) and \(\widecheck{HM}\) versions  of Seiberg-Witten
cohomology defined in \cite{KM} are the same, and the completed
version \(HM^\bullet\) is no different from the plain
\(HM^*\). Hence we denote all these by the same
notation \(HM^*\).

(2) In fact, to describe the dependence of the periodic Floer homology
on \(\ff\), it is more appropriate to regard \(\ff\) as an element in
\(\widetilde{\op{Symp}}(F, w_F)\), the universal covering of the
symplectomorphism group of \((F, w_F)\). This will be clarified in 
Appendix A below.

 \subsection{The isomorphism theorem}\label{sec:1(c)}

We now describe a way to associate a Seiberg-Witten parameter set 
\(\boldsymbol{\mu}_S\) to each Periodic Floer homology parameter set 
\(\boldsymbol{\mu}_M\).

(1) {\it From \(\ff\co (F, w_F)\circlearrowleft\) to \(M\).} 
Recall from the beginning of this section that a closed surface \(F\)
with volume form \(w_F\) together with a volume-preserving map
\(\ff\) on it defines an associated mapping torus \(M\) with
orientation induced from \(w_F\). They also define a closed 2-form
\(w_\ff\) on \(M\). 

(2) {\it From \(\ff\co (F, w_F)\circlearrowleft\) and \(J\) to \(g\).} There is a standard way to
define a metric \(g\) on \(M\) from \(w_F\) and a \(J\in
\mathcal{J}_\ff\). See Section \ref{sec: 3(c)} below.

(3) {\it From \(\Gamma\) to \(\grs\).} There is an isomorphism between the space of \(\Spin^\bbC\)
structures and the homology group \(H_1(M;\bbZ)\) as torsors over 
\(H_1(M;\bbZ)\), as follows: Let \(\grs\) be a \(\Spin^\bbC\) structure 
on \(M\) defined by the equivalence class of the
Riemannian metric \(g\) and a lifting \(\scrF \to  M\) of the
frame bundle ,as described in Section
\ref{sec:1(b)}. Let \(\bbS \) denote the associated
\(\U\, (2)\)-bundle of spinors.  Clifford multiplication by the
1-form \(dt\) gives an anti-Hermitian automorphism of \(\bbS\)
whose square is a negative multiple of the identity.
Write the corresponding orthogonal decomposition of \(\bbS\)
into eigenbundles as \(E \oplus  E\otimes K^{-1}\) where \(E \to  M\) is a
complex, Hermitian line bundle. The convention takes the
left most summand to be that on which Clifford
multiplication by \(dt\) acts as a positive multiple of \(i\).
The assignment \(\grs\mapsto c_{1}(\det\, (E))\) of a cohomology
class to a \(\Spin^{\bbC }\) structure defines a 1-1
correspondence between the set of \(\Spin^{\bbC }\) structures
with \(H^{2}(M; \bbZ )\), and equivalently, with \(H_1(M; \bbZ)\) via
Poincar\'e duality.  Denote by \(\grs_\Gamma\) the \(\Spin^\bbC\)
structure corresponding to the homology class \(\Gamma\in H_1(M;\bbZ)\) 
under this 1-1 correspondence, and by
\(\bbS_\Gamma\) the associated spinor bundle. Note that under this
correspondence, the cohomology class \(c_\Gamma\) in Section
\ref{sec:1(a)}, which is relevant to the periodic Floer homology,
equals the class \(c_1(\bbS_\Gamma)\) in Section \ref{sec:1(b)}, which
in turn is relevant to the Seiberg-Witten-Floer cohomology. In particular, 
this means that the two versions of \(p\), defined differently in the
afore-mentioned two subsections, are actually identical, thus justifying
our use of the same notion for both.

(4) {\it From \(\Gamma\), \(w_\ff\) to \(\varpi\).} Let
\[\varpi_r=2rw_\ff+\wp,\] 
where \(r>0\) and \(\wp\) is a closed
2-form in the cohomology class \(2\pi c_1(\grs_\Gamma)\).
With this choice, 
\[
\mathpzc{p}_{\grs_\Gamma, \varpi_r}=-2r\pi[w_\ff].
\]
Note that with these choices, \(\Gamma\) is positive or negative
monotone with respect to \([w_\ff]\) precisely when 
\(\grs_\Gamma\) and \([\varpi_r]\) are respectively positive or negative monotone
with respect to each other.

(5) {\it The nonlocal perturbation}. Let \(\grq\) be arbitrary.

Because the Seiberg-Witten Floer cohomology only depends on
(\ref{parameter-SW}), the Seiberg-Witten  Floer cohomology associated
to the above choices does not depend on the extraneous data \(r\),
\(\grq\), and the choice of \(\wp\) within its cohomology class.

We are now ready to state the isomorphism theorem under the
monotonicity assumption.

\begin{thm}\label{thm:1}
Let \(\boldsymbol{\mu}_P=\{(F,
w_F), \ff, \Gamma, J\}\) be a periodic Floer parameter set described
in (\ref{PFH-par}), where \(\Gamma\in H_1(M; \bbZ)\) is monotone with respect to
\([w_\ff]\); and associate to it a Seiberg-Witten parameter set
\(\boldsymbol{\mu}_S=\{M, \grs_\Gamma, \varpi_r, g, \grq\}\) according
to the recipe described above. Then 
there is an isomorphism between the two Floer (co)homologies
with coefficient \(\bbZ\):
\[
HP_*(\ff\co (F, w_F)\circlearrowleft, \Gamma )_J\simeq
HM^{-*}(M, \grs  _\Gamma, -\pi[\varpi_r]),
\]
which reverses the relative \(\bbZ/p\bbZ\)-gradings. 
In particular, if the degree \(d_{\Gamma }
<  0\), then both Floer (co)homologies vanish. 
If \(d_{\Gamma } = 0\), then both sides are isomorphic to \(\bbZ \). 
\end{thm}
Since \([w_\ff]\) is invariant under Hamiltonian isotopies of \(\ff\),
and since the Seiberg-Witten Floer cohomology on the right hand side of
the isomorphism is invariant under rescaling of \([w_\ff]\) by
positive numbers, an immediate consequence of the above isomorphism is:
\begin{cor}\label{cor:PFH-inv}
The periodic Floer homology only depends on the orientation of \(F\),
the Hamiltonian isotopy class of \(\ff\), and the monotone homology class \(\Gamma\).
\end{cor}

By working with local coefficients, both the periodic Floer homology
and the Seiberg-Witten Floer cohomology have a more general
definition, where \(\Gamma\) need not be monotone. The above isomorphism
theorem and corollary have corresponding generalizations, see Theorem
\ref{thm:2} and Corollary \ref{cor:PFH-inv2} below.

{\it Remarks.}
(1) The assertions of Theorem \ref{thm:1} when
\(d_{\Gamma } \leq  0\) follow from Theorem 1.3 in the article
\textit{SW\(\Rightarrow\)Gr }from \cite{T5}.  See also the main
theorem in \cite{KT}.  

(2) 
The equality between the 
the Euler characteristics of the two Floer (co)homologies in Theorems
\ref{thm:1} and \ref{thm:2} 
follows from Theorem 1 in the article \textit{SW\(=\)Gr }from \cite{T5}; see also \cite{MT} and \cite{HL}.  

(3) This isomorphism theorem is related to a conjecture of D. Salamon,
see Conjecture 10.1 in \cite{S}. 
Roughly speaking, this conjecture says that the
Seiberg-Witten cohomology of the mapping torus of \(\ff\) is
equivalent to a certain version of the Floer homology of 
symplectic fixed points for certain self map of \(d\)-fold symmetric
product of \(F\) that is defined from \(\ff\).
By viewing the elements \(\{(C_i, m_i)\}_i\) in
\(\mathcal{M}_1 (\Theta_+, \Theta_-)\) as holomorphic (branched)
multi-sections of the surface bundle \(\bbR\times M\to \bbR\times
S^1\), the periodic Floer homology corresponds intuitively to this
latter Floer homology. 
\par

As we noted before, when \(\Gamma\) is positive monotone or
negative monotone, the Seiberg-Witten cohomology on the
right hand side of the isomorphism corresponds respectively to 
the positive monotone or negative monotone version described in
Section 30.1 of \cite{KM}; in particular, in the former case, it
coincides with the ``ordinary Seiberg-Witten cohomology''
\(HM^{-*}(M, \grs_\Gamma)\).
 On the other hand, Taubes proved a sister
version of the above theorem in \cite{T1}--\cite{T4}, which states
that 
\begin{thm*}
Given a closed contact 3-manifold \(M\), and a homology class
\(\Gamma\in H_1(M;\bbZ)\), there is an isomorphism from the 
 associated embedded contact homology \(ECH_*\) 
to the Seiberg-Witten 
Floer cohomology \(\widehat{HM}^{-*}(M, \grs_\Gamma)\) which reverses
the relative gradings.
\end{thm*}
In the above, \(\grs_\Gamma\) is defined in the same way as in Item (3) in the
beginning of this subsection, with \(c_1(K^{-1})\) now denoting the
Euler class of the 2-plane field which defines the contact structure.

Recall from the last remark in Section \ref{sec:1(b)} that
\[\widehat{HM}^{-*}(M, \grs)=HM^{-*}(M, \grs)\]
when \(c_1(\grs)\) is nontorsion; it follows immediately that
\begin{cor}
Let \(M\) be the mapping torus of an orientation-preserving map \(\ff\co
(F, w_F)\circlearrowleft\), and \(\Gamma\in H_1(M; \bbZ)\) be a positive monotone
class with respect to \([w_\ff]\), then 
\[
ECH_*(M, \xi ; \Gamma)\simeq HP_* (\ff, \Gamma+h_\xi),
\]
where the left hand side denotes the embedded contact
homology of \((M, \xi)\); \(\xi\) is a contact structure consistent
with ithe orientation of \(M\), and \(h_\xi\) is Poincare dual to the primary  
obstruction to finding a homotopy between the 2-plane fields \(\xi\) and  
\(\ker \pi_*\). 
\end{cor}

\par

The isomorphism in Theorem \ref{thm:1} may be used to compute one
Floer (co)homology from the other, depending on which is easier to
compute in the context. Here are some sample applications of this
type.

First, note that monotone classes with
respect to a fixed \([w_\ff]\) are specified by their degree
\(d_\Gamma\). In particular, \(d_\Gamma<g-1\) when \(\Gamma\) is
positive monotone, and \(d_\Gamma>g-1\) when \(\Gamma\) is negative monotone.
By the definition of periodic Floer homology, \(d_\Gamma\neq g-1\).

{\it I. From \(HP\) to \(HM\):} 
The degree-1 periodic Floer
homology of the mapping torus of \(\ff\) is 
identical with a certain version of the symplectic
Floer homology of \(\ff\). (See Appendix B for a precise statement).
The latter has been computed for all area-preserving surface
automorphisms, see \cite{C} and references therein.
This implies that the corresponding  
Seiberg-Witten-Floer cohomology of degree \(<2\) is known for all
mapping tori of 
area-preserving surface automorphisms as well. In \cite{HS}, the 
periodic Floer homology of certain Dehn twists has been computed. 
This also provides some new computations of the Seiberg-Witten-Floer cohomology.  

{\it II. From \(HM\) to \(HP\):}
 Corollary 31.5.2 in \cite{KM} gives a long exact sequence relating
the positive monotone version of \(HM_*\), the negative monotone
version of \(HM_*\), and the bar-version of the Seiberg-Witten
homology at the ``balanced perturbation'' \(c_b=-2\pi ^2c_1(\grs)\)
(Cf. Section 31.1 of \cite{KM}). The dual version of this exact
sequence takes the following form:
\begin{equation}\label{ex-sq}\begin{split}
\cdots \stackrel{j^*}{\longrightarrow}HM^* (M, \grs, c_+)
\stackrel{i^*}{\longrightarrow}\overline{HM}^* (M, \grs, c_b)
& \stackrel{p^*}{\longrightarrow}HM^* (M, \grs, c_-)\\
& \quad \stackrel{j^*}{\longrightarrow}HM^* (M, \grs,
c_+)\stackrel{i^*}{\longrightarrow}\cdots , 
\end{split}
\end{equation}
where the first term and the third term
above denote the positive monotone version and 
the negative monotone version of the Seiberg-Witten 
Floer cohomology respectively. Thus, knowledge of the ordinary
Seiberg-Witten cohomology enables one to compute the periodic Floer
homology in degrees larger than \(g-1\).

\begin{cor}[\bf Cotton-Clay]\label{cor:adjunction}
Suppose the genus \(g\) of the surface \(F\) is positive.
Let \(\Gamma\) be a monotone class with respect to \([w_\ff]\) with
degree \(d_\Gamma>2g-2\). Then 
\[
HP_* (\ff, \Gamma)\simeq 
\overline{HM}^{-*}(M, \grs_\Gamma, c_b),
\]
where \(\overline{HM}^*(M, \grs, c_b)\) denotes the
 bar-version of the Seiberg-Witten Floer cohomology at the ``balanced
perturbation'' .
\end{cor}

\pf Let \(\Gamma\) be a monotone class with degree \(d_\Gamma>2g-2\geq
g-1\) as in the statement of the Corollary, and take
\(\grs=\grs_\Gamma\) in the above long exact sequence. By the
equivalence \(HM^* (M, \grs, c_+)\simeq HM^* (M, \grs)\)
(see Theorem 31.1.2 in \cite{KM}) and
the adjunction inequality (see Corollary 40.1.2 in \cite{KM}, or
rather, its dual version), this condition
on degree implies that \(HM^* (M, \grs, c_+)=0\). On the other
hand, by Theorem~\ref{thm:1}, \(HM^* (M, \grs, c_-)\simeq
HP_{-*}(\ff, \Gamma)\), since \(\Gamma\) is negative monotone. 
The claim of the Corollary follows by 
combining these observations with the long exact sequence above. 

{\it Remark.}
The bar-version of Seiberg-Witten Floer cohomology 
is defined solely from the space of flat connections on \(M\) and the
family of Dirac operators parametrized by it. In Chapters 33-35 of
\cite{KM}, Kronheimer and Mrowka showed that the bar-version of
Seiberg-Witten Floer homology can be computed from a
``Coupled Morse chain complex'' on the torus of flat
connections. Given a Morse function on the torus \(H^1(M; \bbR)/H^1(M;
\bbZ)\), this
chain complex is constructed from its moduli spaces of gradient flow
lines of dimensions 1 and 3, together with appropriately chosen
\(\widecheck{\text C}\)ech cocycles in cohomology classes determined by  \(c_1(\grs)\) and the triple cup product
 \(\Upsilon\co \bigwedge^3H^1(M)\to H^3(M)\simeq \bbZ\).
Explicit formulas are given in Chapter 35 of \cite{KM} for \(\overline{HM}_*\) in the
special cases when when \(b_1(M)\leq 3\), or when \(\ff\) is isotopic
to the identity and \(c_1(\grs)\) is proportional to the Poincar\'e
dual of \([F]\). A modification of the computation for the last
example in \cite{KM} finds that under the assumptions that \(|\langle
c_1(\grs), [F]\rangle|>0\) and \(g>0\), 
the homology groups of these complexes are always non-torsion over \(\bbZ\).
Note in contrast that when \(b_1(M)>1\), 
the Euler characteristics of  both Seiberg-Witten Floer cohomology and
periodic Floer homology are trivial
by the adjunction inequality for Seiberg-Witten invariants of
3-manifolds, as the latter is independent of perturbations.

The previous Corollary may be generalized as follows:
By combining the isomorphism theorem \ref{thm:1} with 
the involution in Seiberg-Witten theory, described in the
Floer-homology context in Proposition 25.5.5 of \cite{KM}, one has the
following: 

If \(\Gamma\) is a monotone class of degree \(d>g-1\) with respect to a certain
\([w_\ff]\), then \(HP_*(\ff, \Gamma)\simeq HM^{-*}(M,
\grs_{\Gamma^*}, c_-)\), where \(\Gamma^*\) is the monotone class of
degree \(2g-2-d\) with respect to the  same class \([w_\ff]\).
Thus, it can be computed from \(\overline{HM}^*(M,
\grs_{\Gamma^*}, c_b)\) and \(HM^*(M, \grs_{\Gamma^*}, c_+)\) by
the long exact sequence (\ref{ex-sq}). This means that the periodic
Floer homology is completely determined by the ordinary and
bar-versions of the Seiberg-Witten-Floer cohomologies for monotone
classes of degrees \(0<d<g-1\), and the \(i^*\)-maps between the two
versions of \(HM^*\). In particular, it implies that the periodic
Floer homology with coefficient \(\bbZ\) for all surfaces of \(g\leq
2\), when it is well-defined, is completely determined by
the cohomology of \(M\), except the case when \(g=2\) and
\(e_\Gamma=c_1(K)\). In this last case it is determined by the 
\(i^*\)-map from the relevant ordinary Seiberg-Witten-Floer cohomology, which is \(\bbZ\) in this case.

\subsubsection*{Proving the isomorphism theorems.}
The isomorphism stated in Theorems \ref{thm:1} and \ref{thm:2}
is obtained by using a
version of Seiberg-Witten equations (\ref{(1.7)}), parametrized by
\(r\in \bbR^+\), to define the cochain complex for the
Seiberg-Witten Floer cohomology.  
For suitably chosen \(\ff\) and \(J\) (Cf. Section \ref{sec:2}) and
any sufficiently large \(r\), we
construct an {\it isomorphism} between the periodic Floer homology
complex and the corresponding Seiberg-Witten-Floer cohomology
complex, which reverses the relative \(\bbZ/p\bbZ\)-grading.
Note in constrast that there are no such isomorphisms between the
\(ECH\) and \(HM\) complexes in the proof of
\(ECH=HM\) in \cite{T1}--\cite{T4}. Instead, the latter proof relies on a
filtration on the Floer (co)chain complexes, and
constructed injective chain maps from finitely-generated
sub-complexes of the \(ECH\) chain complex to the \(HM\) cochain complex. An
isomorphism between the respective Floer (co)-{\it homologies} is
obtained after taking a direct limit. This filtration argument is not
needed in this article; this is because the appearance of the non-exact perturbation
two-form \(w_\ff\)  guarantees that the
Seiberg-Witten-Floer cochan complex in our context is finitely generated. 

The proofs of Theorems \ref{thm:1} and \ref{thm:2} have the same 
analytic details, which are in most respects very similar to
those in the proof of Theorem 1 in \cite{T1}, which is in turn based
on the techniques from \cite{T5}.  In particular, most of
the proof of Theorem 1 can be borrowed with only notational
or cosmetic changes from corresponding arguments in
\cite{T1}-\cite{T4}.  This said, the authors have chosen to direct
the reader to such arguments when appropriate with faith
that the reader can make the required changes.  The
detailed arguments in what follows give those parts of the
proof of Theorems \ref{thm:1} and \ref{thm:2} that differ 
substantially from their
\cite{T1}-\cite{T4} analog.  Note in this regard that to show that the
map from \(ECH\) to \(HM\) is onto, the arguments
in \cite{T1} and \cite{T4} use at key points an estimate for the
spectral flow of a certain 1-parameter family of Dirac
operators.  This is the spectral flow estimate given by
Proposition 5.1 of \cite{T7}.  This spectral flow estimate
plays a profound role in \cite{T1} and \cite{T4} and is the key new
idea which was not present in \cite{T5}.  As it turns out,
there is no need here for this spectral flow estimate; judicious use
of the action functional in Seiberg-Witten-Floer theory and 
cohomological considerations of the sort used in the
article \textit{SW\(\Rightarrow\)Gr }from \cite{T5} suffice.  (For example,
none of Section 5 in \cite{T4} is needed here).  Similar arguments are
used in \cite{L} to prove an analog of \textit{SW\(\Rightarrow\)Gr }
in the more general setting of 4-manifold with cylindrical ends
equipped with an asymptotically constant harmonic 2-form, which may not
be nowhere vanishing, and where higher-dimensional moduli spaces are
also considered. 

The proof of Theorem \ref{thm:1} is given in Section \ref{sec:4(d)},
and Section \ref{sec:6(c)} contains the modification needed to
generalize to Theorem \ref{thm:2}.

\medbreak 

\subsection{Notation and conventions}\label{sec:1(d)}

What follows uses \(c_{0}\) to denote a constant that is
greater than 1.  Unless a specific note to the contrary is
present, this constant is independent of other relevant
parameters in the discussion.  The precise value of \(c_{0}\)
can be assumed to increase in subsequent appearances.  The letter
\(\kappa\) is used for a similar purpose. 

The constructions in the subsequent sections use a fixed,
non-increasing function, \(\chi \co [0,\infty ) \to  [0, 1]\) which takes value 1 on
\(\big[0, \frac{5}{16}\big]\) and is zero on \(\big[\frac{7}{16}, \infty \big)\).

On the Seiberg-Witten-Floer theory side, we follow the gauge-theory
literature in which the boundary maps are defined by flows from
the \(s\to -\infty \) side to the \(s\to \infty\) side (and therefore
the {\it co}boundary maps go in the opposite direction). 
 
On the periodic Floer homology side, we follow the symplectic Floer
theory literature in which the boundary maps are defined by flows from
the \(s\to \infty \) side to the \(s\to -\infty\) side. 

This accounts for our notational convention that in the various
corresponding notations in the two Floer theories, 
the positive and the negative ends appear in the opposite order,
such as \(\mathcal{M}_1(\Theta_+, \Theta_-)\) v.s. \(\scrM_1(\grc_-,
\grc_+)\). It is for the same reason that the periodic Floer {\it
  homology} is isomorphic to the Seiberg-Witten-Floer {\it cohomology}.
  
\subsection{Acknowledgements}\label{sec:1(e)}  

The authors owe a debt of gratitude to Michael Hutchings
for sharing his knowledge of periodic Floer homology and
for his thoughts about the correspondence between the
latter and the Seiberg-Witten Floer homology. We also thank Andrew
Cotton-Clay for pointing out Corollary \ref{cor:adjunction} and other
helpful comments. 

Both authors were supported in part by NSF grants.

\section{Deforming the map \(\ff\) and the
  almost complex structure \(J\)}\label{sec:2}
\setcounter{equation}{0}

 The proof of Theorem \ref{thm:1} requires two preliminary steps.  These steps
 are needed in order to use the techniques from \cite{T1}-\cite{T4}.  The first
 step makes a small modification to the diffeomorphism, \(\ff\), that
 defines \(M\) and the associated periodic Floer homology chain complex.
 The second step makes a small modification to the almost complex
 structure on \(\bbR \times M\) that is used to construct the differential
 on this same chain complex.  With a class \(\Gamma \in H_{1}(M; \bbZ )\)
 given, the deformations that are described below can be made so the
 original and the new periodic Floer homology chain complexes have the
 same set of generators and identical differentials.  The deformation
 depends on the choice of a pair \((\delta, N)\), where \(\delta\) is
 a small, positive number, and \(N\) is a positive integer.  

Let \(\Lambda _{N}\) denote the set of periodic orbits of \(\ff\) with period \(N\)
 or less.  If \(\Gamma \in H_{1}(M; \bbZ )\) is any given class, then
 there exists \(N\) such that \(\Gamma \)'s version of the set \(\mathcal{A}\) is defined
 solely from periodic orbits in \(\Lambda _{N}\).  
 \medbreak

\subsection{The deformation of \(\ff\)}\label{sec:2(a)}

A small deformation of \(\ff\) is described here,  It changes \(\ff\) only near
the points that comprise the periodic orbits
in \(\Lambda _{N}\) but fixes these same orbits.    

Fix \(\gamma  \in  \Lambda _{N}\) and let
\(q \in  \{1, \cdots , N\}\) denote its period. A pair of functions \(z = (z_{1},
z_{2})\) defined on a neighborhood of \(\pt  \in   F\) are said to define a
\textit{symplectic} \textit{coordinate chart} centered at \(\pt\) when they
vanish at \(\pt\) and are such
that \(w_{F}\) is given by \(dz_{1}
\wedge  dz_{2}\) near \(\pt\).   Given such a pair of function \(z = (z_{1},
z_{2})\), then \(f^{q}\) appears as a map that can be written as
\begin{equation}\label{(2.1)}
z \mapsto K_{\gamma }\cdot z + \grq \,(z), \quad 
\text{where \(|\grq |\leq  c_{0}|z|^{2}\) and where \(K_{\gamma } \in  \op{Sl}\,(2; \bbR )\).}
\end{equation}
If \(\gamma  \) is elliptic, then the matrix \(K_{\gamma }\) is conjugate in
\(\op{Sl}\,(2; \bbR )\) to a rotation by an angle
\(2\pi \textsc{r}\), where \(\rrr  \in  (0, 1)\).
This \(\rrr\)  is called the
\textit{rotation number} of the periodic orbit.  It does not depend on
the chosen point \(\textup{p} \in   \gamma  \) or on the chosen symplectic coordinate
chart.  Moreover, if \(\gamma\) is nondegenerate, \(\rrr\) is
irrational. If \(\gamma  \) is hyperbolic, then \(K_{\gamma }\) is conjugate in
\(\op{Sl}\,(2; \bbR )\) to a matrix of the form 
\begin{equation}\label{eq:lambda}
\left[\begin{array}{cc}
\lambda  &0\\
0 &\lambda^{-1}\end{array}\right],
\end{equation}
where \(\lambda\) is a real number with \(|\lambda|>1\).  In all cases, the
eigenvalues of \(K_{\gamma }\) do not depend on the choice of the
symplectic coordinate chart.

\begin{defn}\label{(2.2)}
Given \(K_{\gamma }\), define functions \(\nu_\gamma\) and
\(\mu_\gamma \) on \([0, 2\pi ]\) from the type and rotation number of \(\gamma\) as follows:
\begin{itemize}
\item  If \(\gamma \) is elliptic with rotation number \(\rrr\),
then \(\nu_\gamma  = \frac{1}{2}\rrr \) and \(\mu_\gamma  = 0\).
\item  If \(\gamma \) is hyperbolic, then
\(\nu_\gamma =\frac{1}{2}\rrr\) and \(\mu_\gamma = -\frac{i}{4\pi}\ln
|\lambda  |\, e^{2i\rrr \tau}\) with \[
\rrr  = \begin{cases}
0&\text{when
\(\gamma \) is positive hyperbolic;} \\
\frac{1}{2} &\text{ when \(\gamma \) is negative hyperbolic,  }
\end{cases} \]
where \(\lambda\) is as in (\ref{eq:lambda}).
\end{itemize}
\end{defn}

Fix \(z \in  \bbC \), and solve the equation
\begin{equation}\label{(2.3)}
 \frac{i}{2}\partial_\tau\eta+ \nu _{\gamma}\eta+ \mu_\gamma \bar{\eta}= 0
\end{equation}
for a function \(\eta : \bbR  \to  \bbC  \) with \(\eta (0) = z\).  With \(\eta\)
and \(z\) viewed as vectors in \(\bbR ^{2}\), the solution has the form \(\eta
(\tau ) =U(\tau )\cdot z\), where \(U\) 
is a map from \(\bbR\) to \(\op{Sl}\,(2; \bbR )\) with the following properties:
\begin{itemize}
\item \(U(0)\) is the identity. 
\item If \(\gamma  \) is elliptic with rotation number \(\rrr\), then \(U(2\pi )\)
is the rotation through angle \(2\pi  \rrr\).  
\item If \(\gamma  \) is hyperbolic, then \(U(2\pi )\) is the matrix
  in (\ref{eq:lambda}).
\end{itemize}
For \(k \in  \bbZ \), set 
\[
U_{\gamma }^{k/q}=U(2\pi k/q).
\] 
Fix a point \(\pt \in  \gamma \).  As the almost complex structure \(J\) on \(\bbR
\times  M\) must map \(\ker(\pi _{* })\) to itself along \(\bbR
\times  \gamma\); thus it defines an almost complex structure on \(TF|_{\pt}\).  The latter is
denoted by \(\grj _{p}\).
Since \(\grj _{p}^{2} = -1\), it follows that there exist symplectic
coordinate chart \(z = (z_{1},
z_{2})\) for \(F\) near \(\pt\), such that \(\grj_{p}\frac{\p}{\p
  z_1}=\frac{\p}{\p z_2}\). 
Coordinates of this sort are said to be \textit{\(J\)-adapted}.

What follows describes the desired modification of \(\ff\).  

\begin{lemma}\label{lem: 2.1} Fix a
  positive integer \(N\). Then there exists a constant \(\kappa \geq
  1\) such that given
  \(\delta < \kappa ^{-1}\), there exists a diffeomorphism \(\ff'\co F \to
F\) with the following properties:
\begin{itemize}
\item[{\rm (1)}]  The map \(\ff'\) is area-preserving,
   which is to say that \((\ff')^*w_{F} = w_{F}\).

\item[{\rm (2)}]  All periodic orbits of \(\ff'\) are
   non-degenerate.

\item[{\rm (3)}] The maps \(\ff\) and \(\ff'\) have the same set of periodic orbits
   with period \(N\) or less.

   Moreover, corresponding orbits are either both elliptic or both
   hyperbolic.  If both elliptic, then they have the same rotation
   number.  If both hyperbolic, then they are either both positive
   hyperbolic or negative hyperbolic.  In either case, the eigenvalues
   of the corresponding versions of the matrix \(K_{(\cdot )}\) agree.

\item[{\rm (4)}]  Let \(\gamma  \) denote a periodic
orbit of period \(q \leq  N\).  Then there are
\(J\)-adapted, symplectic coordinates about the points that comprise
\(\gamma  \) such that for each \(k \in  \{0, \cdots , q-1\}\) the
map \(\ff'\) near \(\ff^{k}(\pt)\) in these coordinates is given by
\[
z \mapsto U_{\gamma }^{(k+1)/q}\cdot (U_{\gamma }^{k/q})^{-1}z
+ \grp _{k}(z) ,
\]
 where \(\grp _{k}\) is zero in a neighborhood of \(z =0\). 

\item[{\rm (5)}]  The coordinate disks in \(F\) of radius \(\delta\)
  about the points in the orbits from \(\Lambda_N\) are \(J\)-adapted,
  pairwise disjoint; and the maps \(\ff\) and \(\ff\) agree on the complement of their union.
\end{itemize}
\end{lemma}

This lemma is proved momentarily.  
\bigbreak
 
\subsection{The modification of J}\label{sec:2(b)}

A deformation of \(\ff\) to a map \(\ff'\) as described in Lemma~\ref{lem: 2.1} is exploited in later
sections with the help of a corresponding modification of the almost
complex structure \(J\). The definition of this modified almost complex structure requires a preliminary modification of the vector field \(\partial_t\), so as to obtain a vector field whose return map gives \(\ff\)́near the points in \(\Lambda_N\). The latter vector field is denoted by \(\partial_t\), and it is constructed in Part 1 of this subsection. Part 4 describes an allowed set of almost complex structures that can serve as the desired modification of \(J\). The intervening Parts 2 and 3 of this subsection jointly define one of the criteria that are used in Part 4 to describe the allowed set of modified almost complex structures.

\subsubsection*{\it Part 1. The embedding \(\varphi_\gamma\) and the
  modified vector field \(\partial_t'\).}  
It is assumed in what follows that \((\delta, N)\) has been chosen
subject to the conditions in Lemma~\ref{lem: 2.1}.  Let \(\ff'\) denote a
map given by (\ref{(2.1)}) using this data.  
This modification of \(\ff\) is used in what follows to construct a
corresponding modification of the vector field \(\partial_t\).  The new
version of \(\partial_t\) is described below. 

It proves useful to view elements in \(\Lambda _{N}\) now as integral
curves of \(\partial_t\). Fix \(\gamma \in \Lambda _{N}\), now viewed as a
closed integral curve of the vector field \(\partial_t\) on \(M\).  As such,
\(\gamma\) has a tubular neighborhood embedding \(\varphi_\gamma : S^1
\times D \to M\), where \(D \subset \bbC \) is a disk of small radius.  Such
an embedding is described momentarily with certain desirable
properties.

Use \(\tau\) to denote the affine parameter on \(S^1\), and use \(z\) to denote
the complex coordinate for \(\bbC\).  The desired embedding,
\(\varphi_\gamma \), has the following properties:

\begin{lemma}\label{lem: 2.2}
  Given \(N \geq 1\), there exists \(\kappa \geq 1\) with the following significance: 

Fix \(\delta \in (0, \kappa ^{-1})\) and invoke
  Lemma~\ref{lem: 2.1} to construct a map \(\ff'\).
  Fix such an \(\ff'\), and let \(D \subset \bbC \) denote the disk of radius
  \(100\, \delta \).  Then there is a set \(\bigl\{\varphi_\gamma : S^1 \times
  D \to M\bigr\}_{\gamma\in\Lambda_N}\) of embeddings with pairwise
  disjoint images, and such that for each \(\gamma \in \Lambda _{N}\),
  the embedding \(\varphi_\gamma\) obeys:
\begin{itemize}
\item  The composition of \(\varphi_\gamma \)
with the projection \(\pi: M\to S^1\) sends \((\tau , z)\)
to \(q\tau \) with
\(q\) here denoting the period of \(\gamma\).

\item The vector field \(q^{-1}\partial_t\) is the push-forward from
  \(S^1 \times D\) of a vector field that has the form
\[
\partial _{\tau} - 2i\, (\nu z+\mu \bar z+\grr)\, \partial_z+ 2i\, (\nu\bar z+\bar
\mu z+\bar \grr)\, \partial_{\bar z},
\]
where \(\nu\) is an \(\bbR\)-valued function on \(S^1\) and \(\mu\) is
a \(\bbC\)-valued function.  Meanwhile \(\grr \) is such that its
absolute value is bounded by \(\kappa|z|^{2}\) and its first derivatives are
bounded in absolute value by \(\kappa|z|\).

\item  Use
\(\varphi_\gamma \) to
identify the disks
\(\bigcup _{0\leq k\leq q} (\{2\pi k/q\} \times
D)\) with their embedded images via \(\varphi_\gamma \)
in \(F = \pi ^{-1}(0)\).  Granted this identification,
the map \(\ff'\) is obtained by traveling for time \(2\pi /q\)
along the integral curves of the vector field
\[
\op{v}_\gamma:=\partial _{\tau } - 2i \, (\nu_\gamma z +\mu_\gamma \bar z+\grr
')\, \partial _{z} + 2i\, (\nu_\gamma \bar z+\bar \mu _\gamma z+\bar\grr
')\, \partial_{\bar z}
\]
where \(\grr '\) has the following properties: It vanishes where
\(|z|\) is in a neighborhood of \(0\), and it is equal to \(\grr\) where
\(|z| > \delta \).

\item The vector field \(\varphi_{\gamma*}\partial _{z}\) along \(\bbR
\times \gamma \) is such that \(J(\varphi_{\gamma*}\partial _{z})
=\varphi_{\gamma*}(i\partial _{z})\).
\end{itemize}
\end{lemma}

This lemma is proved in Section~\ref{sec:2(d)}.

Define a vector field, \(\partial _{t'}\), on \(M\) as follows:  
\[\partial _{t'}=\begin{cases}
(\varphi_\gamma  )_* \op{v}_\gamma & \text{over the image of
  \(\varphi_\gamma\), for all \(\gamma\in \Lambda_N\),}\\
\partial_t &\text{otherwise.}
\end{cases}
\]
This vector field \(\partial _{t'}\) is the desired
modifcation of \(\partial_t\).

\subsubsection*{\it Part 2. The set of almost complex structures \(\mathcal{J}_{0\ff}.\)}  In this and the next part, we
digress to describe the criteria for membership in the
subset \(\mathcal{J} _{1\ff}\) of the set of almost complex structures on \(\bbR
\times  M\).  The criteria are of two sorts.  The first sort is a
standard genericity condition:  It requires that all irreducible,
pseudoholomorphic subvarieties be non-degenerate in the sense that a
certain associated Fredholm operator has trivial cokernel.  To say more,
suppose that \(C \subset  \bbR  \times  M\) is an irreducible,
pseudoholomorphic subvariety. There exists a \textit{model curve\/} for
\(C\), denoted \(C_{0}\). This is a smooth, complex curve  with a
smooth, pseudoholomorphic map \[
\varphi _{0}: C_{0} \to  M,\]
which is 1-1 on the complement of a finite set, onto \(C\).
The subvariety \(C\) is said to be \textit{immersed\/}
when \(\varphi _{0}\) is an immersion.  

The operator in question has a relatively straightforward description
in the case when \(C\) is immersed.  Supposing that such is the case,
then a normal bundle \(\textup{N} \to C_{0}\) is defined so that it restricts
to any open set \(U \subset C_{0}\) that is embedded by \(\varphi _{0}\)
as the pull-back via \(\varphi _{0}\) of the normal bundle of
\(\varphi _{0}(U)\subset C\) in \(\bbR \times M\).  The almost complex structure \(J\)
endows \(\textup{N}\) with the structure of a complex line bundle.  This said,
view \(\textup{N}\) henceforth as a complex line bundle.  Then there exists on
\(C_{0}\) a canonical, \(\bbR\)-linear Fredholm operator that maps
sections over \(C_{0}\) of \(\textup{N}\) to those of \(\textup{N} \otimes
T^{0,1}C_{0}\).  This operator is denoted as \(\mathcal{D}_{C}\).  Its
action on a section \(\eta\) of \(\textup{N}\) can be written as
\begin{equation}\label{(2.4)}
 \mathcal{D} _{C}\eta
= \bar\partial\eta+ \nu _{C}\eta + \mu_{C}\bar\eta,
\end{equation}
where the notation is as follows: Fix a Hermitian metric on \(\textup{N}\) so as to
give the latter a holomorphic structure; this is used to
define \(\bar\partial\).
Meanwhile, \(\nu _{C}\) and \(\mu _{C}\) are respective sections of
\(T^{0,1}C_{0}\) and \(\textup{N}^{\otimes 2} \otimes T^{0,1}C_{0}\) that are defined using the
1-jet along \(C\) of \(J\).  Note that a change of Hermitian metric is
compensated by a change in what is meant by \(\nu _{C}\).  This is why
(\ref{(2.4)}) is canonical.  The operator \(\mathcal{D} _{C}\) arises
in the following way: Fix \(\varepsilon>0\) and let \(\eta\) denote a
section of \(\textup{N}\). The composition of first
\(\varepsilon\eta\) and then the metric’s exponential map defines a
deformation of the subvariety \(C\). This deformation is
pseudoholomorphic to first order in \(\varepsilon\) if and only if
\(\mathcal{D}_C\eta=0\). 
The assumption that all periodic orbits are
non-degenerate implies that \(\mathcal{D} _{C}\) maps \(L^{2}_{1}(C_{0}; \textup{N})\) to
\(L^{2}(C_{0}; \textup{N} \otimes T^{0,1}C_{0})\) as a Fredholm operator.  The
\textit{kernel\/} of \(\mathcal{D}_{C}\) is used here to denote the space of
sections in \(L^{2}_{1}(C_{0}; \textup{N})\) that are annihilated by \(\mathcal{D}_{C}\).
The \textit{cokernel\/} of \(\mathcal{D}_{C}\) is used here to denote its
\(L^{2}(C_{0}; \textup{N}^{\otimes 2} \otimes T^{0,1}C_{0})\) cokernel.  This is the kernel
of the formal adjoint.  Note also that the kernel and cokernel
elements are smooth.  In any event, the \textit{index} of \(\mathcal{D} _{C}\)
denotes its index as a Fredholm operator from \(L^{2}_{1}(C_{0}; \textup{N})\) to
\(L^{2}(C_{0}; \textup{N} \otimes T^{0,1}C_{0})\).  This Fredholm version of \(\mathcal{D}
_{C}\) is the Fredholm operator that was alluded to at the outset.
There is an analogous Fredholm version of \(\mathcal{D}_{C}\) in the case when
\(\varphi _{0}\) is not an immersion; it is described for example in
Section 4 of \cite{HT1}.

The interest in these operators stems from the following fact: The set
of irreducible, pseudoholomorphic subvarieties has a topology whereby
a neighborhood of any given such subvariety \(C\) is homeomorphic to the
zero locus of a smooth map from a ball in
\(\ker (\mathcal{D}_{C})\) to \(\cok (\mathcal{D}_{C})\).

A by-now standard argument using the Smale-Sard theorem can be used to
prove the following: 
\begin{lemma}
There is a residual subset 
\(\mathcal{J} _{0\ff} \subset   \mathcal{J} _{\ff }\) of almost complex structures with the following property: If
\(J \in \mathcal{J} _{0\ff}\), and if \(C\) is an irreducible, \(J\)-pseudoholomorphic
subvariety, then \(\mathcal{D}_{C}\) has trivial cokernel.  
\end{lemma}
In terms of
this, the first criterion for \(J\) to be in \(\mathcal{J}_{1\ff}\) is: 
\begin{description}
\item[\textbf{[J1]}] 
\(J\in \mathcal{J} _{0\ff}\).
\end{description}

\subsubsection*{\it Part 3. The set of almost complex structures \(\mathcal{J}_{1\ff}\).} This part describes the
second sort of criteria for membership in \(\mathcal{J} _{1\ff}\).  Assume that
\(J \in \mathcal{J} _{0\ff}\).  The criteria here involve the ends of
pseudoholomorphic subvarieties that comprise the set \(\mathcal{M}
_{1}(\Theta , \Theta ')\) for pairs \(\Theta \), \(\Theta '\) in the
set \(\mathcal{A}\).  To say more, suppose that \(\Sigma \in \mathcal{M}
_{1}(\Theta , \Theta ')\). The fact that \(\Sigma\) is in a
1-dimensional moduli space implies that there is a unique \((C,m)\in
\Sigma\) with such that the submanifold \(C\) 
not \(\bbR\)-invariant.  Let \(\mathpzc{E}\subset C\)
denote an end.  Reintroduce the notation \(q_{\mathpzc{E}}\), \(\lambda_{q_{\mathpzc{E}}}\), and
\(\varsigma_{q_{\mathpzc{E}}}\) from the discussion of (\ref{(1.3)}); 
and define a subset \(\op{div}_{\mathpzc{E}} \subset \{1,\cdots , q_{\mathpzc{E}}\}\) as follows: 
 An integer \(q\) is in \(\mbox{div}_{\mathpzc{E}}\) if either of the
 following is true:
\begin{itemize}
\item[(1)] \(q = q_{\mathpzc{E}}\),
\item[(2)] \(q\) is a proper divisor of \(q_{\mathpzc{E}}\), and there is a \((2\pi
q)\)-periodic eigenfunction \(\varsigma_q\) of the operator in (\ref{(1.4)}) with
eigenvalue \(\lambda _{q}\) of the same sign as \(\lambda_{q_\mathpzc{E}}\), such that \(|\lambda _{q_\mathpzc{E}}|>| \lambda
_{q} |> 0\).  
\end{itemize}
Note that the sign in item (2) above is positive or negative
respectively when \(\mathpzc{E}\) is a positive or negative end.
We now state the other criteria for \(J\) to be in \(\mathcal{J}
_{1\ff}\):

For any choice of \(\Theta\), \(\Theta '\), \(\Sigma \), \(C\)
described above,
\begin{description}\label{(2.5)}
\item[\textbf{[J2a]}] \(\mbox{div}_{\mathpzc{E} } =
\{q_{\mathpzc{E}}\}\) for any end \(\mathpzc{E}\) of \(C\).  
\item[\textbf{[J2b]}]  Fix any two distinct ends of \(C\), \(\mathpzc{E}\)
  and \(\mathpzc{E}'\), which are both positive or both negative such
  that their respective constant \(s \in \bbR  \) slices converge as \(|s| \to \infty  \) to the same
periodic orbit in \(M\), and such that \(q_{\mathpzc{E}} =
q_{\mathpzc{E} '}=:q\). Fix respective \(2\pi q\)-periodic
eigenfunctions  \(\varsigma_{q_{\mathpzc{E}}}\),
\(\varsigma_{q_{\mathpzc{E}'}}\) for the \(\mathpzc{E}\) and \(\mathpzc{E}'\)'s versions
  of the operator (\ref{(1.4)}). Then \(\varsigma
  _{q_{\mathpzc{E}}}(t) \neq \varsigma_{q_{\mathpzc{E}'}}(t+2\pi k)\) for any integer \(k\).
\end{description}
\begin{defn}
\(\mathcal{J}_{1\ff}\) is the set of almost complex structures \(J\)
satisfying all three conditions [J1], [J2a], [J2b] above.
\end{defn}

The arguments in Sections 3 and 4 of \cite{HT1} can be adapted with only
notational changes to prove that the set \(\mathcal{J} _{1\ff}\) is a residual
subset of \(\mathcal{J} _{0\ff}\). By way of a parenthetical remark
with regards to [J2a], Section 3 in \cite{HT1} shows that
\(\lambda_{q_{\mathpzc{E}}}\) is the
smallest positive/largest negative eigenvalue of a $2\pi
q_{\mathpzc{E}}$-periodic eigenfunction. Moreover, the condition that
  $\Sigma $ have periodic Floer homology index 1 implies the
  following: Let $q$ denote a proper divisor of $q_{\mathpzc{E}}$, and
  let $\lambda $ denote the eigenvalue of a $2\pi q$-periodic eigenfunction with the sign as \(\lambda_{q_{\mathpzc{E}}}\). Then \(|\lambda|>|\lambda_{q_{\mathpzc{E}}}|\).

\subsubsection*{\it Part 4. The modified almost complex structure \(J'\).} 
The next lemma describes a modified version of the
given almost complex structure, \(J\); this is to go with the modified
version of \(\ff\)
provided by Lemma~\ref{lem: 2.1} and the modified version of
\(\partial_t\) provided by Lemma~\ref{lem: 2.2}.

\begin{lemma}\label{lem: 2.3}
  Given an integer \(N\), there exists \(\kappa > 1\) with the following
  significance: Fix \(\delta < \kappa ^{-1}\) and let \(\ff'\) denote a
  map given by Lemma~\ref{lem: 2.1}.  Define \(\partial _{t'}\) from
  \(\ff'\) as done subsequent to Lemma~\ref{lem: 2.2}.  Then there is an
  almost complex structure, \(J'\), on \(T(\bbR \times M)\) with the
  properties listed below.
\begin{itemize}
\item  \(J'\partial _{s} = \partial _{t'}\)
\item  \(J'\) is tamed by \(ds \wedge dt + w _{\ff}\).
\item  Fix \(\gamma \in \Lambda _{N}\) and
extend Lemma~\ref{lem: 2.2}'s map \(\varphi_\gamma \) as a map 
\[
\hat{\varphi}_\gamma\co \bbR
\times (S^1 \times D) \to \bbR \times M, \quad (\sigma , (\tau , z))
\mapsto (s = q\sigma , \varphi_\gamma (\tau , z)).
\]  
Then \(J'(\hat{\varphi}_\gamma)_*\partial _{z}
=(\hat{\varphi}_\gamma)_*i\partial _{z}\) 
on a neighborhood of the \(z = 0\) locus.
\item  \(J' \in \mathcal{J} _{1\ff'}\).
\item  For any pair of periodic Floer homology generators, there is a canonical 1-1 correspondence
between the \(J\)-pseudoholomorphic curves
that contribute to the periodic Floer homology differential
and those that are \(J'\)-pseudoholomorphic.  This correspondence is
such that partnered curves contribute identical \(\pm 1\) weights to the
respective differentials.
\end{itemize}
\end{lemma}

This lemma is proved in Section~\ref{sec:2(d)}.

With \(N\) given and a choice of \(\delta\) obeying the conditions of
Lemmas~\ref{lem: 2.1}--\ref{lem: 2.3}, fix
\(\ff'\) as in Lemma~\ref{lem: 2.1} and \(J'\) as in Lemma~\ref{lem: 2.3}.  The pair \((\ff', J')\) is said
to be an {\em \((\delta , N)\)-approximation to \((f, J)\)}.  The following is a
direct corollary of Lemmas~\ref{lem: 2.1} and \ref{lem: 2.3}:

\begin{prop}\label{prop: 2.4}
Fix a pair \((\ff, J)\) with \(\ff\co F \to   F\) satisfying
(\ref{f-assumption}) and \(J
\in  \mathcal{J} _{1\ff}\). Use the latter to define the chain complex
and differential for periodic Floer homology.  Then there exists an integer
\(N_* \geq  1\) such that \(\forall N \geq  N_*\), there exists an
\(\delta _{N}\in \bbR\) with the following significance:  

Given any \(\delta
\in  (0, \delta _{N})\) and an \((\delta , N)\)-approximation, \((\ff',
J')\), to \((\ff, J)\), the tautological identification of the
respective generators induces a canonical isomorphism from the \(\ff\)-version of the periodic
Floer chain complex to the \(\ff'\)-version, that intertwines the
corresponding differentials.  In particular, the two versions of
periodic Floer homology are canonically isomorphic.
\end{prop}

{\it Remark.}  The definition of periodic Floer
homology that was given originally in \cite{H1}  restricted the
pair \((\ff, J)\) to those where \(\ff\) obeys the fifth item of
Lemma~\ref{lem: 2.1} and
where \(J\) near each \(\gamma
\in  \Lambda _{N}\) version of \(\bbR
\times  \gamma  \) is constant when written in a suitable coordinate
chart.  This understood, Proposition~\ref{prop: 2.4} can be viewed as the
assertion that the periodic Floer homology for any pair \((\ff, J)\) can be
computed using a pair of the sort described in \cite{H1}.

\subsection{Preliminary constructions}\label{sec:2(c)}         

The proofs of Lemmas \ref{lem: 2.1}-\ref{lem: 2.3} require some
preliminary constructions which will be used in the definitions of
\(\ff'\) and \(\varphi_\gamma\). 
These are given in the subsequent two parts of this subsection.

\subsubsection*{\it Part 1. The tubular neighborhood maps \(\varphi_\gamma\).} Fix \(\gamma
\in  \Lambda _{N}\).  This first part
describes a tubular neighborhood map for \(\gamma\)  that is shown in
Section~\ref{sec:2(d)} 
to satisfy the requirements of Lemma \ref{lem: 2.2}.  In what follows, \(q\) denotes
the period of \(\gamma\) .

Let \(z = (z_{1}, z_{2})\) denote the symplectic coordinates centered on
the corresponding fixed point \(\pt\) of \(\ff\) as described in Lemma
\ref{lem: 2.1}.  Use
these coordinates to identify a neighborhood of \(\pt  \in  F\) with a small
radius disk \(D' \subset  \bbR ^{2}\) about the origin.  This then
gives a smooth map \(\varphi _{0}\co [0, 2\pi ] \times
D' \to  M\) that embeds \([0, 2\pi  )\times
D'\) and is such that
\[
 (\varphi _{0})_*\partial_{\tau } = \partial_t  \quad \text{and} \quad
\varphi _{0}\,(2\pi , z)= \varphi _{0}\, (0, f^{q}(z)).  
\]
The map \(z \mapsto \ff^{q}(z)\) maps \(D'\) to \(\bbR ^{2}\) and fixes the origin.  It
also preserves the area-form \(w_{F}\).  This understood, there
exists a homotopy of
area-preserving maps, \[
\big\{\psi_\tau (z)\co D''\to \bbR^2\big\}_{\tau  \in [0,2\pi ]},\]
where \(D''\) is a disk \(D'' \subset D'\), and \(\psi_\tau\) is such
that:
\begin{itemize}
\item \(\psi_\tau (0) = 0\) for all \(\tau\), 
\item \(\psi_\tau(z) = z\) near
\(\tau = 0\), and 
\item \(\psi_\tau (z) = f^{-q}(z)\) near \(\tau = 2\pi \).  
\end{itemize}
Let \(D^{*}
\subset D''\) denote a small radius disk about the origin such that \(D^*\)
is in \(\psi_\tau (D'')\) for all \(\tau\) .  Define the map
\[
\varphi_1\co [0, 2\pi ]
\times D^* \to M,\quad (\tau, z)\mapsto \varphi
_{0}\, (\tau , \psi_\tau (z)).\]  
By construction, \(\varphi _{1}(0, z) =
\varphi _{1}(2\pi , z)\).  As a consequence, \(\varphi _{1}\) defines a map
from \(S^1 \times D^{*}\) into \(M\) such that \(\varphi _{1}^*dt = q d\tau \).

By assumption, \(J\) maps the kernel of \(\pi _*\) to itself along \(\gamma \).
It is also the case that the quadratic form \(w _{F}( \cdot ,
J(\cdot ))\) is positive on the kernel of \(\pi _{*}\) along \(\gamma\) .  As a
consequence, the action of \(J\) along \(\gamma\) on \({\varphi_1}_*\partial _{z}\) is such that
\begin{equation*}
(\varphi _{1}^{-1})_{*}\,J\, (\varphi _{1})_{*}\partial _{z} = i \,(1+
|\upsilon|^{2})^{1/2}\, \partial _{z} +
\upsilon\, \partial_{\bar{z}},
\end{equation*}
where \(\upsilon\) here is a smooth map from \(S^1\) to \(\bbC\) .  This
understood, we may define:
\begin{defn}
Set \(\varphi_\gamma : S^1 \times D^*
\to M\) to be the map that sends any given point \((\tau , z)\) to \(\varphi
_{1}(\tau , V_{\tau } z)\), where the map 
\(\tau \mapsto V_{\tau }\), from \(S^1\) to \(\op{Sl}(2; \bbR
)\) is chosen so that \(\varphi_\gamma\) satisfies the following
properties: 
\begin{itemize}
\item The vector field \(q^{-1}\partial_t\) is the
push-forward via \(\varphi_\gamma \) of a vector field on \(S^1 \times D^{*}\)
that can be written as
\begin{equation}\label{(2.7)}
\partial _{\tau } - 2i \, (\nu z +
\mu\bar{z}+\grr ) \, \partial _{z} + 2i\,
(\nu \bar{z}+\bar{\mu}z+\bar{\grr})\,\partial_{\bar{z}},
\end{equation}
where \(\nu\) is an \(\bbR\)-valued function on \(S^1\) and \(\mu\) is a \(\bbC\)-valued
function.  Meanwhile \(\grr\) is such that its absolute value is bounded
by \(\kappa|z|^{2}\), and its derivatives are bounded in absolute value by
\(\kappa|z|\) for a positive constant \(\kappa\).  
\item \(J(\varphi_{\gamma})_*\, \partial _{z} =
  (\varphi_{\gamma})_*i\,\partial _{z}\) along \(\gamma\).
\end{itemize}
\end{defn}
In what follows, \(F\) is identified with \(\pi ^{-1}(0) \subset M\) and
\(\varphi_\gamma\) is used to identify \(\{2\pi k/q\}\times D^* \subset S^1 \times D^{*}\) with a disk
in \(F\), for any given \(k \in \{0, \cdots , q\}\).  Granted these identifications, then \((\varphi_\gamma )^*\, f\) on a
small radius, concentric subdisk in the \(\tau = 2\pi k/q\) disk inside
\(S^1 \times D^*\) is obtained from (\ref{(2.7)}) by following the integral curves
of (\ref{(2.7)}) from this subdisk to where they intersect the \(\tau = 2\pi
(k+1)/q\) disk in \(S^1 \times D^{*}\).

\subsubsection*{\it Part 2. The maps \(\uuu_\tau: D\subset F\to F\).}  
This part constructs an \(S^1\)-parametrized family of
area-preserving maps \(\uuu_\tau\) from a certain disk in \(F\) to
\(F\).  This family is used to
construct the map \(\ff'\).  

To start, introduce \((\nu, \mu )\) to denote the pair of
functions that arise when writing \(q^{-1}(\varphi_\gamma
^{-1})_*\partial_t\) as in (\ref{(2.7)}).  There exists a 1-parameter set,
\(
\{(\nu^\lambda, \mu^\lambda)\}_{\lambda\in [0,1]},
\) 
of respectively \(\bbR\)-valued and \(\bbC\)-valued functions
on \(S^1\) satisfying the following:   
\begin{itemize}
\item \((\nu^0, \mu^0)=(\nu_\gamma , \mu_\gamma )\), and \((\nu^1,
  \mu^1)= (\nu , \mu )\).  
\item Let \(K^{\lambda } \in
\op{Sl}\,(2; \bbR )\) be defined so that given \(z \in  \bbR ^{2}\), \(K^{\lambda
}\cdot z\) is the value at \(\tau  =
2\pi  \) of the solution to the equation
\begin{equation*}
\frac{i}{2}\frac{d}{d\tau}\eta+\nu^\lambda\eta+\mu^\lambda\bar{\eta}=0,\quad \eta (0) = z.
\end{equation*}
Then each \(\lambda  \in  [0, 1]\) version of \(K^{\lambda }\) is 
conjugate in \(\op{Sl}\,(2; \bbR )\) to \(K_{\gamma }\).  
\end{itemize}
Fix such a 1-parameter
family.

Reintroduce the function \(\chi\) from Section \ref{sec:1(d)}.  Fix \(\rho_{0} \in (0,
\delta ^{16})\) and set \(\lambda^*\) to be the function on \(\bbR ^{2}\) given
by \(\lambda^*(z) = \chi\,\big(\ln|z|/\ln\rho_{0}\big)\).  Introduce respective \(\bbR\)-valued
and \(\bbC\)-valued functions on \(S^1 \times D\) given at \((\tau , z)\) by
\((\nu^{\lambda^*}, \mu^{\lambda^*})\).
With these function in hand, set
\begin{equation}\label{(2.9)}
\grr ^{(1)} = (\nu^{\lambda^*}-\nu_\gamma)\,
z+(\mu^{\lambda^*}-\mu_\gamma)\, \bar{z} +\lambda^*\grr +\gre,
\end{equation}
where \(\grr\) also comes from (\ref{(2.7)}), and where \(\gre\) is a smooth map from
\(S^1 \times D^*\) to \(\bbC\) with the following properties: 
\begin{itemize}
\item \(\gre\) is
supported where \(|z| \in (\rho_{0}^{7/16}, \rho_{0}^{5/16})\).  
\item 
\(\gre\) and its derivatives to order \(100\) are bounded by \(\rho_{0}^{1000}\).
\end{itemize}
The term \(\gre\) is used subsequently as a source of `suitably
generic' perturbations.

Choose \(\delta \) small enough so that the disk in \(\bbC\) about the origin 
 with radius \(\delta\) lies in \(D^{*}\).  There exists a constant
 \(c_0\) so that over the disk \(D
 \subset D^*\) of radius \(c_{0}^{-1}\delta_*\), the following
 family of maps \(\{\uuu^{(1)}_{\tau }: D \to
 D^*\}_{\tau \in [0,2\pi ]}\) is well-defined: Given \(z \in
 \bbC\), solve the differential equation
\begin{equation}\label{(2.10)}
\frac{i}{2}\frac{d}{d\tau}\eta+ \nu_\gamma\eta +\mu_\gamma \bar{\eta}+
\grr^{(1)}(\eta ) = 0,\quad \eta (0) = z, 
\end{equation}
and set 
\[
\uuu^{(1)}_{\tau }(z): =\eta \,(\tau ).\]  
This map has
the following three properties: 
\begin{itemize}
\item There is a disk \(D^{(1)} \subset
D\) centered at the origin of radius less than \(\rho_{0}^{1/2}\) such that
\(\uuu^{(1)}_{\tau }\) acts on a neighborhood of \(D^{(1)}\) as multiplication
by the matrix \(U(\tau )\) that is constructed in Section \ref{sec:2(a)} from the
solutions to the linear equation depicted in (\ref{(2.3)}).  
\item There is
a disk \(D^{(2)} \subset D\) of radius less than \(\rho_{0}^{1/4}\) that
contains \(D^{(1)}\) and is such that \(\uuu^{(1)}_{\tau }\) acts on a
neighorhood of \(D-D^{(2)}\) as an area-preserving map.  Furthermore, if \(k
\in \{0, \cdots , q-1\}\), then \(\uuu ^{(1)}_{\tau =2\pi k/q}\) on \(D-D^{(2)}\)
is given by \(\varphi_\gamma ^*\, \ff\) on \(\{2\pi k/q\} \times (D-D^{(2)})\).
\item Each \(\uuu ^{(1)}_{\tau }\) is nearly area-preserving in the
sense that its pull-back of the area form \(\frac{i}{2}dz\wedge d\bar{z}\)
can be written as
\begin{equation}\label{(2.11)}
\big(1 + \grw_{\tau}\big) \,\frac{i}{2}dz\wedge d\bar{z},
\end{equation}
where the error term, \(\grw_{\tau }\), is supported in \(D^{(2)}-D^{(1)}\)
and is such that \(|\grw _{\tau }| \leq c_{0}|\ln \rho _{0}|^{-1}\).  
\end{itemize}
Now set 
\[\uuu _{\tau } := \varrho_\tau \circ\uuu^{(1)}_{\tau },\]
where \(\varrho_\tau\) is 
a diffeomorphism of \(D^{*}\), which is the identity on
\(D^{(1)}\) and on \(D-D^{(2)}\) so that the result is
area-preserving.  

The map \(\varrho_{\tau}\) is defined using Moser's trick.  To elaborate, the
map \(\varrho_{\tau}\) is obtained by integrating a \(\tau\)-dependent vector
field on \(D^{(2)}-D^{(1)}\) that vanishes near the boundary.  This is to
say that 
\(\frac{\partial}{\partial_\tau} \varrho_\tau = \op{v}_{\tau}(\varrho_\tau
)\), 
where \(\tau \mapsto \op{v}_{\tau }\) is a vector field
with compact support on \(D^{(2)}-D^{(1)}\).  This vector field
\(\op{v}_{\tau}\) is chosen
so that its contraction with the form in (\ref{(2.11)}) is an anti-derivative,
\(\grb_{\tau }\), of the 2-form \(\frac{i}{2}\grw_{\tau}\, dz \wedge
d\bar{z}\) with
compact support in \(D^{(2)}-D^{(1)}\).  To say more about this
anti-derivative, introduce the cylindrical coordinates \((\rho, \theta )\)
on \(D\) by writing  \(z = \rho \, e^{i\theta}\). The 2-form
\(\frac{i}{2}\grw_{\tau}\, dz \wedge
d\bar{z}\) can be written as \(d\,(h_{\tau } d\theta \,)\), where
\begin{equation}\label{(2.12)}
h_{\tau } (\rho, \theta )= \int_0^\rho  \grw_\tau (x, \theta)\, x\, dx.
\end{equation}
The function \(h_{\tau }\) is zero on \(D^{(1)}\) and it depends only on the
angle \(\theta\) on a neighborhood of \(D-D^{(2)}\).  The fact
that \(\frac{i}{2}\, \grw_{\tau}\, dz \wedge d\bar{z}\) has
integral zero over \(D^{(2)}\) implies that \(h_{\tau }=
\frac{\partial}{\partial \theta} \, u_{\tau }\)
on \(D-D^{(2)}\) where \(u_{\tau }\) is a function only of the angle \(\theta\) .
This understood, let \(\rho_{(2)}\) denote the radius of \(D^{(2)}\), and then
define
\begin{equation}\label{(2.13)}
\grb_{\tau} = h_{\tau } \, d\theta - d\big(\, (1-\chi (|\cdot
|/\rho_{2})\, u_{\tau}\big) .
\end{equation}
This 1-form has compact support in \(D-D^{(2)}\) and is such that
\(d\grb_{\tau } = \frac{i}{2}\grw_{\tau}dz \wedge d\bar{z}\).
Note for future reference that \(|\grb_{\tau }|\) and its \(\tau\)-derivative are
both bounded by \(c_{0} |\ln \rho_{0}|^{-1} |z|\).  In addition,
\(|\nabla\grb_{\tau }|\) and its \(\tau\)-derivative are bounded by
\(c_{0}|\ln \rho_{0}|^{-1}|z|\).  It follows as a consequence that \(\op{v}_{\tau }\) and its corresponding
derivatives obey similar bounds.

\subsection{Proofs of Lemmas \ref{lem: 2.1}, \ref{lem: 2.2}
  and \ref{lem: 2.3}}\label{sec:2(d)}

We now have the ingredients to define the modified map:
\begin{defn}
Given a pair \((\delta, N)\) of a small positive number and an integer, 
let \(\ff'\co F\to F\) be the following map:
\[
\ff'=
\begin{cases}
(\varphi_\gamma ^{-1})^*\uuu_{\tau =2\pi k/q}& \text{on
  \(\varphi_\gamma \big(\{2\pi k/q\} \times D\big)\),} \\
& \quad \text{\(\forall k \in \{0,
  \cdots , q-1\}\) and \(\gamma\in \Lambda_N\)};\\
\quad \ff & \text{on the remainder of \(F\)}.
\end{cases}
\]  
In the above, \(q\) denotes the period of \(\gamma\), and 
the parameters \(\gre\), \(\rho_0\) that enter into the
definition of \(\uuu_{\tau =2\pi k/q}\) are chosen as follows:
\begin{itemize}
\item the perturbation term \(\gre\) in (\ref{(2.9)}) is chosen in a suitably
generic fashion so that all periodic points of \(\ff'\) are
non-degenerate periodic points;
\item let \(\rho_{0}=\delta ^{32}\).
\end{itemize}
\end{defn}

By construction, \(\ff'\) obeys the first two items of 
Lemma~\ref{lem: 2.1}, and it
 obeys the remaining items if all periodic orbits of \(\ff'\) with period
 \(N\) or less are in \(\ff'\)s version of \(\Lambda _{N}\). To prove this, suppose
 that \(\gamma ' \subset F\) is a periodic orbit of \(\ff'\) with period \(q' \in
 \{1, \cdots , N\}\). Let \(\pt' \in \gamma '\). Given that \(|\ff-\ff'| \leq
 c_{0}\, \rho_{0}^{1/4}\), it follows that \(\big|\dist \big(\, \ff'^{q'}(\pt'), \pt'\big)\big| \leq
 c_{0}\, \rho _{0}^{1/4}\) also. Thus, if \(\delta < c_{0}^{-1}\), then there
 is an orbit \(\gamma \subset \Lambda _{N}\) of order a divisor, \(q\), of \(q'\),
 and a point \(\pt  \in \gamma \) with the following property: Each \(k \in
 \{1,\cdots , q'\}\) version of \(\ff^{k}(\pt )\) has distance \(c_{0}\, \rho _{0}^{1/4}\)
 or less from both \(\ff^{k}(\pt')\) and \((\ff')^{k}(\pt')\). This understood, \(\pt'\) and
 its \(\ff'\) iterates are in the image via the tubular neighborhood map
 \(\varphi_\gamma \) of the set \(\bigcup_{1\leq k\leq q} \{2\pi k/q\} \times D\). No
 generality is lost by assuming that \(\pt'\) is in the
 \(\varphi_\gamma \) image
 of a point \((0, z) \in S^1 \times D\). 

Let \(\eta (\tau )\) denote
 the solution to (\ref{(2.10)}) with \(\eta (0) = z\) and set the
 constant 
\[ \lambda = \lambda^*(z).\] 
Given that \(\eta\) obeys (\ref{(2.10)}), it follows that
\begin{equation}\label{(2.14)}
c_{0}^{-1}|z| \leq |\eta| \leq c_{0}\, |z| 
\end{equation}
at all \(\tau \in S^1\). As a consequence, \(|\lambda^*(\eta ) - \lambda|
\leq c_{0}|\ln \delta|^{-1}\) at each \(\tau \in S^1\). This then implies
that the equation for \(\eta\) can be written as
\begin{equation}\label{(2.15)}
\frac{i}{2}\frac{d}{d\tau}\eta +\nu^\lambda\eta + \mu^\lambda\bar{\eta}
+\grx(\eta ) = 0,
\end{equation}
where \(|\grx (\eta)| \leq c_{0}|\ln \delta|^{-1}|\eta|\). Integrating
(\ref{(2.15)}) finds that \[
\big|\, \big(1 - (K^{\lambda})^{q'/q}\big) z\, \big| \leq
c_{0}|\ln \delta|^{-1}|z|.\] 
This last inequality requires that \(z = 0\) when
\(\delta\) is small by virtue of the fact that \(K^{\lambda }\) is conjugate
in \(\op{Sl}\, (2; \bbR )\) to \(K_{\gamma }\) and \(\det \big(1 - (K_{\gamma })^{q'/q}\big) \neq
0\) for all integer multiples \(q'\) of \(q\).

 To obtain Lemma \ref{lem: 2.2}, take \(\grr'\)  to be the version of \(\grr ^{(1)}\)
 given by (\ref{(2.9)}) that is used to construct \(\ff'\). This done, then the
 lemma only summarizes definitions of \(\varphi_\gamma \) and the map \(\ff'\)
 given in the previous section.

 The proof of Lemma \ref{lem: 2.3} differs only in notation from the proof of
 Proposition 2.5 from \cite{T1} given in the latter's Appendix. The reader
 is referred to this same Appendix for the details. The five remarks
 that follows serve only to point out certain key points that enter
 the proof.

\medbreak

\subsubsection*{\it Remark 1.} It is a relatively straightforward
exercise with matrices to construct an almost complex structure \(J'\)
that obeys all but the fifth item of Lemma \ref{lem: 2.3}, and which has the
following additional properties: First, \(|J - J'|\) is zero except in the
image via the maps from \(\{\, \varphi_\gamma \}_{\gamma\in
  \Lambda_N}\) of the set
in \(S^1 \times D\) where \(|z| \leq \rho_{0}^{1/4}\). Second, the difference
between \(J\) and \(J'\) on such an image is such that \(|J - J'| \leq
c_{0}\, |z|\). In general \(|\nabla(J - J')|\) will be \(\mathcal{O}
\, (1)\). This last fact
complicates a direct perturbation theoretic construction of the
desired pairing between \(J\) and \(J'\) pseudoholomorphic curves.

\subsubsection*{\it Remark 2.} To see how this problem is
addressed, keep in mind that the \(\ff\) and \(\ff'\) versions of the set \(\Lambda
_{N}\) agree; as a consequence, so do the two versions of the generating
set \(\mathcal{A}\) for the periodic Floer homology chain complex. This
understood, fix a pair \(\Theta _{-}\) and \(\Theta _{+}\) from \(\mathcal{A}\) so as to
consider the \(J\) and \(J'\) versions of the space \(\mathcal{M}_{1}(\Theta _{+},
\Theta _{-})\). The relatively large size of \(|\nabla (J - J')|\) precludes
a direct comparison between the \(J\) and \(J'\) versions of this set. Even
so, an indirect comparison can be made as follows: Fix some large
integer \(Q\) so as to compare successive versions of \(\mathcal{M} _{1}(\Theta
_{+}, \Theta _{-})\) as defined by an ordered set, \(\{J_{k}\}_{0\leq k\leq
Q}\), of almost complex structures with the following properties:
\begin{itemize}
\item  \(J_{0} = J\), and \(J_{Q}\) is the desired
\(J'\).
\item  \(J_{k} = J\) on the complement of the images
via the maps \(\{\varphi_\gamma\}_{\gamma\in \Lambda_N}\) of where \(|z| \leq \rho_{0}^{1/4}\).
\item  \(|J_{k} - J_{k-1}| \leq c_{0}\, Q^{-1}|z|\)
 on the image of \(\bbR \times S^1 \times D\) of any map
from the set \(\{\varphi_\gamma\}_{\gamma\in \Lambda_N}\).
\item  \(|\nabla(J_{k} - J_{k-1})| \leq c_{0} \, Q^{-1}\).
\item  Each \(\gamma \in \Lambda _{N}\)
version of the cylinder \(\bbR \times \gamma \) is
\(J_{k}\)-pseudoholomorphic. 
\end{itemize}

\subsubsection*{\it Remark 3.} Perturbation-theoretic tools are used
to construct a 1-1 correspondence between the pseudoholomorphic curves
in the respective versions of \(\mathcal{M}_{1}(\Theta _{+}, \Theta _{-}\)) as
defined by the almost complex structures at successive
steps. Likewise, perturbation theory is used to prove that paired
almost pseudoholomorphic curves have equal \(\pm 1\) weights of the sort
used to compute the embedded contact homology differential.

The use of perturbative tools exploits four key facts. Here are the first three:
\begin{itemize}
\item[\bf Fact 1]  The factor \(Q^{-1}\) can be very small.
\item[\bf Fact 2]  Each \(\gamma \in \Lambda _{N}\)
version of the cylinder \(\bbR \times \gamma \) is
\(J_{k}\)-pseudoholomorphic for each \(k\).
\item[\bf Fact 3]  Successive versions of \(|J_{k+1} - J_{k}|\)
can have support in successively smaller neighborhoods of
\(\bigcup _{\gamma\in \Lambda_N} \gamma\); and these neighborhoods can be made as small as desired,
independent of \(Q\) and \(\delta\) .
\end{itemize}
To say something about the fourth key fact, suppose that \(C\) is a
non-\(\bbR\) invariant subvariety that appears in some element from the
\(J_{k}\) version of \(\mathcal{M} _{1}(\Theta _{+}, \Theta _{-})\). Any given end
\(\mathpzc{E} \subset C\) at large \(|s|\) is described by (\ref{(1.3)}) but with the
functions \((\nu , \mu )\) dependent on the index \(k\). This understood, here
is the fourth fact:
\begin{itemize}
\item[\bf Fact 4] There is a lower bound, independent of  \(k\), \(Q\) and
\(\delta\) to the absolute value of the eigenvalue  \(\lambda _{q_\mathpzc{E}}\)
that can appear in (\ref{(1.3)}).
\end{itemize}

\subsubsection*{\it Remark 4.} These four facts are relevant
because any given non-\(\bbR\) invariant, \(J_{k}\)-pseudoholomorphic
subvariety will intersect the small neighborhoods in Fact 3 above
 in only two
ways. First, such a subvariety, \(C\), can intersect this neighborhood
near the points where it intersects a cylinder of the form \(\bbR \times
\gamma \) with \(\gamma \in \Lambda _{N}\). There are but a finite set of
points where \(C\) intersects any given \(\bbR\)-invariant cylinder. In fact,
there is a bound on the number of such points that is essentially
topological and so independent of \(C\) and of the index \(k\). In particular,
if the modification from one step to the next is done in a very small
neighborhood of the cylinders \(\{\bbR \times \gamma\}_{\gamma \in
  \Lambda _{N}}\), then the
region where \(J_{k+1} \neq J_{k}\) will intersect all but the very large
\(|s|\) part of \(C\) in a set of at most \(c_{0}^{-1}\) disks, each with very
small area. This has the following consequences: The subvariety \(C\) is
\(J_{k+1}\)-pseudoholomorphic except at very large \(|s|\), and except in a
union of at most \(c_{0}^{-1}\) disks, each very small area. Moreover, \(C\)
is nearly \(J_{k+1}\)-pseudoholomorphic in each such disk if \(Q\) is large
and \(\delta\) is small, both in a pointwise sense and in an \(L^{2}\)
sense.

The subvariety \(C\) also intersects the region where \(J_{k+1} \neq J_{k}\)
at points far out on its ends, thus at very large values of \(|s|\). The
area of intersection here is, of course, infinite. However, the
relevant versions of (\ref{(1.3)}) are used to prove that this is of no
ultimate consequence. In particular, \(C\) is nearly \(J_{k+1}\)-pseudoholomorphic in a pointwise sense and also in an \(L^{2}\) sense if \(Q\)
is large and \(\delta\) is small.

\subsubsection*{\it Remark 5.} This final remark concerns a most
crucial point in the constructions from the Appendix of \cite{T1}: The
perturbative approach from this appendix will construct a partner to \(C\)
from the \(J_{k+1}\) version of \(\mathcal{M} _{1}(\Theta _{+}, \Theta _{-})\)
provided that the operator \(\mathcal{D}_{C}\) from (\ref{(2.4)}) is invertible. The four
facts stated above, and Fact 4 
in particular, guarantee that such is
the case when \(Q\) is large and \(\delta\) is very small.

\section{Constructing Seiberg-Witten solutions on \(M\)}\label{sec: 3}
\setcounter{equation}{0}

 Define \(M\) from a given area preserving map \(\ff\co F \to F\). The purpose of
 this section is to explain how the constructions in Section 3 of \cite{T2}
 are used to associate solutions to certain versions of the
 Seiberg-Witten equations to finite sets with any given element a pair
 consisting of a closed integral curve of \(\partial_t\) and a positive
 integer. Such an association can also be obtained from the
 constructions in the article \(\op{Gr}\Rightarrow \op{SW}\) from \cite{T5}. The
 constructions from Section 3 of \cite{T2} differ with regard to certain
 details from those in the \(\op{Gr}\Rightarrow \op{SW}\) article, but they are
 identical in spirit.

 The relation between the solutions to the Seiberg-Witten equations
 and sets of closed integral curves is described in Theorem
 \ref{thm:3.2} of the
 upcoming Section \ref{sec: 3(d)}. The intervening subsections set the stage.

\subsection{Vortices}\label{sec: 3(a)}

 This subsection summarizes some of what is said in Section 2 of \cite{T2}
 about the vortex equations on \(\bbC\). These are equations for a pair
 \((A, \alpha )\) where \(A\) is a connection on the product \(\bbC\)-bundle over
 \(\bbC \) and \(\alpha\) is a section of this bundle. These equations read
\begin{equation}\label{(3.1)}
*F_{A} = -i\, (1 - |\alpha|^{2}) \quad \textit{and}\quad
\bar{\partial}_A \alpha = 0.
\end{equation}
Here, \(F_{A}\) denotes the curvature 2-form of \(A\)
and  \(\bar{\partial}_A\) is the d-bar operator defined by \(A\). These equations are augmented with the
constraint that \((1 - |\alpha|^{2})\) be integrable. Solutions \((A, \alpha
)\) and \((A', \alpha ')\) are deemed gauge equivalent when there exists a
map \(u: \bbC \to S^1\) such that \(A' = A - u^{-1}du\) and \(\alpha ' = u\alpha\).

The space of equivalence classes of solutions has path components
\(\{\grC_{m}\} _{m=0,1,\cdots }\). The space \(\grC_{0}\) consists of the gauge
equivalence classes of the pair \((A_{I}, \alpha = 1)\) where
\(A_{I}\) is the product connection on \(\bbC \times \bbC \). Meanwhile, any
given \(m \geq 1\) version of \(\grC_{m}\) consists of the gauge
equivalence classes of solutions to (\ref{(3.1)}) with \(\alpha\)
having \(m\) zeros counting multiplicity. Each such \(\grC_{m}\)
has the structure of a
complex manifold. The complex structure is such that the holomorphic
tangent space to \(\grC_{m}\) at the equivalence class of some given
solution \(\grc = (A, \alpha )\) has a canonical identification with the
\(L^{2}\) kernel of the operator: 
\begin{equation}\label{(3.2)}\begin{split}
\vartheta _{\grc }\co C^{\infty}(\bbC ; \bbC \oplus \bbC ) & \to C^{\infty}(\bbC ; \bbC \oplus \bbC ) , \\
(x, \iota ) & \mapsto \big( \partial x + 2^{-1/2}\bar{\alpha}\iota, \bar{\partial}_A \iota+
2^{-1/2}\alpha x\big) .
\end{split}
\end{equation}
Here, \(\partial\) denotes $\frac{\partial}{\partial z}$.
 The \(L^{2}\) kernel of \(\vartheta _{\grc }\) is denoted by \(\ker\,(\vartheta _{\grc
})\) in what follows. The space \(\grC_{m}\) is biholomorphic to \(\bbC ^{m}\);
and the collection of functions \(\{\sigma _{q}\}_{1\leq q\leq m}\) on
\(\grC_{m}\) given by
\begin{equation}\label{(3.3)}
\sigma _{q} = \frac{1}{2\pi}\int_\bbC z^q (1-|\alpha|^2)
\end{equation}
give a set of holomorphic coordinates on \(\grC_{m}\). 

 The inner product on \(L^{2}(\bbC ; \bbC \oplus \bbC )\) gives \(\grC_{m}\) a
 K\"ahler metric. To be explicit, the induced Hermitian inner product
 on \(T_{1,0}\grC_{m}\) is such that the square of the norm of any given
 element \((x, \iota )\) in \(\ker (\vartheta _{\grc })\) is
\begin{equation}\label{(3.4)}
\frac{1}{\pi} \int_\bbC  (|x|^2+|\iota|^2).
\end{equation}
Except for the \(n = 1\) case, this K\"ahler metric is not the pull-back
of the standard metric on \(\bbC ^{m}\) via the holomorphic identification
that is defined by the functions \((\sigma _{1}, \cdots , \sigma _{m})\).

 The action of \(S^1\) on \(\bbC\) as multiplication by the unit complex
 numbers induces an isometric action of \(S^1\) on \(\grC_{m}\). This action
 has a unique fixed point, this being the unique solution to (\ref{(3.1)})
 whose version of \(\alpha\) vanishes at the origin with degree \(m\). This
 point in \(\grC_{m}\) is called the symmetric vortex.

 Of interest in what follows are certain sorts of dynamical flows on
 \(\grC_{m}\). To say more, fix a pair \((\nu , \mu )\) of maps from \(S^1\) to
 \(\bbR\) and \(\bbC\) . Use the latter to define the real valued function
\begin{equation}\label{(3.5)}
\hhh =\hhh_{(\nu, \mu)}:= \frac{1}{4\pi} \int_\bbC \big( 2\nu
|z|^2+(\mu\bar{z}^2+\bar{\mu}z^2)\big) (1-|\alpha|^2).
\end{equation}
on \(S^1 \times \grC_{m}\). View this as a 1-parameter family of functions
on \(\grC_{m}\). Use the K\"ahler 2-form to define the corresponding 
1-parameter family of Hamiltonian vector fields. Of interest are the
closed, integral curves of this family. These are maps \(\grc : S^1 \to
\grC_{m}\) that obey at each \(\tau \in S^1\) the equation
\begin{equation}\label{(3.6)}
\frac{i}{2}\grc '+ \nabla ^{(1,0)}\hhh|_{\grc } = 0 ,
\end{equation}
where \(\grc '\) is shorthand for the \((1, 0)\) part of \(\grc_*\frac{d}{d\tau}\), and
where \(\nabla^{(1,0)}\hhh\) denotes the \((1, 0)\) part of the gradient of
\(\hhh\).

 A solution to (\ref{(3.6)}) is said to be {\em nondegenerate} when a
 certain auxilliary, self-adjoint operator has trivial kernel. To
 elaborate, suppose for the moment that \(\grc \co S^1 \to \grC_{m}\) is any
 given map. Associate to the latter the symmetric operator
\begin{equation}\label{(3.7)}
 \zeta \mapsto
\frac{i}{2}\nabla_{\tau}\, \zeta + \big(\nabla _{\zeta_\bbR}\nabla^{(1,0)}\hhh \big)\big|_{\grc }
\end{equation}
on \(C^{\infty }(S^1; \grc ^*T_{1,0}\grC_{n})\). Here, \(\nabla _{\tau}\) denotes
the covariant derivative on \(C^{\infty }(S^1; \grc ^*T_{1,0}\grC_{n})\) as
defined by the pull back of the Levi-Civita connection on
\(T_{1,0}\grC_{n}\). Meanwhile,
\((\nabla _{\zeta_\bbR}\nabla^{(1,0)}\hhh )|_{\grc }\) denotes the covariant derivative at \(\grc\) along
the vector defined by \(\zeta\) in \(T\grC_{n}|_{\grc }\) of the vector field
\(\nabla ^{(1,0)}\hhh \in C^{\infty }(\grC_{n}; T_{1,0}\grC_{n})\). A
solution to (\ref{(3.6)}) is {\em non-degenerate} when the operator depicted in
(\ref{(3.7)}) has trivial kernel.

The subsequent lemma plays a key role in what follows. 
\begin{lemma}\label{lem:3.1} 
Fix a positive integer \(m\) and a
pair \((\nu , \mu )\co S^1 \to \bbR \times \bbC  \) whose version of
(\ref{(1.4)}) has trivial \(2\pi q\) periodic kernel for each \(q \in \{1,
\cdots , m\}\). Then the space of solutions to the corresponding
version of (\ref{(3.6)}) is compact; and it is finite if all solutions are
non-degenerate.
\end{lemma}
The proof of this lemma is identical to the proof given in Section 2b
of \cite{T4} of the latter's Proposition 2.1.

All solutions to (\ref{(3.6)}) are known in special cases. What follows
 summarizes what is said in Lemmas 2.1-2.3 of \cite{T2}.
\begin{lemma}\label{(3.8)}
\begin{itemize}
 \item[{\rm (a)}]  Fix a pair \((\nu , \mu )\) with \(\nu\) a real
valued function on \(S^1\) and \(\mu\) a \(\bbC\)-valued function. Use
this pair to define the operator, \(\scrP\co C^\infty (\bbR, \bbC)\to C^\infty (\bbR, \bbC)\)
by the rule \(\eta\mapsto \frac{i}{2}
\frac{d}{d\tau}\eta+\nu\eta+\mu\bar{\eta}\). Assume that 0 is the
only \(2\pi\)-periodic element in the kernel of \(\scrP\).
Then there is a
unique solution to the corresponding \(m = 1\) version of (\ref{(3.6)}),
which is the constant map to the symmetric vortex. This solution is
non-degenerate.
\item[{\rm (b)}]  Fix an irrational number \(\rrr\) and
set \((\nu = \rrr , \mu= 0)\). Then there is a unique solution in \(m \geq 1\)
version of \(\grC_{m}\) to the corresponding version of
(\ref{(3.6)}). This is the constant map from \(S^1 \) to the symmetric
vortex \(\grc _{0} \in \grC_{m}\). This solution is
non-degenerate.
\item[{\rm (c)}]  Fix \(\rrr \in \frac{1}{2}\bbZ\)
 and a positive number \(\zzz\). Set \(\big(\nu=\frac{1}{2}
\rrr , \mu=-\frac{i}{4\pi}\zzz \, e^{2i\rrr\tau}\big)\).
Then
there is a unique solution to the corresponding \(m = 1\) version
of (\ref{(3.6)}). This is the map to the symmetric vortex. This solution is
non-degenerate. On the other hand, there are no solutions to the \(m >
1\) versions. 
\end{itemize}
\end{lemma}

Fix \(m \geq 1\) and suppose that \((\nu , \mu )\) are as described in
Lemma \ref{lem:3.1}. In such circumstance, any given non-degenerate
solution to the corresponding version of (\ref{(3.6)}) can be assigned
an integer degree. To do this, use the given pair \((\nu , \mu )\) to
define the operator \(\scrP\) in Lemma \ref{(3.8)}. Let \(\eta(\tau)\)
denote a non-trivial
element in the kernel of \(\scrP\). Then \(\eta(\tau)\) can be written as \[
\eta (\tau ) = U (\tau)\, (\eta _{0}),\] 
where \(U (\tau )\) is an \(\bbR\)-linear map from
\(\bbC\) to \(\bbC\) that defines a matrix in \(\op{Sl} \, (2; \bbR )\) when
\(\bbC\)  is written as \(\bbR ^{2}\). The pair \((\nu , \mu )\) is
said to be {\em elliptic} with rotation number a given real number
\(\rrr\not\in \bbZ\) when \(U \) on \([0, 2\pi ]\) defines a path in \(\op{Sl} \, (2; \bbR )\) that is homotopic in \(\op{Sl} \, (2; \bbR )\) via a homotopy \(\{U_x(\cdot)\}_{x\in[0,1]}\) with the
following three properties: The matrix \(U_x(0)\) is the identity, the matrix \(U_x(2\pi)\) is
conjugate in \(\op{Sl} \, (2; \bbR )\) to \(U(2\pi)\), and the path \(U_{x=1}(\cdot)\) is the rotation with smooth angle
parametrization \(\tau\mapsto2\pi\rrr\tau\). It is assumed in what follows that \(\rrr\) is irrational when
\((\nu, \mu)\) are elliptic.
The pair \((\nu , \mu )\) is said to be
{\em hyperbolic} when the absolute value of the trace of \(U(2\pi )\) is
greater than 2. In this case, there is a 1-parameter family of pairs
of function \[
\{(\nu _{x}, \mu _{x})\co S^1 \to \bbR \times \bbC \}_{x\in [0,1]}\] 
with the following three properties: 
\begin{itemize}
\item  \((\nu _{0}, \mu _{0}) = (\nu , \mu )\). 
\item Each pair \((\nu
  _{x}, \mu _{x})\) has a corresponding \(U_{x} (2\pi )\) whose trace
  has absolute value at least 2. 
\item The path \(\tau\mapsto U_{x=1}(\tau)\) as defined by the pair
  \((\nu_1, \mu_1)\) is the rotation with smooth angle parametrization
  \(\tau\mapsto2\pi\rrr\tau\) with rotation number \(\rrr \in 
\frac{1}{2}\bbZ \). 
\end{itemize}
Granted the preceding, let \[
\grc _{0}\co S^1 \to \grC_{m}\] 
denote the
constant map to the symmetric vortex in \(\grC_{m}\). Associate to \(\grc
_{0}\) the operator 
\[
\scrR \co C^{\infty }(S^1; T\grC_{m}|_{\grc_0})\to C^{\infty }(S^1; T\grC_{m}|_{\grc_0}),
\] 
that is defined as follows: 
\begin{itemize}
\item If \((\nu , \mu )\) is elliptic with rotation number
\(\rrr\), this operator is the version of (\ref{(3.7)}) with
\(\hhh=\hhh_{(\frac{1}{2}\rrr, 0)}\)
\item If \((\nu , \mu )\) is hyperbolic with rotation number \(\rrr\), the
operator \(\scrR\) is the version of (\ref{(3.7)}) with 
\(\hhh=\hhh_{(\frac{1}{2}(\rrr -\rrr '), 0)}\), where 
\(\rrr' > 0\) is an irrational
number such that \(m\rrr ' \ll 1\). 
\end{itemize}
Lemma 2.4 in \cite{T2} asserts that the
operator \(\scrR\) has trivial kernel in these cases.
\begin{defn}
Let \((\nu , \mu )\) be as in Lemma 3.1, and suppose that \(\grc \co S^1 \to
\grC_{m}\) is nondegenerate in the sense that (\ref{(3.7)}) has trivial
kernel. Then {\it the degree of \(\grc\)}, denoted \[
\deg_{\grC}(\grc)=\deg_{\grC}\, (\grc, (\nu, \mu))\] 
in the following, is the spectral flow from the operator \(\scrR\) to \(\grc\)'s
version of 
(\ref{(3.7)}).
\end{defn}
 For the definition of this spectral flow, see e.g. Section 2 of
 \cite{T8}). 
\begin{lemma}\label{lem:vortex-degree}
Let \((\nu, \mu)\) be the pair \((\nu_\gamma, \mu_\gamma)\) that
appears in Definition \ref{(2.2)}. Fix \(m\in\{1, 2, \ldots\}\) when
\((\nu, \mu)\) are elliptic and set \(m=1\) when \(\gamma\) is
hyperbolic. Then the map to the symmetric vortex
\(
\grc_0\co S^1\to \grC_m
\)
has degree
\[
\deg_\grC\, (\grc_0, (\nu_\gamma, \mu_\gamma))=\begin{cases}
0 &\text{when \(\gamma\) is elliptic or negative hyperbolic,}\\
1 & \text{when \(\gamma\) is positive hyperbolic.}
\end{cases}
\] 
\end{lemma}
\pf When \(\gamma\) is elliptic, this is true by definition. When
\(\gamma\) is hyperbolic and \(m=1\), the operator in (\ref{(3.7)})
when written using the (\ref{(3.3)})'s coordinate \(\sigma_1\) to identify
\(\grC_1\) with \(\bbC\) is the operator \(\scrP\) that appears in
Lemma \ref{(3.8)}. This understood, the observation about the spectral flow in
the negative hyperbolic case follows from the fact that the version of
\(\scrP\) with \((\nu, \mu)=(\pi, 0)\) has trivial \(2\pi\)-periodic
kernel. Meanwhile, the calculation of the assertion about the spectral
flow in the positive hyperbolic case follows by direct computation
using the \([0, 1]\)-parametrized path \(x\mapsto \scrP _x\) with
\(\scrP _x\) defined by \((\nu_x , \mu_x)=(-(1-x)\delta, -i\delta)\) with \(\delta\in (0,1)\).
\medbreak

\subsection{Periodic orbits and vortices}\label{sec: 3(b)} 

 Return now to the context where \(F\) is a surface with area form \(w
 _{F}\) and \(\ff\co F \to F\) is an area-preserving map. Let \(M\) denote the
 resulting 3-manifold. Suppose that \(\Theta\) is a finite set whose elements are described next. A given element is a pair \((\gamma , m)\) with \(\gamma\) a closed integral
 curve of \(\partial_t\) and \(m\) a positive integer. Assume that distinct
 pairs from \(\Theta\) have distinct integral curves. Note that
 \(\Theta\) can be taken to be the empty set.

Let \((\gamma , m) \in \Theta \). Fix a tubular neighborhood map
\(\varphi_\gamma : S^1 \times D_{\gamma } \to M\) with \(D_{\gamma } \subset
\bbC \) a disk about the origin. The map \(\varphi_\gamma \) should be chosen
to obey the conditions set forth by the first two items and the
fourth item of Lemma \ref{lem: 2.2}. Use this map \(\varphi_\gamma \)
to define
corresponding functions \((\nu , \mu )\); and 
let \(\grC_{(\gamma ,m)}\) denote the set of solutions to the corresponding
version (\ref{(3.6)}) with \(\hhh=\hhh_{(\nu, \mu)}\). Let \[
\grC\Theta :=\times _{(\gamma
,m) \in  \Theta } \grC_{(\gamma ,m)}.\] 
Use
\(\grC\Theta^*\subset \grC\Theta \) to
denote the subset of elements with the property that constituent maps
\(\{\grc _{\gamma }\co S^1 \to \grC_{m}\}_{(\gamma ,m) \in  \Theta }\)
are non-degenerate.

Now suppose that a class \(\Gamma \in H_{1}(M; \bbZ )\) has been
specified. As in the introduction, use \(d_{\Gamma }\) to denote the
intersection number between \(\Gamma\) and a typical fiber of the
projection map \(\pi\) . Let \(p\) denote the greatest integer divisor of the
class \(c_\Gamma \in H^{2}(M; \bbZ )\). 

Let \(\mathcal{Z}\) denote the set whose elements are described next. A given element is a set \(\Theta\) as just
described, but constrained so that the integral chain \(\sum_{(\gamma
,m) \in  \Theta } m\gamma \) represents the class \(\Gamma\). Note that
\(\mathcal{Z} = \emptyset \) if \(d_{\Gamma } \leq 0\) and if \(\Gamma
=0\), then \(\mathcal{Z}\) has only the element \(\Theta=\emptyset\). Hutchings specifies in
Definition 1.5 of \cite{H1} a relative \(\bbZ /p\bbZ \)-degree for each element
\(\Theta \in \mathcal{Z} \). This degree is relative in the following sense: A
reference element \(\Theta^* \in \mathcal{Z}\) is assigned degree 0
when \(\mathcal{Z}\neq\emptyset\), then each
\(\Theta \in \mathcal{Z}\) has a well defined degree assignment in \(\bbZ
/p\bbZ \). This degree is denoted by \(\deg_{pFh}(\Theta )\). If
\(\Theta '\) is a second element in \(\mathcal{Z}\), then the difference
\(\deg_{pFh}(\Theta ) - \deg_{pFh}(\Theta ')\) is the \(\bbZ /p\bbZ
\) reduction
of what was denoted in Section \ref{sec:1(a)} by \(I\, (\Theta , \Theta '; Z)\), where \(Z\)
is a 2-chain in \(M\) with \(\partial Z = \sum_{(\gamma ,m)\in \Theta} m\gamma - \sum_{(\gamma ,m) \in  \Theta '} m\gamma\).

 The relative \(\bbZ /p\bbZ \)-degree for a given \(\Theta \in \mathcal{Z}\)
 can be used to assign a relative \(\bbZ /p\bbZ \)-degree to each \(x =
 \{\grc _{\gamma }\}_{(\gamma ,m)\in  \Theta } \in \grC\Theta^*\),
 this being \[
\deg_{\grC \mathcal{Z} }(x ):=
\deg_{pFh}(\Theta ) + \sum_{(\gamma ,m) \in  \Theta }\deg_{\grC}(\grc
_{\gamma }). 
\]

\subsection{The Seiberg-Witten equations}\label{sec: 3(c)} 

 This subsection sets the stage for the definition of the relevant
 versions of the Seiberg-Witten equations. 

\subsubsection*{The almost complex structure \(J\).}
 Suppose that a non-degenerate \(\ff\) is given to define \(M\), and
 that \(J\in \mathcal{J}_\ff\) is an \(\bbR\) -invariant
 almost complex structure on \(\bbR \times M\). Recall that by
 definition, such a \(J\) satisfies:
\begin{itemize}
 \item[\rm (1)]  \(J\partial _{s} = \partial_t\) where
\(\partial_t\) is a vector field on \(M\) with \(\langle
dt, \partial_t\rangle = 1\) whose first return is \(\ff\).
\item[\rm (2)]  \(J\) is tamed by the form \(\Omega_F = ds \wedge  dt + w _{F}\).
\end{itemize}

Because of  Item (1) above,
the almost complex structure \(J\) is
determined by its action on \(\ker (\pi _*)\). 
Its action here determines a
pair \(  (\grj , a )\) where \(\grj\) is an endomorphism of \(\ker(\pi _*)\)
 with square \(-1\) and \(a\) is a section of \(\Hom\, (\ker(\pi _{*}); \bbR
)\). These are defined so that
\begin{equation}\label{(3.10)}
 Jv = \grj v + a (v)\, \partial _{s} + a (\grj
 v)\, \partial_t\quad 
\text{for \(v \in \ker \, (\pi _*)\).}
\end{equation}
The assertion that \(J\) is tamed by \(w _{0}\) is equivalent to the
assertion that \[
|a (v)|^{2} < 4w_{F}(v, \grj v).\]

It is important for what follows to note that
\begin{equation}\label{(3.11)}
 J\, (v + a (v)\partial_t) = \grj v + a (\grj v)\, \partial_t .
\end{equation}

As a consequence, \(J\) preserves the subbundle
\begin{equation}\label{(3.12)}
 K^{-1} = \big\{v + a (v)\, \partial_t\co v \in \ker\, (\pi _*)\big\} .
\end{equation}

By design, this subbundle has the structure of a complex line bundle
over \(M\). The bundle is of course isomorphic to the kernel of \(\pi _*\)
with its complex structure defined by \(\grj\). The bundle \(K^{-1}\) is the
kernel of the 1-form
\begin{equation}\label{(3.13)}
 \textsf{a} = dt - a 
\end{equation}
where \(a\)  is viewed here and subsequently as a 1-form on \(M\) that annihilates \(\partial_t\). 

\subsubsection*{The metric \(g\).}
 The almost complex structure \(J\) and the area form \(w_F\) together
 define a metric \(g\) on \(M\) as
 follows: 
\begin{itemize}
\item The
 vector field \(\partial_t\) has norm 1 and is orthogonal to
 \(K^{-1}\). 
\item The inner product on \(K^{-1}\) is given by \(w
 _{F}(\cdot , \grj (\cdot ))\). 
\end{itemize}
With this definition, \(g\) is such that \[
|\textsf{a}| = 1\quad  \text{and} \quad 
*\textsf{a} = w_{F}.\] 
Use this metric to define the oriented orthonormal frame
 bundle of \(M\) and to define the corresponding \(\Spin^{\bbC }\) lifts of
 this frame bundle.

\subsubsection*{The \(\Spin^\bbC\) structure.}
 Fix a \(\Spin^{\bbC }\)-bundle \(\pi : \scrF \to M\), and let \[
\bbS = \scrF \times _{U(2)} \bbC \to M\] 
denote the corresponding 
\(\bbC ^{2}\)-bundle. Clifford multiplication on \(\bbS\) by the 1-form \(\textsf{a}\) from (\ref{(3.13)})
 has eigenvalues \(\pm i\), and so splits \(\bbS\) as a direct sum of the
 corresponding two eigenbundles. This splitting is written as
\begin{equation}\label{(3.14)}
\bbS = E \oplus  E\otimes K^{-1}, 
\end{equation}
where the convention has \(\cl (\textsf{a})\) acting as \(i\) on the left most
summand. The assignment of \(c_{1}(E)\) to the \(\Spin^{\bbC }\)-bundle defines
a 1-1 correspondence between the set of equivalence classes of
\(\Spin^{\bbC }\)-bundles over \(M\) and the elements in \(H^{2}(M; \bbZ )\). The
canonical \(\Spin^{\bbC }\)-structure is one whose spinor bundle
is splits as \(\bbS_{I} = I_{\bbC } \oplus K^{-1}\), where \(I_{\bbC }\) here
denotes a complex line bundle that is isomorphic to the product bundle
\(M \times \bbC \).
\bigbreak

\subsubsection*{\(\Spin^\bbC\) connections and Dirac operators.}
 As explained in Section 1c of the article \(\op{SW}\Rightarrow \op{Gr}\) from \cite{T5},
 there is a unique connection, \(A_{K}\), on \(K^{-1}\) which is characterized
 as follows: Fix a unit norm section, \(1_{\bbC }\), of the bundle
 \(I_{\bbC}\), and let \(A_{I}\) denote the unique Hermitian connection on \(I_{\bbC }\)
 that makes \(1_{\bbC }\) covariantly constant. Use the pair \(A_{K}\) and
 \(A_{I}\) to define a connection on \(\det\, (\bbS)\). Then the section \(\psi_{I} = (1_{\bbC }, 0)\) obeys the Dirac equation.

Because of (\ref{(3.14)}), the connection \(A_{K}\) and a Hermitian
connection on \(E\) together define a
connection \(\bbA\)
on \(\det\, (\bbS) = E^{\otimes 2}\otimes K^{-1}\). Conversely, a connection on
the latter bundle defines with \(A_{K}\) a Hermitian connection on \(E\). This
understood, the Seiberg-Witten equations are viewed in what follows as
equations for a pair \((A, \psi )\) of connection on \(E\) and section of
\(\bbS\) . The corresponding Dirac operator is denoted now as \(D_{A}\). The
space of smooth, Hermitian connections on \(E\) is denoted in what follows
by \(\op{Conn}(E)\).

\subsubsection*{The 3-dimensional Seiberg-Witten equations.}
The versions of the Seiberg-Witten equations of interest are
parametrized by a number \(r > 1\). Fix such an \(r\). With the metric
and \(\Spin^\bbC\)-structure chosen previously, the corresponding
equations for a pair \((A, \psi )\) from the space \(\op{Conn}(E) \times
C^{\infty }(M; \bbS)\) read
\begin{equation}\label{(3.15)}
 B_{A} = r\,  (\psi^{\dag }\tau\psi - i\textsf{a})
-\frac{1}{2} B_{A_K}-\frac{i}{2}*\wp,\quad 
\text{and} \quad D_{A}\psi = 0.
\end{equation}
Here, \(B_{A}\) denotes the Hodge star of \(A\)'s curvature 2-form,
and \(\wp\) is a closed 2-form in the cohomology class \(2\pi c_1(\grs)\). Solutions
\((A, \psi )\) and \((A', \psi ')\) are gauge equivalent if there is a smooth
map, \(u\), from \(M\) to \(S^1\) such that \[
A' = A - u^{-1}du\quad  \text{and}\quad  \psi ' = u\psi.
\]
Because \(*\textsf{a} = w_\ff\), the equations in (\ref{(3.15)}) correspond to a
version of (\ref{(1.7)}) with \[
\varpi=2r w _{\ff}+\wp, \quad \Psi=(2r)^{1/2}\psi,  
\]
and \(\bbA\) defined from \(A\) and \(K\) as above. 
When \(r\) is fixed, we shall use 
\((\bbA, \Psi)\) and \((A, \psi)\) interchangably to denote a 
{\em configuration}, and write the configuration space as 
\(\op{Conn}(E)\times C^\infty(M, \bbS)\). 
Use \(\scrC^{r}\) in what follows to denote the space of gauge equivalence classes
of solutions to (\ref{(3.15)}).

 \subsubsection*{The relative grading.}
Each configuration \(\grc = (A, \psi )\) 
has an associated elliptic,
 symmetric operator 
\[
\grL_\grc \co C^{\infty }(M; iT^*M \oplus \bbS \oplus i\bbR)\to C^{\infty }(M; iT^*M \oplus \bbS \oplus i\bbR),
\] 
which sends any given section
 \((b, \eta , \phi )\) of \(iT^*M \oplus \bbS \oplus i\bbR \) to the section
 whose respective \(iT^*M\), \(\bbS\) and \(i\bbR \) components are
\begin{equation}\label{(3.16)}
\begin{cases} 
&*db - d\phi - 2^{-1/2} r^{1/2} (\psi ^{\dag}\tau\eta + \eta ^{\dag }\tau\psi ),\\
& D_{A}\eta + 2^{1/2}r^{1/2} (\cl (b)\, \psi + \phi\psi ),\\
& *d*b - 2^{-1/2 }r^{1/2} (\eta ^{\dag }\psi - \psi ^{\dag }\eta ).
\end{cases}
\end{equation}
(Note that despite its appearance, \(\grL_\grc\) depends only on the configuration
\(\grc=(\bbA, \Psi)\), not on \(r\)).
A configuraton  \(\grc\) 
is said to be
{\em nondegenerate} when \(\grL_\grc\)
has trivial cokernel. Let \(\grc_-\), \(\grc_+\) be nondegenerate configurations.
Given a path \(\grd(s)\) of configurations from \(\grc_-\) to
\(\grc_+\), let \(\scrI(\grd)\) denote the spectral flow of the
associated family
of operators from \(\grL_{\grc_-}\) to \(\grL_{\grc_+}\), defined as
in \cite{T8}.
For \(\grc_+=u\cdot\grc_-\), one has \begin{equation}\label{I-period}
\scrI(\grd)=\langle c_1(\grs),
  [u]\rangle,
\end{equation}
where \(u\in C^\infty(M; S^1)\) is a gauge
transformation and \([u]\) denotes  its cohomology class:
\[
[u]:=u^{-1}du/(2\pi i).
\]
Thus, as promised in Section \ref{sec:1(b)}, 
we may write \(\scrI(\grc_-, \grc_+; \grh)=\scrI(\grd)\), which
depends only on the gauge equivalence classes of \(\grc_-, \grc_+\),
and the relative homotopy class \(\grh\) of \(\grd\).
Fix a nondegenrate reference  configuration \(\grc_E:=(A_{E}, \psi
_{E})\). Given a nondegenerate configuration \(\grc\), let
\[
\deg_{SW}(\grc ):=\scrI(\grc, \grc_E)\in \bbZ/p\bbZ.
\]
This depends only on the gauge equivalence class of \(\grc\).

\bigbreak

\subsection{Periodic orbits and Seiberg-Witten solutions} 
\label{sec: 3(d)}

 Fix a class \(\Gamma \in H^{1}(M; \bbR )\). Let \(\grs_\Gamma\) denote the
 \(\Spin^{\bbC }\) structure whose spinor bundle splits as
 in (\ref{(3.14)}) with the first Chern class of the bundle \(E\) equal to the
 Poincar\'e dual of \(\Gamma\) . The theorem that follows refers to the
 versions of (\ref{(3.15)}) as defined for this particular \(\Spin^{\bbC }\)
 structure.

\begin{thm}\label{thm:3.2}
Let \(\pmb{\mu}_P=\{(F, w_F), \ff, \Gamma, J\}\) and
\(\pmb{\mu}_S=\{M, \grs_\Gamma, \varpi_r, g, \grq\}\) be corresponding
periodic Floer homology and Seiberg-Witten-Floer cohomology data sets
as in the statement of Theorem \ref{thm:1}, and let 
\(r \geq 1\). Use \(\scrC^{r}\) to denote
the space of gauge equivalence classes of solutions to the corresponding version of
(\ref{(3.15)}). Then: 
\begin{itemize}
\item[{\rm (1)}]  If \(d_{\Gamma } \leq 0\) and \(\Gamma
\neq 0\), then \(\scrC ^{r} = \emptyset  \) for all sufficiently large \(r\). 
\item[{\rm (2)}]  If \(\Gamma = 0\), then \(\scrC^r\)
consists of a single element for all sufficiently large
\(r\).
\item[{\rm (3)}]  Suppose that \(d_{\Gamma } > 0\). Let 
\[\grX=
\begin{cases}
 \bigcup _{\Theta  \in \mathcal{Z} } \grC\Theta  &\text{if 
 \, \(\grC\Theta ^{*} = \grC\Theta  \) for all \(\Theta \in \mathcal{Z}\)},\\ 
\text{a finite subset of  \(\bigcup _{\Theta  \in  \mathcal{Z}}\grC\Theta^*\)} &\text{otherwise}.
\end{cases} 
\]
There exists,
for all sufficiently large \(r\), a 1-1 map \[\Phi ^{r}: \grX \to
\scrC^r\] 
with the following properties:
\begin{itemize}
\item[{\rm (a)}]  The image of \(\Phi ^{r}\) consists of
non-degenerate solutions to (\ref{(3.15)}).
\item[{\rm (b)}]  If \(x_{-}, x_{+} \in \grX \), then \[\deg_{SW}(\Phi ^{r}(x_{+})) -
\deg_{SW}(\Phi ^{r}(x _{-})) =\ deg_{\grC \mathcal{Z} }(x _{-}) -
\deg_{\grC\mathcal{Z} }(x _{+}).\]
\item[{\rm (c)}]  If \(\grC\Theta ^{*} = \grC\Theta \)
for all \(\Theta \in \mathcal{Z}  \) and \(\grX = \bigcup _{\Theta  \in \mathcal{Z} } \grC\Theta \), then \(\Phi ^{r}\) is also onto
\(\scrC^r\).
\end{itemize}
\end{itemize}
\end{thm}

\bigbreak

\subsubsection*{Proof of Theorem \ref{thm:3.2}:} 
There are five parts to the
proof of this theorem. Part 1 constructs the solution whose existence
is asserted by the second item. Part 2 constructs the map \(\Phi ^{r}\)
from \(\grX\) into \(\scrC^r\) and proves that it is one to one. Part 3 proves
Items (3a) and
(3b) of the assertions. Part 4 establishes some basic a priori bounds
for solutions of (\ref{(3.15)}). These bounds are then used in Part 5 to prove Item (1), the uniqueness assertion of Item (2), and
Item (3c).

\subsubsection*{\it Part 1.} The solution whose existence is asserted by
Item (2) can be constructed by copying almost verbatim what is done in
the first parts of Section 2d of \cite{T9}, with the following
replacements:
\BTitem\label{HP-replace}
\item Replace the contact form \(\textsf{a}\) in \cite{T9} by the
  1-form \(\textsf{a}\) introduced in (\ref{(3.13)}).
\item Replace \(B_{A_K}\) in \cite{T9} by 
\[
B_0:=B_{A_K}+i*\wp.
\]
\ETitem
These replacements do not affect the arguments in Section 2d of
 \cite{T9}; and these prove the following:
\begin{lemma}\label{lem:3.3} 
Suppose that the \(\Spin^{\bbC}\)
structure is the canonical one. Then there exists a constant \(\kappa > 1\)
such that for all \(r \geq \kappa \), there is a solution,
\((A_{0}, \psi _{0})\), to (\ref{(3.15)}) with the following properties:
\begin{itemize}
\item \( |\psi_{0}| \geq 1 - \kappa r^{-1/2}\) on the whole of \(M\).
\item  If \((A, \psi )\) is a solution to
(\ref{(3.15)}) and if \(|\psi| \geq 1 - \kappa ^{-1}\) on the whole of
\(M\), then \((A, \psi )\) is gauge equivalent to \((A_{0}, \psi _{0})\).
\end{itemize}
\end{lemma}

\subsubsection*{\it Part 2.} The map \(\Phi^r\) is constructed by copying what
is done in Section 3 of \cite{T2}. The latter constructs the analog of
the map \(\Phi^r\) in the following context: First the orbits that
appear in any given set from \(\mathcal{Z}\) are closed integral curves of
the Reeb vector field of a contact 1-from on \(M\). Second, that same
contact 1-form is used as the 1-form \(\textsf{a}\) to define
(\ref{(3.15)}). 
With the same sort of replacements as (\ref{HP-replace}), 
these constructions can be applied in the present context
with only cosmetic changes. The proof that the current version of
\(\Phi^r\) is injective is obtained by copying almost verbatim what is
said in Section 3g of \cite{T2} to prove the latter's Theorem
1.1. Note in this regard that the function \(\textsc{e}\) that appears
in Equation (1.12) and Therorem 1.1 of \cite{T2} has an analog here;
the definition is the same: 
\begin{equation}\label{(3.17)}
\smE (A)=i \int_M \textsf{a}\wedge * B_A,
\end{equation}

What is denoted by \(\ell_\gamma\) in \cite{T2}
corresponds here to \(2\pi q\), where \(q\) is the period of \(\gamma\). As a
consequence, \[
\sum_{(\gamma ,m)\in \Theta }m \ell_\gamma=2\pi d_{\Gamma }.\] 
More is said in the upcoming Part 4 of the proof of
Theorem \ref{thm:3.2} about the values
of the function \(\textsc{E}\) in (\ref{(3.17)}) on the space \(\scrC^r\).

\subsubsection*{\it Part 3.} 
What is asserted by Items (3a) and (3b) of Theorem \ref{thm:3.2}
has its analog in Theorem 1.1 of \cite{T3}. 
Again, with the same sort of substitution as (\ref{HP-replace}), 
the arguments in Section 2 of \cite{T3} can be applied here
with only cosmetic changes to prove items (3a) and (3b) of Theorem \ref{thm:3.2}.

\subsubsection*{\it Part 4.} This part of the proof establishes the key a priori
bounds on the behavior of solutions to (3.15). The first of these is
an analog of Lemma 2.3 in \cite{T4}:
\begin{lemma}\label{lem:3.4} 
There exists a constant \(\kappa > 1\) with the
following significance: Suppose that \(r \geq \kappa  \) and that
\((A, \psi = (\alpha , \beta ))\) is a solution to (\ref{(3.15)}). Then:
\begin{itemize}
\item[{\rm (1)}]  \(|\alpha| \leq 1 + \kappa r^{-1}\).
\item[{\rm (2)}] \( |\beta|^{2} \leq \kappa r^{-1} (1 - |\alpha|^{2}) + \kappa ^{2} r^{-2}\).
\item[{\rm (3)}]  \(|\nabla _{A}\alpha|^{2} \leq \kappa r (1 - |\alpha|^{2}) + \kappa ^{2}\) .
\item[{\rm (4)}]  \(|\nabla  _{A}\beta|^{2} \leq \kappa (1 - |\alpha|^{2}) + \kappa ^{2} r^{-1}\).
\end{itemize}
In addition, for each integer \(q \geq 1\), there exists a
constant \(\kappa _{q}\) such that when \(r \geq \kappa \),
then 
\begin{itemize}
\item[{\rm (5)}] \(  |\nabla _{A}^{q}\alpha| +
  r^{1/2}|\nabla_{A}^{q}\beta| \leq \kappa _{q} \, r^{q/2}\) .
\end{itemize}
\end{lemma}

\bigbreak

\subsubsection*{Proof of Lemma \ref{lem:3.4}:} The first two item and the
final item are proved by copying what is done in the proofs of
Lemmas 2.2-2.4 of \cite{T7}. The third and fourth items are proved by
copying what is done in the proof of the analogous assertions of Lemma
2.3 in \cite{T4}.

 The next a priori bound concerns the integral of\( r (1 - |\alpha|^{2})\)
 over \(M\). Note that a bound on this integral gives a bound on the value
 of the function \(\smE\) from (\ref{(3.17)}). This follows by using (\ref{(3.15)}) and the
 first two items in Lemma \ref{lem:3.4} to write
\begin{equation}\label{(3.18)}
*B_{A} = -i\, \big(r \, (1 - |\alpha|^{2}) + \gre _{1}\big) w_{\ff} + \gre _{2},
\end{equation}
where \(\gre _{1}\) and \(\gre _{2}\) are such that:
\[\begin{split}
|\gre _{1}| &  \leq c_{0} \, (1 - |\alpha|^{2}) + c_{0}^{2} r^{-1}; \\
 |\gre _{2}| & \leq r^{1/2}|\, (1 - |\alpha|^{2})\, | +
c_{0} \quad \text{and \(\textsf{a} \wedge \gre _{2} = 0\).}
\end{split}
\]
 
\begin{lemma}\label{lem:3.5}
 There exists a constant \(\kappa > 1\) with the
following significance: Suppose that \(r \geq \kappa  \) and that
\((A, \psi )\) is a solution to (\ref{(3.15)}). Then
\[
r \int_M \big| 1-|\alpha|^2\big|\leq   \kappa .  
\]
\end{lemma}

\subsubsection*{Proof of Lemma \ref{lem:3.5}:}
Since \(\frac{i}{2\pi} F_{A}\)
represents the first Chern class of the bundle \(E\) in de Rham
cohomology, it follows that
\begin{equation}\label{(3.19)}
i\int_M dt\wedge * B_A= 4\pi^{2} d_{\Gamma } .
\end{equation}
Meanwhile, the first two items in Lemma \ref{lem:3.4} together with 
(\ref{(3.18)}) imply that the
integral on the left hand side is no less than
\begin{equation}\label{(3.20)}
(1 - c_{0}^{-1}) \, r\int_M (1-|\alpha|^2) -c_0.
\end{equation}
This bounds the integral over \(M\) of \(r (1 - |\alpha|^{2})\) by an \(r\) and
\((A, \psi )\) independent constant. Together with the first item in Lemma
\ref{lem:3.4}, such a bound supplies the same sort of bound on the
integral over \(M\) of \(|(1 - |\alpha|^{2})|\).

\subsubsection*{\it Part 5.} 
This part proves assertions (1), (3c), and the uniqueness assertion in
item (2) 
of Theorem \ref{thm:3.2}.
To start, consider a sequence \(\{(r_{n}, (A_{n}, \psi _{n} = (\alpha _{n}, \beta
_{n})))\}_{n=1,2,\cdots }\) such that \(\{r_{n}\}_{n=1,2,\cdots }\) is
unbounded, and such that \((A_{n}, \psi _{n})\) obeys the \(r = r_{n}\)
version of (\ref{(3.15)}). If \(c_{1}(E) = 0\), assume in addition that \((A_{n},
\psi _{n})\) is not the solution \((A_{0}, \psi _{0})\) given by Lemma \ref{lem:3.3}.

Given the bound supplied by Lemma \ref{lem:3.5}, the argument from Section 6d of
\cite{T7} (or part I of \cite{T5}) 
can be used with only notational modifications to find a finite
set \(\Theta = \{(\gamma , m)\}\) with the following properties: 
\begin{itemize}
\item
\(\gamma\) is a closed integral curve of \(\partial_t\) and \(m\) is a positive
integer. 
\item \(\sum_{(\gamma ,m)\in \Theta } m\gamma\)
defines a chain that represents \(\Gamma\) in \(H_{1}(M; \bbZ )\). 
\item The
loci \(\alpha _{n}^{-1}(0)\) converges geometrically to 
\(\sum_{(\gamma ,m)\in \Theta } m\gamma \) in the following sense:
Fix \(\delta > 0\). If \(n\) is large, then \(|\alpha _{n}| \geq 1 -
\delta \) at distances greater than \(c_{0}r^{-1/2}\) from \(\bigcup
_{(\gamma ,m)\in \Theta } \gamma \). Meanwhile, for any given
\((\gamma , m) \in \Theta \), the restriction of \(\alpha _{n}/|\alpha _{n}|\)
to the boundary of the \(\varphi_\gamma \) image of a disk \(\{\tau \} \times D\)
has degree \(m\).
\item \(\frac{i}{2\pi}F_A\) converges as a current to \(\sum_{(\gamma ,m)\in \Theta } m\gamma \).
\end{itemize}
These last conclusions imply that \(d_{\Gamma } > 0\) if
\(\Theta\neq\{\emptyset\}\) and thus
assertion (1) of Theorem \ref{thm:3.2}. If \(\Theta=\{\emptyset\}\), these conclusions with Lemma \ref{lem:3.3}
imply the uniqueness assertion (2) of Theorem \ref{thm:3.2}.

Granted \(d_{\Gamma } > 0\), then Item (3c) of Theorem
\ref{thm:3.2} is proved by copying almost verbatim the arguments that are given
in Section 2a of \cite{T4}. The latter prove an analog for the case when
the 1-form \(\textsf{a}\) in (\ref{(3.15)}) is a contact 1-form. This contact condition
plays no essential role in the arguments in Section 2a of 
\cite{T4}.

\section{Pseudoholomorphic curves and Seiberg-Witten solutions on
\(\bbR \times M\)}\label{sec: 4}
\setcounter{equation}{0}

Fix a nondegenerate, area-preserving map \(\ff\co F \to F\) so as to
construct the manifold \(M\), and a monotone class
\(\Gamma \in H_{1}(M; \bbZ )\) with \(d_{\Gamma } > 0\). With \(\Gamma\) given,
define \(\mathcal{Z}\) as in Section \ref{sec: 3(b)}. 
Fix next a \(J \in \mathcal{J} _{1\ff}\). 
By Lemmas \ref{lem: 2.1}-\ref{lem: 2.3}, we may assume without loss of
generality that the pair \((\ff, J)\)
obeys the following constraint: Let \((\gamma , m)\)
denote any given element from \(\mathcal{Z}\) and let \(q\) denote the period of
\(\gamma\) . There is a disk \(D \subset \bbC \) about the origin and a tubular
neighborhood map \(\varphi_\gamma \co  S^1 \times D \to M\) such that:
\BTitem\label{(4.1)}
\item  The curve \(\gamma  \) is \(\varphi_\gamma (\cdot , 0)\).
\item  The vector field \(q^{-1}\partial_t\)
is the pushforward from \(S^1 \times D\) of 
\[
\partial _{\tau } - 2i \, (\nu_\gamma z +\mu_\gamma
\bar{z}+\grr)\, \partial_z +2i \, (\nu_\gamma\bar{z}+\bar{\mu}_\gamma
z+\bar{\grr})\, \partial_{\bar{z}},
\]
where \((\nu_\gamma , \mu_\gamma )\) are given in (\ref{(2.2)})
and where \(\grr\) is zero on a neighborhood of \(S^1 \times \{0\}\).
\item  The almost complex structure \(J\) is
such that 
\(J\cdot \varphi_{\gamma*}\partial _{z} = \varphi_{\gamma*}i\partial _{z}\)
on the  image of a neighborhood of
\(S^1 \times \{0\}\) under \(\varphi_\gamma  \).
\ETitem

 Recall from Section \ref{sec:1(a)} that
\(\mathcal{A}\subset \mathcal{Z}\) denotes the subset that consists of those
 elements \(\Theta \) that do not partner hyperbolic integral curves of
 \(\partial_t\) with numbers greater than 1. Given (\ref{(4.1)}), the following is
 now a consequence of Lemma \ref{(3.8)}:
\BTitem\label{(4.2)}
\item  \(\grC\Theta^* = \grC\Theta  \) for all \(\Theta
\in \mathcal{Z}\).
\item  \(\grC\Theta \) is a single point if \(\Theta \in
\mathcal{A}\), and \(\grC\Theta =\emptyset\) otherwise.
\ETitem

Use the almost complex structure \(J\) to define the metric on \(M\) as done
in the preceding section. Use the \(\Spin^{\bbC }\)-structure on \(M\) whose
spinor bundle splits as in (\ref{(3.14)}) with \(c_{1}(E)=e_\Gamma\). 
For \(r \geq 1\), let \(\scrC^r\) denote the corresponding space
of gauge equivalence classes to (\ref{(3.15)}). Because of (\ref{(4.2)}), all
sufficiently large \(r\) versions of Theorem \ref{thm:3.2}'s map \(\Phi ^{r}\) supply a
1-1 correspondence between the set \(\mathcal{A}\) and the set \(\scrC^r\).

The body of this section concerns a map \(\Psi^r\) that associates an instanton
solution to versions of the Seiberg-Witten equations on \(\bbR \times M\)
to a weighted, pseudoholomorphic subvariety of the sort that is used
to define the differential in periodic Floer homology. This map is
described in the upcoming Section \ref{sec:4(c)}. 
The intervening sections set the
stage for what is said in Section \ref{sec:4(c)} about this map.

 \subsection{Pseudoholomorphic subvarieties in \(\bbR
 \times M\)}\label{sec: 4(a)}

This subsection discusses certain issues concerning the differential for
periodic Floer homology. To start,  recall that the periodic Floer 
homology chain complex is the free \(\bbZ\)-module generated by
equivalence classes of pairs of the form \((\Theta , \gro)\) where \(\Theta
\in \mathcal{A}\) and \(\gro\) is an ordering of the set of pairs in \(\Theta\) of the
form \((\gamma , 1)\) with \(\gamma\) a hyperbolic Reeb orbit with even
rotation number. The equivalence relation has \((\Theta , \gro) \sim
(-1)^{\sigma } (\Theta , \gro')\) where \(\sigma\) is 0 or 1, which is the
parity of the permutation that takes \(\gro\) to \(\gro'\).

Let \((\Theta _{-}, \gro_{-})\) and \((\Theta _{+}, \gro_{+})\) denote
generators of the periodic Floer chain complex. Hutchings
observed that the extra data given by \(\gro_{-}\) and \(\gro_{+}\) can be
used to associate a weight from the set \(\{\pm 1\}\) to each component of
\(\mathcal{M} _{1}(\Theta _{+}, \Theta _{-})\); this is the weight
\(\sigma\) that is used
to define via (\ref{(1.5)}) the periodic Floer homology differential. The
definition of this \(\pm 1\) weight is identical in all respects to an
analogous definition that is used to define embedded contact
homology. The latter is given in Section 9 of \cite{HT1} and also described
in Section 3b of \cite{T3}. What follows immediately is a two part summary
of how the \(\pm 1\) weight is defined.

\subsubsection*{\it Part 1.} Let \(C \subset \bbR \times M\) denote an embedded,
pseudoholomorphic curve. Then \(C\) has a canonical complex structure, and
\(C\)'s normal bundle has a complex structure and, as a complex line
bundle, a canonical holomorphic structure. Let \(\upN\) denote the latter
bundle. Associated to \(C\) is an operator \(\mathcal{D} _{C}: C^{\infty }(C; \upN)
\to C^{\infty }(C; \upN \otimes T^{0,1}C)\) depicted in (\ref{(2.4)}).
This operator defines a Fredholm operator
from the space of \(L^{2}_{1}\) sections along \(C\) of \(\upN\) to the space of
\(L^{2}\) sections of \(\upN \otimes T^{0,1}C\). The assumption that \(J \in \mathcal{J}
_{1\ff}\) guarantees that this Fredholm incarnation of (\ref{(2.4)}) in the case
when \(C\) is not \(\bbR\)-invariant is surjective.

 Suppose next that \(\Theta \in \mathcal{A}\) and that \((\gamma , m) \in \Theta\)
 . Associated to \(\gamma \) is the \(\bbR\)-invariant, pseudoholomorphic
 cylinder \(C = \bbR \times \gamma \subset \bbR \times M\). The vector
 field \(\varphi_{\gamma*}\partial _{z}\) along \(\bbR \times \gamma\)
 identifies the normal bundle of \(C\) with the product bundle \(C \times
 \bbC \). Meanwhile, \(\varphi_\gamma \) identifies \(C\) with \(\bbR \times
 S^1\). This understood, what is written in (\ref{(2.4)}) can be viewed as a
 differential operator on \(C^{\infty }(\bbR \times S^1; \bbC )\). In this
 guise, the sections \(\nu _{C}\) and \(\mu _{C}\) become the functions
 \(\nu_\gamma \) and \(\mu_\gamma \) that appear in ({\ref{(4.1)}). What is written in
 (\ref{(2.4)}) can be lifted to any
given \(p\) -fold cover of \(\bbR \times S^1\) and so act on
\(C^{\infty }(\bbR \times (\bbR /(2\pi p \bbZ )); \bbC )\). Of interest
are the cases when \(p\leq m\). This lifted version of (\ref{(2.4)})
defines a Fredholm operator mapping the
 space of \(L^{2}_{1}\)-functions on \(\bbR \times (\bbR /(2\pi p\bbZ
 ))\) to the space of \(L^2\)-functions. This lifted version of (\ref{(2.4)})
 is also denoted by \(\mathcal{D}_C\) when \(C=\bbR\times \gamma\). Note
 that such a lifted version of \(\mathcal{D}_C\) has trivial kernel and cokernel when \(\gamma\) is non-degenerate.
\bigbreak

\subsubsection*{\it Part 2.} Suppose that \(\Theta _{-}\) and \(\Theta _{+}\) come from
\(\mathcal{A}\). Recall from Section \ref{sec:2(a)} that the assumption that \(J \in \mathcal{J} _{1\ff}\) implies the following:
\begin{itemize}
\item  Each component of \(\mathcal{M} _{1}(\Theta _{+},
\Theta _{-})\) is a smooth, 1-dimensional manifold with a free
\(\bbR\)-action that is induced by the action of \(\bbR\) 
on \(\bbR \times M\) as the constant translations along the \(\bbR\)
 factor.
\item  Let \(\Sigma \in \mathcal{M} _{1}(\Theta _{+},
\Theta _{-})\). There is precisely one pair from \(\Sigma\)
 whose pseudoholomorphic curve component is not \(\bbR\)-invariant. 
\item  The latter pair has the form \((C, 1)\)
with \(C\) an embedded pseudoholomorphic subvariety which
has no points in common with the \(\bbR\)-invariant subvarieties
from \(\Sigma\). 
\item  The corresponding version of \(\mathcal{D}
_{C}\) has index 1 and trivial cokernel.
\end{itemize}

Granted the above, 
Quillen's ideas \cite{Q} about determinant line bundles identify the real line 
\begin{equation}\label{(4.5)}
 \det\, (\oplus _{(C,m)\in \Sigma } \ker (\mathcal{D} _{C}))
\end{equation}
as the fiber over \(\Sigma \in \mathcal{M} _{1}(\Theta _{+}, \Theta _{-})\) of
the orientation sheaf, \(\Lambda (\Theta _{+}, \Theta _{-})\), of \(\mathcal{M}
_{1}(\Theta _{+}, \Theta _{-})\). Definition 9.9 in \cite{HT1} can be used
almost verbatim in the present context to define the notion of a
{\em coherent system of orientations} for the collection \(\{\Lambda
(\Theta _{+}, \Theta_-)\}_{\Theta_-, \Theta_+\in \mathcal{A}}\). This
means the following thing:
Each \(\Theta\) from \(\mathcal{A}\) has a canonically associated
\(\bbZ /2\bbZ \)-module, which we denote by \(\Lambda (\Theta )\). Moreover, there is a canonical
isomorphism from \(\Lambda (\Theta _{+}, \Theta _{-})\) to \(\Lambda (\Theta
_{+})\Lambda (\Theta _{-})\). This understood, a set of orientations
\(\{\gro(\Theta _{+}, \Theta _{-}) \in \Lambda (\Theta _{+}, \Theta_-)\}_{\Theta_-, \Theta_+\in \mathcal{A}}\)
is said to be coherent when there exists a corresponding set of orientations
\(\{\gro_{\Theta } \in \Lambda (\Theta )\}_{\Theta \in \mathcal{A}}\) such
that any given \(\gro(\Theta _{+}, \Theta _{-})\) is equal to \(\gro(\Theta
_{+})\, \gro(\Theta _{-})\).

Replace the phrase `embedded contact homology' with `periodic Floer
homology' in Section 9.5 of \cite{HT1} and Parts 3-6 in Section 3b of \cite{T3}
to see how a coherent system of orientations is canonically defined by
assigning to each set \(\Theta \in \mathcal{A}\) an ordering of the subset of
pairs whose periodic orbit is hyperbolic and positive. Such
orientation is denoted by \(\gro _{pFh}\) in Theorem \ref{thm:4.1} below.

The orientation \(\gro _{pFh}\) is used to define the \(\pm 1\)
contributions from each component of \(\mathcal{M} _{1}(\Theta _{+}, \Theta
_{-})\) as follows: The generator of the \(\bbR\)-action on \(\mathcal{M} _{1}(\Theta
_{+}, \Theta _{-})\) orients any given component. Either this
orientation agrees with \(\gro_{pFh}\) or not. If so, assign the component
\(+1\), if not assign the component \(-1\).

\subsection{Instantons on \(\bbR \times M\)}
\label{sec: 4(b)} 
 With the class \(\Gamma \in H_{1}(M; \bbZ )\) chosen, use \(E\) in what
 follows to denote a complex line bundle over \(M\) whose first Chern
 class is Poincar\'e dual to \(\Gamma\). The \(\Spin^{\bbC }\) structure used
 in what follows is such that its spinor bundle \(\bbS\) splits as in
 (\ref{(3.14)}).

Fix \(r \geq 1\). The relevant versions of the Seiberg-Witten equations on
\(\bbR \times M\) constitute a system of equations for a pair \((A, \psi )\co
\bbR \to Conn(E) \times C^{\infty }(M; \bbS)\); these being
\begin{equation}\label{(4.6)}
\begin{cases}
& \frac{\partial }{\partial s} A+B_{A} -
r\, (\psi ^{\dag }\tau\psi - i \textsf{a}) + \frac{1}{2} B_{0}=0,\\
& \frac{\partial}{\partial s} \psi +D_A\psi =0,
\end{cases}
\end{equation}
where \(B_0\) is as in (\ref{HP-replace}).

Of interest in what follows are the instanton solutions; those whereby
\((A, \psi)|_{s}\) converges as \(s \to \pm \infty \) to respective solutions
to (\ref{(3.15)}). The rest of this subsection concerns the solutions to
(\ref{(4.6)}). The discussion is in two parts.

 \subsubsection*{\it  Part 1.} Fix an
instanton solution to (\ref{(4.6)}) and denote it by \[s \mapsto \grd
(s) = (A, \psi )|_{s}.\] 
Associated to \(\grd\) is the differential
operator on \(C^{\infty }(\bbR \times M; iT^*M\,  \oplus \, \bbS 
\oplus \, i\, \bbR
)\) that sends a section \((b, \eta , \phi )\) to the section with
respective \(iT^*M\), \(\bbS\) and \(i\bbR \) components
\begin{equation}\label{(4.7)}
\begin{cases}
 & \frac{\partial}{\partial s} b + *db -
d\phi - 2^{-1/2}r^{1/2} (\psi ^{\dag }\tau\eta + \eta ^{\dag }\tau\psi
) ,\\
& \frac{\partial}{\partial s}\eta +
D_{A}\eta + 2^{1/2}r^{1/2} (\cl\, (b)\psi + \phi\psi ) ,\\
& \frac{\partial}{\partial s} \phi+ *d*b - 2^{-1/2}r^{1/2} (\eta ^{\dag }\psi - \psi ^{\dag }\eta ) .
\end{cases}
\end{equation}

Because of Items (3a) and (3c)  of Theorem \ref{thm:3.2}, 
this operator defines a
Fredholm operator 
\[
\grD _{\grd}\co L^{2}_{1}( \bbR\times M; iT^*M \oplus \, \bbS \oplus
\, i\, \bbR  ) \to L^{2} ( \bbR\times M; iT^*M \oplus \, \bbS \oplus
\, i\, \bbR  ).\]

The instanton \(\grd\) is said
to be {\em non-degenerate} when \(\grD _{\grd}\) has trivial
cokernel.

 Fix a pair, \(\grc _{-}\) and \(\grc _{+}\), of solutions to
 (\ref{(3.15)}). Recall from Section \ref{sec:1(b)} that \(\EuScript{M}_{1}(\grc _{-}, \grc _{+})\) denotes the space
 of instanton solutions to (\ref{(4.6)}) with the two properties: The
 limit as \(s \to -\infty \) is \(\grc _{-}\) and whose limit as \(s
 \to \infty \) is gauge
 equivalent to \(\grc _{+}\). In addition, the corresponding version of
 (\ref{(4.7)}) has Fredholm index 1. Note that the space \(\EuScript{M}_1(\grc _{-}, \grc
 _{+})\) is canonically associated to the pair of respective gauge
 equivalence classes of \(\grc _{-}\) and \(\grc _{+}\). This is because there
 is a canonical identification between \(\EuScript{M}_{1}(\grc _{-}, \grc _{+})\)
 and \(\EuScript{M}_{1}(u\grc _{-}, \grc _{+})\) for any given map
 \(u: M \to S^1\), sending a given instanton \(\grd = (A, \psi )\) to the instanton
 \(u\cdot \grd = (A - u^{-1}du, u\, \psi )\).

The space \(\scrM _{1}(\grc _{-}, \grc _{+})\) has the structure of a
smooth, 1-dimensional  manifold in a neighborhood of any non-degenerate
instanton. This space also enjoys a free \(\bbR\)-action, this induced
from the action of \(\bbR\) on \(\bbR \times M\) by translation along the \(\bbR\)
factor. As the assigment \(\grd \to \grD _{\grd}\) is \(\bbR\)-equivariant, so it follows that a component in \(\EuScript{M}_{1}(\grc _{-},
\grc _{+})\) of a non-degenerate instanton consists solely of
non-degenerate instantons. In particular, such a component is a smooth
manifold with a diffeomorphism to \(\bbR\) that intertwines the
afore-mentioned \(\bbR\)-action with the translation action of \(\bbR\) on
itself.

\subsubsection*{\it Part 2.} Let \(\grc _{-}\) and \(\grc _{+}\) for the moment denote a
given pair of non-degenerate elements in \(\op{Conn}(E) \times C^{\infty }(M;
\bbS)\). Reintroduce the notation \(\grP = \grP (\grc _{-}, \grc
_{+})\) from Section \ref{sec:1(b)}, now regarded as a space of
piecewise differentiable maps from \(\bbR\) to \(\op{Conn}(E) \times C^{\infty
}(M; \bbS )\).
Each \(\grd \in \grP \) has its corresponding version of (\ref{(4.7)}). Given
that both \(\grc _{-}\) and \(\grc _{+}\) are non-degenerate, \(\grd\)'s
version of (\ref{(4.7)}) defines a Fredholm operator, \(\grD_\grd\), mapping
the space of \(L^{2}_{1}\) sections over \(\bbR \times M\) of the
bundle \(i\, T^*M
\oplus \bbS \, \oplus i\, \bbR \) to the space of \(L^{2}\) sections of this same
bundle. Quillen \cite{Q} showed (in a somewhat different context) how such
operators define a real line bundle, \(\det\, (\grD _\grd) \to \grP \). In
particular, if \(\grd \in \grP \) is such that either the kernel or
cokernel of \(\grD _{\grd}\) is non-trivial, then the fiber of \(\det\, (\grD_\grd
)\) at a given \(\grd \in\grP  \) has a canonical identification with \(\bigwedge
^{\max}(\ker\, (\grD_\grd )) \times _{\bbR } (\bigwedge
^{\max}\cok \, (\grD_\grd ))\). Introduce \(\Lambda (\grc _{-}, \grc _{+})\)
to denote the orientation sheaf of \(\det\, (\grD )\).

 Take \(r \gg  1\) so that Theorem \ref{thm:3.2} can be used to
 conclude that all solutions to (\ref{(3.15)}) are
 non-degenerate. Granted this non-degeneracy, Chapter 20 of \cite{KM}
 associates to each gauge equivalence class \(\grc \in \scrC^r\) an
 associated \(\bbZ /2\bbZ \)-module, \(\Lambda (\grc )\), with the
 following property: If \(\grc _{-}\) and \(\grc _{+}\) are any two
 solutions to (\ref{(3.15)}), then there is a canonical isomorphism
 between the modules \(\Lambda (\grc _{-}) \otimes_{\bbZ /2\bbZ }
 \Lambda (\grc _{+})\) and \(\Lambda (\grc _{-}, \grc _{+})\). In
 particular, a choice of \(\gro(\grc _{-}) \in \Lambda (\grc _{-})\)
 and \(\gro(\grc _{+}) \in \Lambda (\grc _{+})\) defines a unique element
 in \(\Lambda (\grc _{-}, \grc _{+})\). A collection of orientations
 \(\{\gro(\grc _{-}, \grc _{+}) \in \Lambda (\grc _{-}, \grc _{+})\}_{\grc _{-}, \grc _{+} \in \scrC^r}\) 
is said to be {\em coherent} if there
exists a corresponding set of orientations \(\{\gro(\grc ) \in \Lambda  (\grc )\}\)
such that any given \(\gro(\grc _{-}, \grc _{+})\) is equal to \(\gro(\grc _{-})\gro(\grc _{+})\). 

The relevance of this orientation business to Seiberg-Witten Floer
cohomology is as follows: Suppose that \(\grc _{-}\) and \(\grc _{+}\) are
solutions to (\ref{(3.15)}), and suppose that \(\grd \in \EuScript{M}_{1}(\grc _{-},
\grc _{+})\) is nondegenerate. Then the restriction of \(\Lambda (\grc
_{-}, \grc _{+})\) to the component of \(\EuScript{M}_{1}(\grc _{-}, \grc
_{+})\) containing \(\grd\) is canonically isomorphic to the latter's orientation
sheaf. With this understood, fix orientations \(\{\gro(\grc ) \in \Lambda
(\grc )\}_{\grc\in \scrC^r}\) so as to
define a collection of coherent orientations for \(\{\Lambda (\grc _{-},
\grc _{+})\}_{\grc _{-}, \grc _{+} \in \scrC^r}\).
Use these orientations to define the orientation for the components of
\(\{\EuScript{M}_{1}(\grc _{-}, \grc_{+})\}_{\grc _{-}, \grc _{+} \in \scrC^r}\)
with nondegenerate instantons. 
This orientation is denoted in what follows
by \(\gro_{SW}\). Let \(\EuScript{M} \subset \EuScript{M}_{1}(\grc _{-}, \grc _{+})\) denote
a component with nondegenerate instantons. The generator of the \(\bbR\)-action on \(\EuScript{M}\) also orients \(\EuScript{M}\). This orientation is denoted by
\(\gro_{\bbR }\). Now view \(\grc _{-}\) and \(\grc _{+}\) as generators of the
Seiberg-Witten Floer cohomology complex. Then \(\EuScript{M}\) contributes \(+1\) to
the sum that defines via (\ref{(1.9)}) the multiple of \(\grc _{-}\) in the
coboundary of \(\grc _{+}\) when \(\gro_{SW} = \gro_{\bbR }\). Otherwise,
\(\EuScript{M}\) contributes \(-1\).


\subsection{Pseudoholomorphic curves and instantons}
\label{sec:4(c)}
Assume as usual that \(\Gamma\) is monotone. Therefore, both
\(\mathcal{M}_1(\Theta_+, \Theta_-)\) and \(\scrM_1(\grc_-, \grc_+)\)
have finitely many components. Moreover,
\begin{equation}\label{mono-balanced}
[w_\ff]=\kappa \, c_1(\grs) \quad \text{for some \(\kappa\in \bbR\)}
\end{equation}
because \(\Gamma\) is monotone, and (\ref{mono-balanced}) with
(\ref{I-period}) imply the following: Given any \(\grd(s)=(A, \psi)\in
\grP(\grc_-, \grc_+)\) for a fixed pair \(\grc_-, \grc_+\), 
a bound on its spectral flow \(\scrI(\grd)\) gives rise to
a bound on the integral
\begin{equation}\label{(4.9)}
\frac{i}{2\pi}\int_{\bbR\times M} \frac{\partial A}{\partial
  s} \, ds  \wedge w_\ff .
\end{equation}
With this said, 
the following is Section \ref{sec: 4}'s main theorem.

 \begin{thm}\label{thm:4.1}  
Let \(\pmb{\mu}_P=\{(F, w_F), \ff, \Gamma, J\}\) and
\(\pmb{\mu}_S=\{M, \grs_\Gamma, \varpi_r, g, \grq\}\) be 
as in the statement of Theorem \ref{thm:1}.
Then there exists a constant \(\kappa > 1\), 
and given \(N>\kappa\), \(\delta< \kappa^{-1}\) and a \((\delta,
N)\)-approximation \((\ff ', J')\) to \((\ff, J)\), there exists
\(\kappa'\geq \kappa\) such that if \(r\geq\kappa'\), then the following is
true: Use \((\ff ', J')\) to define \(\mathcal{A}\) and \(\scrC^r\).
\begin{itemize}
\item[\rm (0)] 
Theorem \ref{thm:3.2} can be invoked so as to provide
the identification, \(\Phi ^{r}\), between \(\mathcal{A}\)  and
\(\scrC^r\). 
\item[{\rm (1)}]  Fix a pair \(\Theta _{-}\) and \(\Theta _{+} \) in \(\mathcal{A}\); and let \(\grc _{-}\) and \(\grc
_{+}\) denote solutions to (\ref{(3.15)}) whose respective gauge
equivalence classes are \(\Phi ^{r}(\Theta _{-})\) and \(\Phi
^{r}(\Theta _{+})\). Let \(\mathcal{M}_1(\Theta_+, \Theta_-)\)
and \(\scrM_1(\grc_-, \grc_+)\) be the versions defined using \((\ff', J')\).
\begin{itemize}
\item[{\rm (a)}]  The space \(\EuScript{M}_{1}(\grc _{-}, \grc
_{+})\) has a finite set of components and each component
consists of non-degenerate solutions to (\ref{(4.6)}).
\item[{\rm (b)}]
 There is a smooth \(\bbR\)-equivariant diffeomorphism \[\Psi^{r}\co \mathcal{M} _{1}(\Theta _{+}, \Theta _{-}) \to \EuScript{M}_{1}(\grc _{-}, \grc _{+}).\] 
\end{itemize}
\item[{\rm (2)}]  There are coherent systems of orientations
for sheaves \newline
\(\{\Lambda (\Theta _{+}, \Theta_-)\}_{\Theta_-, \Theta_+\in \mathcal{A}}\) and 
\( \{\Lambda (\grc _{-}, \grc_{+})\}_{\grc_-, \grc_+\in \scrC^r}\)
such that any given \((\Theta _{-}, \Theta _{+})\) version of
\(\Psi^{r}\) maps the orientation \(\gro_{pFh}\) to the
orientation \(\gro_{SW}\).
\end{itemize}
\end{thm}

The proof of this theorem is given below modulo certain facts about
solutions of (\ref{(4.6)}) that are derived in the upcoming Section
\ref{sec:5}.

\subsubsection*{\it \textbf{Proof of Theorem \ref{thm:4.1}, Part 1:}}

 The \(\bbR \)-equivariant map \[
\Psi^{r}\co \mathcal{M} _{1}(\Theta _{+}, \Theta_{-}) \to
\EuScript{M}_{1}(\grc _{-}, \grc _{+})\] 
is constructed by copying
 the construction given in Sections 4-7 of \cite{T2}. The latter sections
 of \cite{T2} construct an analogous map for the versions of
 (\ref{(3.15)}) and (\ref{(4.6)}) that take the 1-form \(\textsf{a}\) to be a contact 1-form. The fact that
 the \cite{T2} versions of (\ref{(3.15)}) and (\ref{(4.6)}) use a contact 1-form does not
 play a role in the constructions. By the same token, the arguments
 from Section 7k of \cite{T2} can be borrowed almost verbatim to prove
 that \(\Psi^{r}\) is injective when \(r\) is large.

The assertion that \(\Psi^{r}\) maps \(\mathcal{M}_{1}(\Theta _{+}, \Theta _{-})\)
onto a union of components of \(\EuScript{M}_{1}(\grc _{-}, \grc _{+})\) that
contain solely non-degenerate instantons can be had by copying the
arguments given in Section 3a of \cite{T3} that prove the analogous
assertion for when (\ref{(3.17)}) and (\ref{(4.6)}) take the 1-form
\(\textsf{a}\) to be a contact form. Assertion (2) of Theorem
\ref{thm:4.1} are
proved using almost verbatim the arguments from Section 3b of \cite{T3}
that prove the analog for the case when the 1-form \(\textsf{a}\) is a contact
1-form.

 \medbreak

\subsubsection*{\bf \textit{Proof of Theorem \ref{thm:4.1}, Part 2:}} 
It remains only to
prove that \(\Psi^{r}\) maps \(\mathcal{M}_{1}(\Theta _{+}, \Theta _{-})\) onto
\(\EuScript{M}_{1}(\grc _{-}, \grc _{+})\). The arguments for this borrow
heavily from what is done in Sections 3-7 of \cite{T4} to prove an
analogous result in the case when the 1-form a that appears in (\ref{(3.15)})
and (\ref{(4.6)}) is a contact form. As in \cite{T4}, the proof that \(\Psi^{r}\) maps
\(\mathcal{M} _{1}(\Theta _{+}, \Theta _{-})\) onto \(\EuScript{M}_{1}(\grc _{-}, \grc
_{+})\) has three distinct parts which, for lack of better terms, will
be called {\em estimation, convergence}, and
{\em perturbation}.

\subsubsection*{\it (a) The estimation part:} Here and in
\cite{T4}, the global part of the proof establishes certain properties of
instanton solutions to the respective large \(r\) versions of (\ref{(4.6)}). In
\cite{T4}, the properties in question are summarized by the various results
in Section 3 of \cite{T4} and in Lemmas 5.2 and 5.3 of \cite{T4}. Each of the
relevant results in \cite{T4} has an analog here. Corresponding results are
as follows: Lemmas 3.1, 3.6, 3.8, 3.9 and 3.11 in \cite{T4} supply
pointwise a priori bounds for instanton solutions. These lemmae
respectively correspond to the upcoming Lemmas \ref{lem:5.2},
\ref{lem:5.6}, \ref{lem:5.10}, \ref{lem:5.11}
and \ref{lem:5.12}. The pointwise bounds given by Equation (3.35) from \cite{T4}
correspond to those given in the upcoming Equation (\ref{(5.52)}). Lemma 3.10
in \cite{T4} asserts a crucial monotonicity result; its analog here is
Lemma \ref{lem:5.9}. Lemmas 5.2 and 5.3 in \cite{T4} supply crucial a priori bounds
for a certain integral; Lemma \ref{lem:5.8} below plays a similar
role.

Note however that the statements of Lemmas \ref{lem:5.3}--\ref{lem:5.10} and Lemma
\ref{lem:5.12}  below differ from their analogs in \cite{T4} 
 in the following aspect: They share the common condition
 (\ref{I-bdd}), which asks for a bound on the integral (\ref{(4.9)}),
 while in \cite{T4}, this condition is substituted by the requirement
 that the spectral flow of the instantion is bounded. 
However, as observed in the beginning of this subsection, the former
condition implies the latter when (\ref{mono-balanced}) holds. This is
true both in the 
 monotone case under discussion, and in the context of \cite{T4},
   where \(w_\ff\) is replaced by \(d\textsf{a}\).

 \subsubsection*{\it (b) The convergence part:} This part
 uses the global results to assign an element in \(\mathcal{M} _{1}(\Theta
 _{+}, \Theta _{-})\) to a given instanton in \(\EuScript{M}_{1}(\grc _{-}, \grc
 _{+})\). The arguments here borrow heavily from what is done in
 Sections 5-7 in the article \(\op{SW}\Rightarrow\op{Gr}\) from \cite{T5}. There is
 necessarily more, in order to deal with the fact that \(\bbR \times M\) is
 non-compact. In any event, this part of the argument in \cite{T4} occupies
 the latter's Section 4 with the results summarized by Proposition 4.1 in \cite{T4}. 
In Lemma \ref{lem:5.8} below, we show that the condition of 
Proposition 4.1 in \cite{T4} is met in our case, and this 
replaces Proposition 5.1 in \cite{T4}.
The rest of the convergence part of the proof here can be
 copied from Section 4 of \cite{T4} with only cosmetic and notational
 changes. This understood, no more will be said here about this part
 of the proof.

 \subsubsection*{\it (c) The perturbation part:} This part
of the proof establishes that the correspondence that is established
in the convergence part of the proof is given by the map
\(\Psi^{r}\). This is to say the following: Let \(\grd \in \EuScript{M}_{1}(\grc
_{-}, \grc _{+})\) and let \(\Sigma \in \mathcal{M} _{1}(\Theta _{+}, \Theta
_{-})\) denote its partner from part (b) above. Then \(\grd\) is the image via \(\Psi^{r}\) of
\(\Sigma\) . This part of the proof in \cite{T4} occupies the latter's Sections
6 and 7. This part of the proof also calls on the global results
mentioned above. The analogous perturbation part of the proof here can
be copied from Sections 6 and 7 of \cite{T4}, again with only cosmetic and
notational changes. As a consequence, no more will be said about this
part of the proof.

\subsection{The proof of Theorem \ref{thm:1}}\label{sec:4(d)}

 The periodic Floer homology is defined once given a pair \((\ff, J)\) with
 \(\ff\co F \to F\) an area preserving diffeomorphism and with \(J \in \mathcal{J}
 _{1\ff}\). With this pair chosen, the periodic Floer homology is defined
 for a chosen monotone class \(\Gamma \in H_{1}(M; \bbZ )\). The
 assertions in Theorem \ref{thm:1} for the cases \(d_{\Gamma } < 0\) and \(d_{\Gamma }
 = 0\) follow directly from assertions (1) and (2) in
 Theorem \ref{thm:3.2} respectively. This understood, assume henceforth that \(d_{\Gamma } > 0\).

Fix some very large integer \(N\), chosen so that \(d_{\Gamma } \ll N\). This
guarantees the following: Reintroduce the set \(\mathcal{Z}\) from
Section \ref{sec: 3(b)}, whose typical
element, \(\Theta \), consists of pairs of the form \((\gamma , m)\) with
\(\gamma\) a periodic orbit and \(m\) a positive integer. Require in addition
that no two pairs from \(\Theta\) have the same periodic orbit, and that
\(\sum_{(\gamma ,m)\in  \Theta } m\gamma \) represent the class
\(\Gamma\). Let \(\Theta \in \mathcal{Z}\) and let \((\gamma , m) \in \Theta \). Then
\(\gamma\) has period much less than \(N\). Now fix a small enough
\(\delta > 0\) to invoke Proposition \ref{prop: 2.4}. This proposition finds a \((\delta ,
N)\)-approximation \((\ff', J')\) to the original pair \((\ff, J)\). As Proposition
\ref{prop: 2.4} finds a canonical isomorphism between the respective \((\ff, J)\) and
\((\ff', J')\) versions of periodic Floer homology, it is sufficient to
suppose at the outset that the pair \((\ff, J)\) is such that the conditions
in (\ref{(4.1)}) hold. This is assumed in what follows.

For each \(\Theta \in \mathcal{Z}\), introduce as in Section \ref{sec:
  3(b)} the sets
\(\grC\Theta \) and \(\grC\Theta^*\). It is a consequence of Lemma \ref{(3.8)} that these
two sets are the same. 
Moreover, \(\grX = \bigcup
_{\Theta \in \mathcal{Z} } \grC\Theta \) has a canonical identification
with the set \(\mathcal{A}\) that generates the periodic Floer homology chain
complex. This understood, then 
all \(r \geq c_{0}\) versions of the map
\(\Phi ^{r}\) from Theorem \ref{thm:3.2} supply an 
identification between the set \(\mathcal{A}\) and the set \(\scrC^r\) of gauge
equivalence classes of solutions to (\ref{(3.15)}). Theorem
\ref{thm:3.2} (3a) asserts 
that each equivalence class in \(\scrC^r\) consists of
non-degenerate solutions to (\ref{(3.15)}). As a consequence, the set \(\scrC
^{r}\) provides a set of generators for the Seiberg-Witten Floer
cohomology cochain complex. 
Note in this regard that both \(\mathcal{A}\) and \(\scrC^r\) are
finite sets.
Recall that \(p\neq 0\) denote the divisibility of \(c _{\Gamma }= c_{1}(\grs_\Gamma)\) in
\(H^{2}(M; \bbZ )/\op{Tors}\). 
Since by assumption, \((\ff, J)\) is of the modified form described in
Sections \ref{sec:2(a)}, \ref{sec:2(b)}, we may invoke Lemma
\ref{lem:vortex-degree} which asserts that with this choice, \(\deg_{pFh}(\Theta)=\deg_{\grC\mathcal{Z}}(x)\).
Thus, by Theoren \ref{thm:3.2} (3b), the map \(\Phi^r\) reverses the
relative \(\bbZ /p\bbZ \)-grading of these generators. 
The \(\bbZ\)-linear extension of Theorem
\ref{thm:3.2}'s identification from the generators to the respective chain and
cochain complexes defines an isomorphism between these complexes.

Consider now the ramifications of Theorem \ref{thm:4.1}. Let \(\Theta _{-}\) and
\(\Theta _{+}\) denote any given pair from \(\mathcal{A}\), and let
\(\grc _{-}=\Phi^r (\Theta_-)\) and
\(\grc _{+}=\Phi^r(\Theta_+)\) denote solutions to (\ref{(3.15)}) as described in Theorem
\ref{thm:3.2}. According to assertion (1a) of Theorem \ref{thm:4.1}, each
component of \(\EuScript{M}_{1}(\grc _{-}, \grc _{+})\) is non-degenerate. As a
consequence, the collection \(\{\EuScript{M}_1(\grc _{-}, \grc
_{+})\}_{\grc_-, \grc_+ \in \scrC^r}\)
 can be
used to define the differential on the Seiberg-Witten Floer cochain
complex via the formula (\ref{(1.9)}).

Let \(\Theta _{-}\), \(\Theta _{+}\) and \(\grc _{-}\), \(\grc _{+}\)
be as the
above. Then assertion (1b) of Theorem \ref{thm:4.1} supplies an \(\bbR\)-equivariant diffeomorphism between \(\mathcal{M} _{1}(\Theta _{+}, \Theta
_{-})\) and \(\EuScript{M}_{1}(\grc _{-}, \grc _{+})\). According to
assertion (2) of Theorem 4.1, this diffeomorphism can be assumed to preserve
the respective coherent orientations. This implies that the integer
\(\sigma (\Theta _{+}, \Theta _{-})\) that is used in (\ref{(1.5)}) to define the
periodic Floer homology differential is identical to the corresponding
integer \(\sigma (\grc _{-}, \grc _{+})\) that is used in (\ref{(1.9)}) to define
the Seiberg-Witten Floer cohomology differential. It follows directly
from this that the identification between the respective generating
sets given by Theorem \ref{thm:3.2} intertwines the action on these generators
of the respective differentials. As a consequence, the isomorphism
provided by Theorem \ref{thm:3.2} between the respective chain and cochain
complexes descends to give a relative \(\bbZ /p\bbZ \)-grading reversing
isomorphism between the periodic Floer homology and the Seiberg-Witten
Floer cohomology. This is the isomorphism asserted in Theorem \ref{thm:1}.

\section{The global arguments: Properties of instantons}\label{sec:5}
\setcounter{equation}{0}

This section derives the properties of solutions to (\ref{(4.6)})
that constitute the global part of Part 2 of the proof of Theorem
\ref{thm:4.1}. As noted in the preceding subsection, everything done here has an
analog in either Section 3 or Section 5 of \cite{T4}. In addition, proofs
of analogous assertions are very similar; and in some cases the
similarity is such as to warrant directing the reader to the
corresponding \cite{T4} argument.

\subsection{Special functions on
 \(\textrm{Conn}(E)\times C^{\infty }(M;\bbS)\)}\label{sec:5(a)}

 This first subsection introduces some functions of pairs \((A, \psi )\)
 that play a central role in the proof of Theorem \ref{thm:4.1}. 
To this end,  reintroduce the the reference Hermitian connection \(A_{E}\) on \(E\) that was used in
 Section \ref{sec: 3(c)} to define the degree of a non-degenerate element in
 \(\op{Conn}(E) \times C^{\infty }(M; \bbS)\). This done, write \(A = A_{E} +
 \op{\hat{a}}_{A}\) with \(\op{\hat{a}}_{A}\in \Omega^1(M; i\bbR)\).
This done, introduce the following important functionals:
\begin{eqnarray}\label{(5.1)}
&& \op{\grc\grs} \, (A) =-\int_M \op{\hat{a}}_A\wedge \, d\op{\hat{a}}_A
-2\int_M\op{\hat{a}}_A\wedge  *\big(B_E +\frac{1}{2} B_{A_K}\big);\\
& &Q_{F}\, (A)  = i\int_M \op{\hat{a}}_A\wedge  w_\ff ; \label{(5.2)}\\
& & \gra   \, (A, \psi)= \frac{1}{2}\op{\grc\grs }- \frac{i}{2}\int_M \op{\hat{a}}_A\wedge  \varpi _r+ r\int_M \psi^\dag D_A\psi.
\label{(5.3)}
\end{eqnarray}
The first functonal above is the Chern-Simons functional, and the
third is the {\em action} functional in the version of Seiberg-Witten
Floer theory relevant to us: 
Any given pair \((A, \psi ) \in \op{Conn}(E) \times C^{\infty }(M; \bbS)\) is
a solution to (\ref{(3.15)}) if and only if \((A, \psi )\) is a critical point of
\(\gra\). Meanwhile, (\ref{(4.6)}) are the formal (downward) gradient flow equations for
\(\gra\) when the latter's gradient is defined using the \(L^{2}\) inner
product on the tangent space to \(\op{Conn}(E) \times C^{\infty }(M; \bbS)\).

 In general, none of these functionals are fully gauge-invariant. They
 vary with the cohomology class of the gauge transformation
as follows:
\begin{equation}\label{action-period}\begin{split}\op{\grc\grs }(A -u^{-1}du) -\op{\grc\grs
  }(A) &=\langle 4\pi^2c_1(\grs)\, [u]\rangle,\\
Q_F(A-u^{-1}du)-Q_F(A)& =\langle 2\pi[w_\ff], [u]\rangle\\
\gra\, (A-u^{-1}du)-\gra\, (A)& =-\langle 2\pi r[w_\ff], [u]\rangle.
\end{split}
\end{equation}
However, given two configurations \(\grc_-,
\grc_+\), and and an element \(\grd (s)\in \grP(\grc_-, \grc_+)\) from
\(u_-\cdot\grc_-\) to \(u_+\cdot\grc_+\), where \(u_-, u_+\in
C^\infty(M; \U(1))\), the value
\[
\Delta\gra (\grd):=\gra (u_+\cdot \grc_+)-\gra (u_-\cdot \grc_-)
\]
depends only on the  gauge equivalence classes of \(\grc_-, \grc_+\),
and the relative homotopy class of \(\grd(s)\). Thus, similarly to the
grading function \(\scrI  (\grc_-, \grc_+; \grh)\) introduced in Section
\ref{sec:1(b)}, we shall write
\[\Delta \gra (\grc_-, \grc_+; \grh)=\Delta\gra (\grd)\] for any
\(\grd(s)\in \grP(\grc_-, \grc_+)\) in the relative homotopy class
\(\grh\). Use the same notation convention for the functionals \(\grc\grs \)
and \(Q_F\) as well. Note that \(\Delta Q_F\) is \(2\pi\) times the
integral (\ref{(4.9)}).
The lemma that follows discusses the values of these functions for
solutions of the Seiberg-Witten equations (\ref{(3.15)}).

\begin{lemma}\label{lem:5.1} 
There exists \(\kappa > 1\) with the
following significance: Suppose that \(r \geq \kappa  \) and that
\(\grc=(A, \psi )\) is a solution to (\ref{(3.15)}) on \(M\). Then
there exists a gauge transformation \(u\in C^\infty (M; \U(1))\)
such that 
\begin{eqnarray*}
|\op{\grc\grs}\, (u\cdot \grc)| & \leq  &\kappa r^{1/2}, \\
|\, Q _F\, (u\cdot \grc)|  & \leq & \kappa  ,\\ 
|\, \gra \, (u\cdot \grc)| & \leq & \kappa r . 
\end{eqnarray*}
\end{lemma}

\subsubsection*{\it \textbf{Proof of Lemma 5.1:}} 
First, note that by (\ref{(3.18)}) and 
the first two items in Lemma \ref{lem:3.4}, we have the pointwise
bound
\[
|F_A|\leq c_0 \, r.
\]
Choose to work in a gauge
\[
d^*\op{\hat{a}}_A=0.
\]
In this gauge, we may write
\[
\op{\hat{a}}_A=a'+a_h,
\]
where \(a'\) is coexact and \(a_h\) is harmonic. By performing a
further gauge transformation, we may require that
\[
|a_h|^2\leq  b^1(M).
\]
On the other hand, the above pointwise bound on \(F_A\) implies via a
standard rescaling and elliptic regularity argument (see e.g. the
proof of Lemma \ref{lem:5.6} below) that we have
\[
|a'|\leq c_0 r^{1/2}.
\]
Thus, by Lemma \ref{lem:3.5}, with this choice of gauge
\begin{equation*}
\grc\grs \, (A)\leq c_{0} \, (|\op{\hat a}_A| + 1)\Big(\int_M |F_{A}|
+1\Big)\leq c_0' r^{1/2}.
\end{equation*}

Consider next the bound for \(|Q _{F}|\). To this end, fix \((\gamma
 ,m) \in \Theta \) and use the fact that \(H^{2}(S^1 \times D; \bbR ) = 0\)
 to find a smooth 1-form \(y_{\gamma }\) with support on the image of
 \(\varphi_\gamma \) and such that \(w _{\ff} - dy_{\gamma } = 0\) on the
 \(\varphi_\gamma \) image of the points in \(S^1 \times D\) with distance
 \(c_{0}^{-1}\) from \(S^1 \times \{0\}\). 
On the other hand, one may pointwisely estimate \(|a'|\) by multiplying
the Green's function on \(F_A-F_0\) and integrate over \(M\). Together
with Lemmas \ref{lem:3.4}, \ref{lem:3.5}, and the Seiberg-Witten
equations (\ref{(3.15)}), this shows that \(|a'|\leq c_0\) outside the
images of \(\varphi_\gamma\)'s. 

Set \[w _{\ff}' = w _{\ff} -
 \sum_{(\gamma ,m)\in \Theta } dy_{\gamma }.\]
Use this 2-form to write
\begin{equation*}
 Q _{F}\, (A) = i \int_M \hata_A \wedge w_\ff' +i\sum_{(\gamma, m)\in
   \Theta}\int_M F_A\wedge y_\gamma,
\end{equation*}
and examines the two terms on the right hand side separately.

The integrand in the first term above is supported away from the
images of \(\varphi_\gamma\)'s, where \(|\op{\hat{a}}|< c_0\) with our
choice of gauge, as we have just seen. This shows that the first term
is bounded by a constant. 

Meanwhile, the second term is also bounded by a constant, again by
combining  Lemmas \ref{lem:3.4}, \ref{lem:3.5}, and the Seiberg-Witten
equations (\ref{(3.15)}).

Lastly, to bound \(\gra\), note that by (\ref{(3.15)}), \(\int_M\psi
^\dag D_A\psi=0\). This said, the asserted bound on \(\gra\) follows
from the above bound for \(\op{\grc\grs}\) and \(Q_F\), with the
closed form \(w_\ff\) in the definition of \(Q_F\) replaced by \(r^{-1}\varpi _r\).

\subsection{The behavior of \(\alpha \), \(\beta \), and the
curvature: Part 1.}\label{sec:5(b)}

 This subsection derives various a priori bounds on the components of
 an instanton solution to (\ref{(4.6)}). 
The desired bound on the spinor are
 summarized in the following lemma. This lemma states the analog here
 of Lemma 3.1 in \cite{T4}.

\begin{lemma}\label{lem:5.2}
 There exists \(\kappa \geq 1\) with
the following significance: Suppose that \(r \geq  \kappa\)
 and that \((A, \psi = (\alpha , \beta ))\) is an
instanton solution to (\ref{(4.6)}). Then 
\begin{itemize}
\item[{\rm (1)}]  \(|\alpha| \leq 1 + \kappa r^{-1}\).
\item[{\rm (2)}]  \(|\beta|^{2} \leq \kappa r^{-1} (1 - |\alpha|^{2}) + \kappa ^{2} r^{-2}\).
\end{itemize}
\end{lemma}

The proof of this lemma uses the Bochner-Weitzenb\"ock formula for the
Dirac operator on \(\bbR \times M\). The version that follows holds for any map
\(s \mapsto (A, \psi )|_{s}\) from \(\bbR\) to \(\op{Conn}(E) \times C^{\infty }(M;
\bbS)\). The relevant Dirac operator is 
\[\scrD _{A}=\frac{\partial}{\partial s}+ D_{A};\]
and the Bochner-Weitzenb\"ock formula reads
\begin{equation}\label{(5.9)}
 \scrD _{A}^{\dag }\scrD _{A} \psi = \nabla _{\bbA}^{\dag }\nabla
_{\bbA}\psi - \cl\Big(\frac{\partial}{\partial s} A + B_{A}\Big)\psi + \frac{1}{4}
R_g \psi,
\end{equation}
where the notation used on the right hand side is as follows: First,
\(\psi\) is viewed as a section over \(\bbR \times M\) of the pull-back of the
spinor bundle \(\bbS\) . Second, \(\nabla _{\bbA}\) here denotes the
covariant derivative on this pull-back bundle that is defined by
viewing \(A\) and the canonical connection on \(K^{-1}\) as connections over
\(\bbR \times M\) for the respective pull-backs of \(E\) and \(K^{-1}\). Finaly,
\(R_g\) denotes the scalar curvature of the metric \(g\) on \(M\). The left hand
side of the preceding equation is zero when \((A, \psi )\) is an
instanton.

 \medbreak

\subsubsection*{\textit{Proof of Lemma 5.2.}} 
View \(\bbR \times M\) as a four-dimensional Riemannian manifold so as to view the instanton equations
as the 4-dimensional Seiberg-Witten equations on \(X=\bbR\times M\). 
Then, copy the arguments in
Sections 2a-c of the article \(\op{SW} \Rightarrow \op{Gr}\) from
\cite{T5}. The latter
use (\ref{(5.9)}) with zero on the left hand side and with 
\(\frac{\partial}{\partial s} A + B_{A}\)
given by (\ref{(4.6)}) to derive differential equalities for \(w = (1 -
|\alpha|^{2})\) and \(|\beta|^{2}\). In particular, these are of the
following sort:
\begin{eqnarray}\label{(5.10)}
&& \frac{1}{2}\, d^{\dag}dw + r|\alpha|^{2}w - |\nabla_A\alpha|^{2} + \gre _{w} = 0;\nonumber\\
&& \frac{1}{2} \, d^{\dag}d|\beta|^{2} + r|\alpha|^{2}|\beta|^{2} + r
\, (1 +
|\beta|^{2})\, |\beta|^{2} + |\nabla _{A}\beta|^{2} + \gre _{\beta } = 0,
\end{eqnarray}
where 
\[\begin{split}
|\gre _{w}| & \leq c_{0}\, \big(|\alpha|^{2} + |\nabla_{A}\beta|^{2} +
|\beta|^{2}\big) \quad \text{and}\\
 |\gre _{\beta }|  & \leq c_{0} \, \big(|\beta|^{2} + |\beta| |\alpha| +
 |\beta| |\nabla_{A}\alpha|\big). 
\end{split}
\]
Here, \(\nabla _{A}\) denotes the
covariant derivative on sections over \(\bbR \times M\) of the pull-backs
of \(E\) and \(E \otimes K^{-1}\); it is defined by viewing \(A\) and the
canonical connection on \(K^{-1}\) as connections on the respective
pull-backs. As in Sections 2a-c of the article \(\op{SW} \Rightarrow \op{Gr}\) from
\cite{T5}, the maximum principle is used with these inequalities to derive
the bounds asserted by the lemma. The use of the maximum principle
in this non-compact setting requires Lemma \ref{lem:3.4} to guarantee that the
bounds given by the lemma hold as \(s \to \pm \infty \).

The next lemma is the analog here of Lemma 3.2 in \cite{T4}, as it states
bounds on \(A\) and \(B_{\bbA }\) for
an instanton \(\grd = (A, \psi )\). This lemma refers to
\begin{equation}\label{(5.11)}
\begin{split}
\scrK  & = r\sup_{s\in \bbR} \op{\underline{M}}(s), \quad \text{where}\\
\op{\underline{M}}(s)& =\int_{[s,s+1]} \big| 1-|\alpha|^2\big|.
\end{split}
\end{equation}

\begin{lemma}\label{lem:5.3} 
Fix a positive real number \(L\). 
Then there exists a constant \(\kappa > 1\) depending on \(
L\), 
with the following significance: For any \(r \geq \kappa\) and any
 solution \(\grd=(A, \psi)\) to the \(r\)'s version of instantion equation
(\ref{(4.6)}) with 
\begin{equation}\label{I-bdd}
\Delta Q_F (\grd)\leq  2\pi L,
\end{equation}
one has:
\[
\Big|\frac{\partial  }{\partial s} A- B_{A}\Big|
\leq r \, (1 - |\alpha|^{2}) + \kappa \big(r^{3/4} + (\scrK r)^{1/2}\big). 
\]
\end{lemma}

This lemma is used in conjunction with the next lemma, which gives,
among other things, a preliminary bound on \(\scrK\) .

\begin{lemma}\label{lem:5.4} 
Fix a positive real number \( L\). 
Then there exists a constant \(\kappa > 1\) depending on \(L\), 
with the following significance: For any \(r \geq \kappa\) and any
 solution \(\grd=(A, \psi)\) to the \(r\)'s version of instantion equation
(\ref{(4.6)}) satisfying (\ref{I-bdd}), one has
\[
\begin{split}
\scrK  & \leq \kappa r^{1/2};\\
\int_{[s, s+1]\times M} \big(|B_A|^2+r|\nabla_A \psi|^2\big)
& \leq \kappa r^{3/2}\quad \forall s \in \bbR  
.
\end{split}
\]
\end{lemma}

The proof of Lemma
\ref{lem:5.4} in turn requires an a priori \(L^{2}\)-bound on \(\bbR \times
M\) for \(\frac{\partial }{\partial s} A\), \(\frac{\partial }{\partial s}\psi \), \(D_A\psi \) and
\begin{equation}\label{(5.12)}
\grB _{A} := B_{A} - r (\psi ^{\dag }\tau ^{k}\psi - i \textsf{a})
+\frac{1}{2} B_{0}, 
\end{equation}
where \(B_0\) is as in (\ref{HP-replace}), as follows:

\begin{lemma}\label{lem:5.5} 
Fix a positive real number \( L\). 
Then there exists a constant \(\kappa > 1\) depending on \(
L\), 
with the following significance: For any \(r \geq \kappa\) and any
 solution \(\grd=(A, \psi)\) to the \(r\)'s version of instantion equation
(\ref{(4.6)}) satisfying (\ref{I-bdd}), one has
\[
\int_{\R\times M} \Big( \Big|\frac{\partial }{\partial s}A\Big|^2+|\grB_A|^2+2r \Big|\frac{\partial }{\partial s}\psi \Big|^2+2r |D_A\psi |^2\Big)\leq  \kappa r.
\]
\end{lemma}

\subsubsection*{\textit{Proof of Lemma \ref{lem:5.5}:}} 
Fix \(s > s' \in \bbR \). The fact
that (\ref{(4.6)}) are the gradient flow equations for \(\gra\) imply that
\begin{equation}\label{(5.13)}
\begin{split}
\int_{[s,s']\times M}\Big( \Big|\frac{\partial }{\partial
  s}A\Big|^2+|\grB_A|^2+2r \Big|\frac{\partial }{\partial s}\psi \Big|^2+2r
|D_A\psi |^2\Big) & =-\Delta\gra\, (\grd |_{[s,s']}),\\
\int_{\R\times M} \Big(\Big |\frac{\partial }{\partial
  s}A\Big|^2+|\grB_A|^2+2r \Big|\frac{\partial }{\partial s}\psi \Big|^2+2r
|D_A\psi |^2\Big)& = -\Delta\gra\, (\grd).
\end{split}
\end{equation}
Combining Lemma \ref{lem:5.1} and (\ref{action-period}), 
we see that 
\[
-\Delta\gra\, (\grd)\leq r \, (\, c_0+L) \quad \text{for some constant
  \(c_0\) depending only on \(\grc_-, \grc_+\)}.
\]
The assertion of the lemma follows by substituting this into the
previous equality.

\subsubsection*{\textit{Proof of Lemma \ref{lem:5.4}:}} 
Use the first equation in (\ref{(4.6)}) with
the first two equations of Lemma \ref{lem:5.2} to see that
\begin{equation}\label{(5.14)}
i\int_{\{s\}\times M} dt\wedge *\frac{\partial}{\partial s} A +i \int_{\{s\}\times M} dt \wedge  *
(B_A+\frac{1}{2} B_{A_K}) = r\int_{\{s\}\times M} \big| 1-|\alpha|^2\big| +\gre,
 \end{equation}
where \(|\gre | \leq \frac{r}{100} \int_M \big| 1-|\alpha|^2\big| +c_{0}\). 
The second integral on the left hand side of (\ref{(5.14)}) is \(\pi\)
times the pairing between \(c_{1}(\grs)\) and the class represented
by \(F\) in \(H_{2}(M; \bbR )\). This understood, integrate both sides of
(\ref{(5.14)}) over \([s, s+1]\) and appeal to Lemma \ref{lem:5.5} to bound the
\(L^{1}\)-norm of \(\frac{\partial }{\partial s} A\) supplies 
the bound for \(\scrK\) asserted by Lemma \ref{lem:5.4}.

The integral bound on \(|B_{A}|^{2}\) asserted by Lemma \ref{lem:5.4} 
follows from the bound on \(\scrK\), Lemma \ref{lem:5.2}, the
definition of \(\grB_{A}\) and Lemma \ref{lem:5.5}'s bound on \(\int_M|\grB_A|^2\). The integral bound for \(r|\nabla
_{A}\psi|^{2}\) asserted by Lemma \ref{lem:5.4} follows using the
bound for \(\scrK\), Lemma \ref{lem:5.2},
Bochner-Weitzenb\"ock formula for the square the Dirac operator on
\(M\), and the bound on \(\int_M |D_A \psi|^2\) provided in Lemma \ref{lem:5.5}.

\subsubsection*{\textit{Proof of Lemma \ref{lem:5.3}:}} 
With (\ref{(5.10)}), the arguments
used to prove Lemma 2.5 in the article \(\op{SW}\Rightarrow \op{Gr}\) in \cite{T5} can be
borrowed in a verbatim fashion to derive \((A, \psi )\) and \(r\) independent
constants \(z_{1}\), \(z_{2}\), \(z_{3}\) and \(z_{4}\), such that the following is
true: Introduce the functions
\begin{eqnarray}
 &q_{0} &= r \, (1 + r^{-1}z_{1}) \, (1 - |\alpha|^{2}) - z_{2}r|\beta|^{2} + z_{3},
\label{(5.15)}\\
& \mathpzc{s} &=\Big|\frac{\partial}{\partial s}A -B_{A}\Big|,\nonumber\\
& q & = \max\, (\mathpzc{s} - q_{0}, 0).\nonumber
\end{eqnarray}
The function \(q\) obeys
\begin{equation}\label{(5.16)}
 d^{\dag }dq + 2r|\alpha|^{2}q \leq z_{4 }(\mathpzc{s}  + r (1 - |\alpha|^{2})).
\end{equation}
As noted in the proof of Lemma 3.2 from \cite{T4}, the preceding equation
for \(q\) implies that 
\[
d^{\dag }dq - z_{4}q \leq c_{0} \, q_{0} +
c_{0}. 
\]
Granted this, mimic what is done in the proof of Lemma 3.2 in
\cite{T4} with the Dirichlet Green's function for the operator \(d^{\dag }d
- z_{4}\) on small radius balls in \(\bbR \times M\) to conclude the
following: If \(\rho  \in (0, c_{0}^{-1})\), then

\begin{equation}\label{(5.17)}\begin{split}
 q(x)  & \leq c_{0}\rho ^{-4}\left(\int_{\dist (x, \cdot)\leq \rho}
 \Big(\Big |\frac{\partial}{\partial s} A\Big|^2+ |B_A|^2\Big)
\right)^{1/2}\\
&\qquad +
c_{0} r\int_{\dist (x, \cdot)\leq \rho} \frac{1-|\alpha|^2}{\dist (x,
  \cdot)^2}  + c_{0} .
\end{split}
\end{equation}
This uses the fact that the Green's function for \(d^{\dag }d - z_{4}\) with pole at \(x\) obeys
\begin{equation}\label{(5.18)}
0 \leq G(\cdot , x) \leq c_{0} \dist (\cdot, x)^{-2}\quad \text{ and}
\quad  |\, dG(\cdot , x)| \leq c_{0} \dist (\cdot, x)^{-3}.
\end{equation}
Use the integral bounds for \(B_{A}\) in Lemma \ref{lem:5.4} and those
for \(\frac{\partial A}{\partial s}\) from Lemma \ref{lem:5.5} to bound the first term on the right side of (\ref{(5.17)}) by
\(c_{0} r^{3/4}\). Fix \(c \in (0, \rho)\) so as to consider the contribution
to the right most integral in (\ref{(5.17)}) from the part of the integration
domain where \(\dist(x, \cdot ) > c\). This contribution is bounded by
\(c_{0} \scrK c^{-2}\). Meanwhile, the contribution from where \(\dist(x,
\cdot ) \leq c\) is bounded by \(c_{0} r c^{2}\). This understood, take \(c =
(\scrK r)^{-1/4}\) to bound the right most integral in (\ref{(5.17)}) by \(c_{0}
(\scrK r)^{1/2}\). Thus, \[
q(x) \leq c_{0}\, \big(r^{3/4} + (\scrK r)^{1/2}\big).
\]

\subsection{The behavior of \(\alpha \), \(\beta\), and the curvature: Part 2.}\label{sec:5(c)}

 The first lemma below concerns the size of the covariant derivatives
 of \(\alpha\) and \(\beta\) . It is the analog of Lemma 3.6 in \cite{T4}. This
 lemma and the subsequent discussions use \(\nabla _{A}\) to denote the
 covariant derivative of a section of a bundle over \(\bbR \times M\) as
 defined by viewing \(A\) as a connection on the pull-back of \(E\) over \(\bbR
 \times M\). In particular, \(\nabla _{A}\) has a component that
 differentiates along the \(\bbR\) factor of \(\bbR \times M\).
 
\begin{lemma}\label{lem:5.6}  
Fix a positive real number \( L\). 
Then there exists a constant \(\kappa > 1\) depending on \(
L\), 
with the following significance: For any \(r \geq \kappa\) and any
 solution \(\grd=(A, \psi)\) to the \(r\)'s version of instantion equation
(\ref{(4.6)}) satisfying (\ref{I-bdd}), one has
\begin{itemize}
\item[{\rm (1)}]  \(|\nabla _{A}\alpha|^{2} \leq \kappa r \). 
\item[{\rm (2)}]  \(|\nabla _{A}\beta|^{2} \leq \kappa \).
\end{itemize}
In addition, for each integer \(q \geq 1\), there exists
a constant \(\kappa _{q}\) that depends only on \(L\) and is such
that when \(r \geq \kappa \), then
\begin{itemize}
\item[{\rm (3)}]
\(|\nabla_{A}^{q}\alpha| +r^{1/2}|\nabla_{A}^{q}\beta| \leq \kappa _{q} r^{q/2}\).
\end{itemize}
\end{lemma}

\subsubsection*{\textit{Proof of Lemma \ref{lem:5.6}:}} 
These claims are local in nature
and are proved by rescaling the Seiberg-Witten equation as written in
Gaussian normal coordinates about any given point. Thus, the
coordinate functions \(\{x^{\upsilon }\}_{\upsilon=\{1,2,3,4\}}\) for the
Gaussian coordinate chart are written as \(x^{\upsilon } =
r^{-1/2}y^{\upsilon }\). Uniform bounds on the curvature in the rescaled
coordinates follow from Lemma \ref{lem:5.3}. With the curvature bounded in the
rescaled coordinate, uniform bounds on the covariant derivatives of
rescaled sections can be obtained using standard elliptic regularity
techniques. Undoing the rescaling gives the asserted bounds.

  The next lemma refines Lemma \ref{lem:5.3}'s bound on the curvature. It is the
 analog here of Lemma 3.7 in \cite{T4}.

\begin{lemma}\label{lem:5.7}  
Fix a positive real number \( L\). 
Then there exists a constant \(\kappa > 1\) depending on \(
L\), 
with the following significance: For any \(r \geq \kappa\) and any
 solution \(\grd=(A, \psi)\) to the \(r\)'s version of instantion equation
(\ref{(4.6)}) satisfying (\ref{I-bdd}), one has
\[
\Big|\frac{\partial}{\partial s}A - B_{A}\Big|
\leq r \, \big(1 + \kappa \scrK ^{1/2} r^{-1/2}\big) \, (1 - |\alpha|^{2})+ \kappa
.\] 
\end{lemma}

\subsubsection*{\textit{Proof of Lemma \ref{lem:5.7}:}} 
Reintroduce the function \(q\) from
the proof of Lemma \ref{lem:5.3}. 
It follows from (\ref{(5.16)}) and from Lemma \ref{lem:5.3} that
a constant \(z_{5} \geq 1\) can chosen independently of \((A, \psi )\) and \(r\) so
that the function \[q_{1} = \max \,(q - z_{5}, 0)\] obeys
\begin{equation}\label{(5.19)}
d^{\dag }dq_{1} + 2r|\alpha|^{2}q_{1} \leq c_{0} r \, (1 - |\alpha|^{2})
\end{equation}
in the weak sense.

Given (\ref{(5.10)}), what is done in the proof of Lemma 2.7 in the article
\(\op{SW}\Rightarrow \op{Gr}\) from \cite{T5} when repeated here finds constants, \(\delta_{1} > 0\) and \(\delta _{2} > 0\), both independent of \(r\) and of \((A, \psi )\),
and such that the function
\[
\op{v}_{1} = 1 - |\alpha|^{2} + r^{-1}\delta _{1} - \delta
_{2}|\beta|^{2}
\] 
has the following properties when \(r \geq c_{0}\):
\begin{equation}\label{(5.20)}
 \begin{split}
 &\op{v}_{1}\geq  r^{-1}\delta;\\
& \op{v}_{1} \geq |1 - |\alpha|^{2}|\\
&  d^{\dag }d\op{v}_{1} + 2r\, |\alpha|^{2}\op{v}_{1} \geq 0 .
\end{split}
\end{equation}
Set \[
\varepsilon = r^{-1/2}\quad \text{and} \quad \op{v}_{2} =
\op{v}_{1}^{1-\varepsilon}. 
\]
It follows from the second inequality in (\ref{(5.20)})
that the function \(\op{v}_{2}\) also
obeys \[
\op{v}_{2} \geq \big|1 - |\alpha|^{2}\big|. 
\]
This and the third inequality in
(\ref{(5.20)}) imply that
\begin{equation}\label{(5.21)}
d^{\dag }d\op{v}_{2} + 2r|\alpha|^{2}\op{v}_{2} \geq 2r^{1/2}|\alpha|^{2}(1 - |\alpha|^{2}) .
\end{equation}
Granted the latter equation, it then follows from (\ref{(5.19)}) that 
there is
a constant \(z_{6} \geq 0\) with the following properties: First, \(z_{6}\) is
independent of both \(r\) and \((A, \psi )\). To state the second, introduce
the function 
\[
\mathpzc{v} = q_{1} - z_{6}\, (r^{1/2} + ||q_{1}||_{\infty })\, 
\op{v}_{2}.\] 
Then this function obeys:
\begin{itemize}
\item[{\rm (1)}] \(\mathpzc{v}<0\) at points where \(|\alpha|^{2}
= \frac{1}{2}\);
\item[{\rm (2)}] \(d^{\dag }d\mathpzc{v} + 2r|\alpha|^{2}\mathpzc{v} \leq 0\) 
at points where \(|\alpha|^{2} \geq \frac{1}{2}\).
\end{itemize}
This last inequality with the maximum principle finds 
\[
\mathpzc{v} \leq 0 \quad \text{where
\(|\alpha|^{2} \geq \frac{1}{2}\).}
\]
This implies the assertion of Lemma \ref{lem:5.7} at points where 
\(|\alpha|^{2} \geq \frac{1}{2}\).
The assertion at points where \(|\alpha|^{2}\leq \frac{1}{2}\)
follows directly from Lemma \ref{lem:5.3}.

\subsection{Bounds for the integral of \(r  (1-|\alpha|^{2})\)}
\label{sec:5(d)}

 Lemma \ref{lem:5.4} bounds what (\ref{(5.11)}) calls \(\scrK\) by \(c_{0}r^{1/2}\). The next
 lemma asserts that \(\scrK\) has an \(r\)-independent bound.

\begin{lemma}\label{lem:5.8}
Fix a positive real number \( L\). 
Then there exists a constant \(\kappa > 1\) depending on \( 
L\), 
with the following significance: For any \(r \geq \kappa\) and any
 solution \(\grd=(A, \psi)\) to the \(r\)'s version of instantion equation
(\ref{(4.6)}) satisfying (\ref{I-bdd}), one has
\[
\scrK \leq \kappa . 
\]
\end{lemma}

 \subsubsection*{\textit{Proof of Lemma \ref{lem:5.8}:}} 
As noted in Section \ref{sec: 3(c)}, the fact
that \(J\) is tamed by the symplectic form \(ds \wedge dt + w _{\ff}\)
requires that the 1-form \(a\) that appears in (\ref{(3.13)}) 
have norm less than 2. This understood, introduce 
\[
\delta = 2 - \sup_{M}|a |. 
\]

Recall from (\ref{(5.11)}) that \(\scrK\) is defined as the supremum
of \(r\op{\underline{M}}(s)\). Without loss of generality, suppose that
this supremum is attained at a point \(s \in \bbR \). 
Otherwise, if \(\op{\underline{M}}(s)\neq  r\scrK\) \(\forall s\), then
\[
\scrK=\lim_{s\to \infty} r \op{\underline{M}}(s) \quad \text{or}\quad
\lim_{s\to -\infty} r\op{\underline{M}}(s).
\]
Since \((A(s), \psi(s))\) converges to a solutions of (\ref{(3.15)})
when \(s\to \infty\) or \(-\infty\), the assertion of the lemma then
follows from Lemma \ref{lem:3.5}.

Recall the function \(\chi\) from Section \ref{sec:1(d)}, and
define the function \(\chi _{s}\) on \(\bbR\) by
\[
\chi_s(\cdot)=\chi\left(\frac{1}{4}|s -(\cdot )|\right). 
\]
In particular, this function has compact support in the interval \([s - 2,
s + 2]\), but it is equal to 1 on \([s, s+1]\). As will be
explained momentarily, for any \(\varepsilon \in (0, 1)\), there exists a constant \(c_{\varepsilon }
\geq 1\) that is independent of \(r\),\(\delta\), and \((A, \psi )\), so that the following
holds:
\begin{equation}\label{(5.23)}
\delta\scrK \leq   2\int_{\bbR\times M}\chi_s \, \big(r |\alpha|^2(1-|\alpha|^2)-|\nabla_A\alpha|^2\big)+
c_{\varepsilon} + c_{0}\, \varepsilon\scrK   \quad \text{when \(r \geq c_{\varepsilon}\).}
\end{equation}

Note that the integral in (\ref{(5.23)}) is bounded by \(c_{0}\) via
the following argument:
Multiply both sides of the first equation in (\ref{(5.10)}) by the function \(\chi_{s}\) and integrate the
resulting equation over \(\bbR \times M\). Then the claimed bound
follows from integration by parts,
Lemmas \ref{lem:5.2} and \ref{lem:5.6}.

This understood, then any \(\varepsilon < c_{0}^{-1}\delta \) version of (\ref{(5.23)})
supplies an \(r\)-independent bound on \(\scrK\) . 

The derivation of (\ref{(5.23)})
is given next in five steps.

\subsubsection*{\it Step 1:} In this step \(\nabla_\bbA\) or \(\nabla_A\) denotes the
covariant derivative on bundles over \(\bbR\times M\) as in (\ref{(5.9)}), and
use \(\nabla_\bbA^M\) or \(\nabla_A^M\) to denote the corresponding
covariant derivative along the \(M\)-direction. 
 The identity 
(\ref{(5.9)}) has an analog for any given pair \((A, \psi )\in
\op{Conn}(E) \times C^{\infty }(M; \bbS )\), with \(\scrD _{A}\) replaced by \(D_{A}\), with
\(\nabla _{\bbA}\) replaced by \(\nabla_\bbA^M\), and
with \(\frac{\partial}{\partial s} A\) set to zero there. 
Integrate both sides of this identity over \(M\); then integrate by
parts on both sides to find 
\begin{equation}\label{(5.24)}
\begin{split}
\int_M \Big( |B_A-r (\psi^\dag \tau\psi -i \textsf{a})|^2  +2r &
|D_A\psi|^2  \Big)= \\
\int_M\Big( |B_A|^2+ r^2 |\psi^\dag \tau\psi -i \textsf{a}|^2  &-2i r
\textsf{a} \wedge * B_A +2r
\big(|\nabla_A^M\psi|^2+\frac{R_g}{4} |\psi|^2\big)\Big).
\end{split}
\end{equation}

Now take the pair \((A, \psi )\) above to be the restriction to a slice \(\{s'\}
\times M\) of the instanton solution to (\ref{(4.6)}), and regard
either side of the identity as a function on \(s'\in \bbR\). Multiply
this function by \(\chi_s(s')\) and integrate over \(s'\in \bbR\) to
see the equality of the left hand sides of Equations (\ref{(5.25)}) and (\ref{(5.26)}) below. 

Noting that 
\(B_{A} - r(\psi ^{\dag }\tau\psi - i\textsf{a}) = \grB _{A}
- \frac{1}{2} B_{0}\) (see (\ref{(5.12)})), use
Lemma \ref{lem:5.5} to see that 
\begin{equation}\label{(5.25)}\begin{split}
\int_{\bbR\times M} \chi_s & \Big( |B_A-r (\psi^\dag \tau\psi -i
\textsf{a})|^2   +2r 
|D_A\psi|^2  \Big) \\
&\leq \int_{\bbR\times M} \chi_s \big( |\grB _A|^2+2r |D_A\psi|^2\big) +c_{0} r^{1/2}.
\end{split}
\end{equation}

Write \(\textsf{a}\) as \(dt - a\) as in (\ref{(3.13)}), and 
recall that \(a\) is in the space of 1-forms that annihilate
\(\partial_t\). Denote the latter space by \(\Omega^\perp\). Use
Lemmas
\ref{lem:5.2} and \ref{lem:5.5} to find 
\begin{equation}\label{(5.26)}\begin{split}
\int_{\bbR\times M}\chi_s & \Big( |B_A|^2+ r^2 |\psi^\dag \tau\psi -i \textsf{a}|^2 -2i r
\textsf{a} \wedge * B_A +2r
\big(|\nabla_A^M\psi|^2+\frac{R_g}{4} |\psi|^2\big)\Big)\\
& >\int_{\bbR\times M} \chi_s\big(|B_A|^2+r^2 \big|1-|\alpha|^2\big|^2 +2r |\nabla_A
\alpha|^2-2r |a|\, |B^\perp|\big)-c_0r,
\end{split}
\end{equation}
where \(B^{\perp}\) denotes the orthogonal projection of \(B_{A}\) to
\(\Omega^\perp\).
Note that the derivation of (\ref{(5.26)}) uses the
fact that
\[\int_M dt \wedge i*B_{A}=4\pi ^{2} d_{\Gamma }.
\]

 \subsubsection*{\it Step 2:}
 Write 
\(\frac{\partial}{\partial s} A\) and \(B_A\) as
\begin{equation}\label{(5.27)}
\begin{split}
\frac{\partial}{\partial s} A &= -i \, (1 - \upsilon ) \big(r (1 -
|\alpha|^{2}) + \grz \big)\,  \textsf{a} + \grx  + \grX \quad \text{and}\\
B_{A}  & = -i\upsilon\, \big(r (1 - |\alpha|^{2}) + \grz \big)\,  \textsf{a} + \grx  - \grX,
\end{split}
\end{equation}
where: 
\begin{itemize}
\item[--] \(\upsilon\) is a function on \(\bbR \times M\);
\item[--] \(\grz = r|\beta|^{2} -\frac{i}{2} *(\textsf{a}  \wedge
*B_{A_K})\);
\item[--] both \(i\grX\) and \(i\grx\) are in \(\Omega^\perp\).
\end{itemize}
It is a consequence of Lemma \ref{lem:5.2} that 
\(|\grz | \leq c_{0}\). Meanwhile, the form \(\grx\) is constrained by
the first equation 
in (\ref{(4.6)}). In particular, the latter and Lemma \ref{lem:5.2} 
require that
\begin{equation}\label{(5.28)}
|\grx | \leq c_{0} \, \big(r^{1/2}\big|1 - |\alpha|^{2}\big|^{1/2} + 1\big) .
\end{equation}
Let \(w\) denote \(( 1 - |\alpha|^{2})\). The norm of \(\grX\) is constrained
using Lemmas \ref{lem:5.4} and \ref{lem:5.7} to obey
\begin{equation}\label{(5.29)}
4|\grX |^{2} + (1 - 2\upsilon )^{2} r ^{2}w^{2} \leq r^{2} \big(1 +
c_{0}\, \scrK ^{1/2}r^{-1/2}\big)^{2}w^{2} + c_{0} \, (r|w| + 1).
\end{equation}

Given that Lemma \ref{lem:5.4} finds \(\scrK  \leq c_{0}r^{1/2}\), this last inequality requires
\begin{equation}\label{(5.30)}
 \begin{split}
 |\grX |^{2} \leq  & r^{2}\upsilon\, (1 - \upsilon ) w^{2}+
c_{0} (r^{7/4} w^{2} + r|w| + 1) ;
\\
 &-c_{0} \, r^{-1/4} \leq \upsilon  \leq 1 + c_{0}\, r^{-1/4}.
\end{split}
\end{equation}

\subsubsection*{\it Step 3:} For the rest of the argument we shall
frequently use the fact that \(|w|\leq 1\) by Lemma \ref{lem:5.2} to
bound \(w^2\) by \(|w|\) without further mention. 

Let \(\varepsilon \in (0, 1)\). It follows
from (\ref{(5.24)})--(\ref{(5.26)}), (\ref{(5.28)}), (\ref{(5.30)}), and Lemma~\ref{lem:5.5} that
\begin{equation}\label{(5.31)}
\begin{split}
 \int_{\bbR\times M}  & \chi_s \left(2\varepsilon^{-1}\Big|\frac{\partial}{\partial s}
 A\Big|^2+|B_A|^2+r^2 w^2 +2r |\nabla_A\alpha|^2-r^2 |a|\,
 |w|\right)\\
 & \leq
c_{0}r \, (\varepsilon ^{-1} + \varepsilon\scrK )
\end{split}
\end{equation}
when \(r \geq c_{\varepsilon }\) with \(c_{\varepsilon } \geq 1\) a constant that
depends solely on \(\varepsilon\) .

\subsubsection*{\it Step 4:} 
Use (\ref{(5.27)}) to see that
\begin{equation}\label{(5.32)}
 \int_{\bbR\times M}\chi_s \Big|\frac{\partial}{\partial s}
 A\Big|^2=\int_{\bbR\times M}\chi_s \big((\upsilon-1)^2r^2 w^2+|\grX|^2\big)
+ \gre_{A},
\end{equation}
where \(|\gre _{A}| \leq c_{0} \, r (\varepsilon ^{-1} + \varepsilon \scrK )\). Likewise,
\begin{equation}\label{(5.33)}
 \int_{\bbR\times M}\chi_s |B_A|^2=\int_M \chi_s
 \big(\upsilon^2r^2w^2+|\grX|^2\big) \, + \gre_{B},
\end{equation}
where \(|\gre _{B}| \leq c_{0} r \, (\varepsilon ^{-1} + \varepsilon
\scrK )\). Thus, 
\begin{equation}
\begin{split}
 \int_{\bbR\times M} & \chi_s \left(2\varepsilon^{-1}\Big|\frac{\partial}{\partial s}
 A\Big|^2+|B_A|^2-r^2 w^2 \right)\\
& \geq \int_{\bbR\times M}\chi_s
 \Big((1+2\varepsilon^{-1})\upsilon^2-4\varepsilon^{-1}\upsilon
 +2\varepsilon ^{-1}-1)r^2 w^2
\Big)
+ 2\varepsilon ^{-1}\gre_{A}+\gre_B\\
& \geq
-c_{0}\, r \, (\varepsilon ^{-2} + \varepsilon\scrK )
\end{split}
\end{equation}
for all sufficiently large \(r\), since \[(1+2\varepsilon^{-1})\upsilon^2-4\varepsilon^{-1}\upsilon
 +2\varepsilon ^{-1}-1\leq 4(1-\varepsilon^{-1})(2\varepsilon
 ^{-1}+1)^{-1}<0.\]
Combine this inequality with (\ref{(5.31)}) to get
\begin{equation}\label{(5.36)}
\begin{split}
 \int_{\bbR\times M}  & \chi_s \left(2r^2 w^2 +2r |\nabla_A\alpha|^2-r^2 |a|\,
 |w|\right)\\
 & \leq
c_{0}r \, (\varepsilon ^{-2} + \varepsilon\scrK )
\end{split}
\end{equation}
when \(r \geq c_{\varepsilon }\) for different values of \(c_0\) and \(c_{\varepsilon }\).

 \subsubsection*{\it Step 5:} 
Recall the notation  \(\sup_{M}|a | = 2 - \delta \); it follows immediately from
 (\ref{(5.36)}) that
\begin{equation}\label{(5.37)}
 \delta  \int_{\bbR\times M} \chi_s r^2 |w| +2r \int_{\bbR\times M}
 \chi_s \big( rw^2 +|\nabla_A \alpha |^2 -r|w|\big)
\leq c_{0}r \, (\varepsilon ^{-2} + \varepsilon \scrK ).
\end{equation}
Use the first inequality  in Lemma \ref{lem:5.2} to see that this inequality also holds
with \(|w|\) replaced by \(w + c_{0}r^{-1}\) in the right most integrand on
the left side of this inequality. Make this replacement. The resulting
inequality directly implies (\ref{(5.23)}).

\subsubsection*{\it Remark.} By imposing additional conditions on the
choice of the almost complex structure \(J\), the proofs of Lemmas
\ref{lem:5.3}--\ref{lem:5.8} may be simplified. 

When the degree of \(\Gamma\) is larger than \(g-1\) or less than 2, 
one may set \(a=0\) (i.e. \(J\) is admissible) yet
still have the transversality and compactness results required for
the definition of the periodic Floer homology, see the arguments in
\cite{H1}. 
In this case, Lemma \ref{lem:5.8} follows directly from taking the squares of
both sides of the Seiberg-Witten equations (\ref{(4.6)}) and the
monotonicity of the action functional, without making use of Lemmas
\ref{lem:5.3}, \ref{lem:5.4}, and the somewhat involved argument in
the above proof. 

In general, \(a\) may be taken to be arbitrarily small for the
purpose of the aforementioned transversality and compactness
results. Assuming this,  the final two steps of the preceding argument
may be omitted. 
  
 \subsection{Monotonicity for the integral of \(r (1 -
|\alpha|^{2})\)}\label{sec:5(e)}

 The lemma that follows plays a key role in the proof of the
 surjectivity of  \(\Psi^{r}\). Denote by
 \(\smm\) the
 function from \((\bbR \times M) \times (0, \infty ) \to (0, \infty )\) that
 assigns to a given pair \((\rho, x)\) the number
\begin{equation}\label{(5.38)}
 \smm (x, \rho) = r\, \int_{\dist (x, \cdot) \leq \rho}
 \big|1-|\alpha|^2\big| .
\end{equation}
The upcoming lemma gives upper and lower bounds for this function. It
is the analog of Lemma 3.10 in \cite{T4} and Proposition 3.1 in the
article \(\op{SW}\Rightarrow\op{Gr}\) from \cite{T5}.

\begin{lemma}\label{lem:5.9} 
Fix a positive real numbers \( L\). 
Then there exists a constant \(\kappa > 1\) depending on \( 
L\), 
with the following significance: For any \(r \geq \kappa\) and any
 solution \(\grd=(A, \psi)\) to the \(r\)'s version of instantion equation
(\ref{(4.6)}) satisfying (\ref{I-bdd}), one has
\begin{itemize} 
\item[{\rm (1)}]  If \(\rho_{1} > \rho _{0} \in (r^{-1/2},
\kappa ^{-1})\), then for any \(x\in \bbR\times M\), \[\smm \, (x, \rho_{1}) \geq \kappa
^{-1}\rho _{1}^{2}/\rho_{0}^{2} \, \smm \, (x,\rho _{0}).\]
\item[{\rm (2)}]  Suppose that \(|\alpha (x)|
\leq  \frac{1}{2}\) at \(x\). If \(\rho \in (r^{-1/2}, \kappa ^{-1})\),
then \[
\kappa ^{-1}\rho^{2} \leq \smm \, (x, \rho) \leq \kappa\rho^{2}.\]
\end{itemize}
\end{lemma}

\subsubsection*{\textit{Proof of Lemma 5.9:}} 
To start, let \[
\omega = ds \wedge \textsf{a} + w _{\ff}. 
\]
View \(A\) as a connection on the pull-back of the bundle \(E\)
to \(\bbR \times M\) and use \(F_{A}\) to denote its curvature 2-form. The
latter is equal to \(ds \wedge \frac{\partial}{\partial s}A + *B_{A}\)
where it is understood that \(*\) here refers to the Hodge star for the
metric on \(M\). Granted this notation, it then follows using
(\ref{(4.6)}) with
Lemma \ref{lem:5.2} that
\begin{equation}\label{(5.39)}
 \smm \, (x, \rho) =\int_{\dist (x, \cdot)\leq \rho} i\, \omega\wedge F_A + \gre_{1},
\end{equation}
where \(|\gre _{1}| \leq c_{0}\rho ^{4}\). Now fix Gaussian coordinates
\((y_{1}, \cdots , y_{4})\) centered at \(x\), chosen so that
\begin{equation}\label{(5.40)}
\omega= dy_{1} \wedge  dy_{2} + dy_{3} \wedge dy_{4} + \mathcal{O}(|y|) . 
\end{equation}
Set \(\omega_{0} = dy_{1} \wedge dy_{2} + dy_{3} \wedge dy_{4}\). It then
follows from (\ref{(5.39)}) with Lemmas \ref{lem:5.2} and \ref{lem:5.8} that
\begin{equation}\label{(5.41)}
 \smm (x, \rho) = \int_{\dist (x, \cdot)\leq \rho} i\, \omega_0 \wedge F_A
+ \gre _{2},
\end{equation}
where \(|\gre _{2}| \leq c_{0} \big(\rho^{4} + \rho\, \smm \, (x, \rho)\big)\). Next
note that \(\omega _{0} = d\theta_{0}\), where
\begin{equation}\label{(5.42)}
\theta _{0} = \frac{1}{2}
\big(y_{1}dy_{2} - y_{2}dy_{1} + y_{3}dy_{4} - y_{4}dy_{3}\big).
\end{equation}
Note in particular that 
\(|\theta_{0}|= \frac{1}{2} |y| + \gre_{3}\) where \(|\gre _{3}| \leq
c_{0}|y|^{2}\). 
Write the integral in (\ref{(5.41)}) as
\begin{equation}\label{(5.43)}
\int_{\dist (x, \cdot)\leq \rho} i\theta_0 \wedge (F_A)^T,
\end{equation}
where \((F_{A})^{T}\) denotes the restricton of \(F_{A}\) to the tangent space
of the sphere of radius \(\rho\) into \(\bbR \times M\).

 To continue, note that 
\begin{equation}\label{(5.44)}
|\theta \wedge  (F_{A})^{T}| \leq \frac{1}{2}\rho r\,  (1 +
c_{0}r^{-1/2}) \, (1 + \rho)\, \big|1 - |\alpha|^{2}\big| + c_{0}\rho  .
\end{equation}

Indeed, this follows from the bound \[
|\theta_{0}|\leq \frac{1}{2}|y| +c_{0}|y|^{2}\] 
with the following two facts: First,
\[
\Big|\frac{\partial}{\partial s}A + B_{A}\Big|\leq r\, (1 - |\alpha|^{2})
+ c_{0}, 
\]
which is a consequence of (\ref{(4.6)}) and
Lemma \ref{lem:5.2}. Second, 
\[
\Big|\frac{\partial}{\partial s}A - B_{A}\Big|\leq  r \, (1 +
c_{0}r^{-1/2})\, (1 - |\alpha|^{2}) + c_{0}, 
\]
which is a consequence of Lemmas \ref{lem:5.7} and \ref{lem:5.8}.

Granted these bounds, it follows from (\ref{(5.41)}) and
(\ref{(5.43)}) that the
function \(\rho \mapsto\smm (x, \rho)\) obeys the differential inequality
\begin{equation}\label{(5.45)}
 \smm \, (x, \rho) \leq \frac{1}{2}
\rho \, \big(1 + c_{0}\rho +
c_{0}r^{-1/2}\big)\frac{\partial\smm}{\partial \rho }(x,\rho) + c_{0}\rho^{4}.
\end{equation}

Integration of (\ref{(5.45)}) gives the what is asserted in Item (1)
of Lemma \ref{lem:5.9}. The left hand inequality of Item (2) follows
from Item (1), since \[
\smm \, (x, r^{-1/2}) \geq c_{0}^{-1}r^{-1}
\quad \text{if \(|\alpha (x)| \leq \frac{1}{2}\).}
\]
Indeed
such an upper bound follows directly from the first inequality in Lemma
\ref{lem:5.6}. The right hand inequality of Item (2) follows from Item
(1) and Lemma \ref{lem:5.8}.

\subsection{The behavior of \(\alpha \), \(\beta \) and the
curvature: Part 3.}\label{sec:5(f)}

 The bounds given below refine the bounds given in Lemmas \ref{lem:5.6} and
 \ref{lem:5.7}. The first lemma below states the analog here of what is asserted
 by Propositions 2.8 and 4.4 of the article \(\op{SW}\Rightarrow\op{Gr}\) in
 \cite{T5}. It is the analog of Lemma 3.8 in \cite{T4}.

\begin{lemma}\label{lem:5.10} 
Fix a positive real number \( L\).
Then there exists a constant \(\kappa > 1\) depending on \( 
L\), with the following significance: Fix an \(r \geq \kappa\) and a
 solution \(\grd=(A, \psi)\) to the \(r\)'s version of instantion equation
(\ref{(4.6)}) satisfying (\ref{I-bdd}). 
Set \(X_*\) to denote the subset in \(\bbR \times M\) where \(1 -
 |\alpha| \geq \kappa ^{-1}\). Then 
\begin{itemize}
\item[{\rm (1)}]  \(|\nabla _{A}\alpha|^{2} + r\, |\nabla_{A}\beta|^{2}
\leq \kappa r \, (1 - |\alpha|^{2}) + \kappa ^{2}\);
\item[{\rm (2)}]  \(r \, (1 - |\alpha|^{2}) + |\nabla _{A}\alpha|^{2} +
r\, |\nabla_{A}\beta|^{2} \leq \kappa \big(r^{-1} +r e^{-\sqrt{r} \dist
  (\cdot, X_*)/\kappa}\big)\);
\item[{\rm (3)}] \( |\beta|^{2} \leq \kappa \, (r^{-2} + r^{-1}e^{-\sqrt{r} \dist
  (\cdot, X_*)/\kappa}\big)\).
\end{itemize}
\end{lemma}

As in Lemma \ref{lem:5.6}, what is written as \(\nabla _{A}\) refers to the
covariant derivative over \(\bbR \times M\) as defined by viewing \(A\) as a
connection on the pull-back of the bundle \(E\) to \(\bbR \times M\).

\subsubsection*{\textit{Proof of Lemma 5.10:}} 
The proof differs in one place
from the proof of Lemma 3.8 in \cite{T4}. Here is the
argument in our case: Use Lemma \ref{lem:5.7} with the manipulations done in Step 2 from the
proof of Proposition 4.4 in the article \(\op{SW}\Rightarrow\op{Gr}\) in \cite{T5} to
obtain an \(r\) and \((A, \psi )\) independent constant \(z_* \geq 1\)
such that \[
y= |\nabla _{A}\alpha|^{2} + r|\nabla_{A}\beta|^{2} - z_* 
\]
obeys the differential inequality
\begin{equation}\label{(5.46)}
 d^{\dag }dy + 2r|\alpha|^{2} y \leq c_{0} r \, (1 - |\alpha|^{2})^{ } y +
c_{0} y
\end{equation}
at points on \(\bbR \times M\) where \(1 - |\alpha|^{2} \leq
c_{0}^{-1}\). Meanwhile, use (\ref{(5.10)}) with Lemma \ref{lem:5.1}
to obtain an \((A, \psi)\) and \(r\) independent constant \(z_{**}\)
such that the function \[
\mathpzc{w} =(1 - |\alpha|^{2}) - z_{**}|\beta|^{2} 
\]
obeys
\begin{equation}\label{(5.47)}
-c_{0} + y \leq d^{\dag }d\mathpzc{w} + 2r|\alpha|^{2}\mathpzc{w} \leq
c_{0}\, 
(y + 1) .
\end{equation}

Fix \(c_{* } \geq 1\) and then use (\ref{(5.46)}) and the left most inequality
in (\ref{(5.47)}) to see that \[
u = \max \Big(y - c_*\big( ||y||_{\infty} + r + 1\big) \mathpzc{w} -
c_{0}c_*, 0\Big) 
\]
obeys the differential inequality
\begin{equation}\label{(5.48)}
 d^{\dag }d u +\frac{r}{64} u \leq 0
\end{equation}
on the domain \(U \subset \bbR \times M\) where \(1 - |\alpha|^{2} \leq
c_{0}^{-1}\). If \(c_* \geq c_{0}\), then \(u\) is negative on the boundary of \(U\)
and it has compact support, the latter being a consequence of Lemma
\ref{lem:3.4}. The
maximum principle demands \(u = 0\) which proves Item (1) of the 
lemma for points in \(U\). Meanwhile, Lemma \ref{lem:3.4} and
the fact that \(\mathpzc{w} > c_{0}^{-1}\) on the complement of \(U\)
imply Item (1) on the complement of \(U\).

To obtain the assertion of Item (2), use (\ref{(5.46)})
with the left hand inequality in (\ref{(5.47)}) to see that 
\[
u' = \max (y + c_{0}^{-1}r\mathpzc{w}- c_{0}, 0) 
\]
obeys 
\[
d^{\dag }du'+\frac{r}{64} u' \leq 0 \quad \text{in \(U\)}. 
\]
Keeping this in mind, let \(c_{M} > 0\) denote a constant that is
much less than the injectivity radius of \(M\). Let \(x \in X_*\) denote a
point with \(s(x) \in [s_{0} - R, s_{0} + R]\); and let \(B \subset X_*\)
denote the ball with center \(x\) and radius equal to half of the minimum
of \(c_{M}\) and \(\dist(x, X_*)\). Use \(\rho\) to denote the radius of the ball
\(B\). Let \(G(\cdot , x)\) denote the Green's function for the operator
\(d^{\dag }d + \frac{r}{64}\) with pole at \(x\). This operator obeys the bounds
\begin{equation}\label{(5.49)}
\begin{split}
0 \leq G(\cdot , x)  & \leq c_{0} \dist (\cdot , x)^{-2} e^{-\sqrt{r}
  \dist (\cdot ,x)/c_0} \quad 
\text{and}\\
|dG(\cdot , x)| &  \leq c_{0} \dist (\cdot , x)^{-3} e^{-\sqrt{r}
  \dist (\cdot ,x)/c_0} . 
\end{split}
\end{equation}

 Multiply both sides of the inequality 
\[d^{\dag }du'+\frac{r}{64} u' \leq 0 \]
by \(\chi (\dist(\cdot , x)/\rho) G(\cdot , x)\) and integrate by over
\(B\). Given that \(|u'| \leq c_{0} r\), integration by parts finds 
\[
u'|_{x} \leq c_{0} r\, e^{-\sqrt{r} \rho /c_0}.
\]
This
implies what is asserted by the second item of the lemma. The third
item follows from the second using Lemma \ref{lem:5.2}.

 The next lemma refines the bounds given by Lemma \ref{lem:5.7}
for the curvature. This lemma is the analog of Proposition 3.4 in the article
 \(\op{SW}\Rightarrow \op{Gr}\) from \cite{T5}, and it is the analog of Lemma 3.9 in \cite{T4}.

\begin{lemma}\label{lem:5.11}
Given 
\(\scrK \geq
1\), there exists \(\kappa \geq 1\) with the following
significance: Suppose that \(r \geq \kappa \), and that \((A, \psi
= (\alpha , \beta ))\) is an instanton solution to (\ref{(4.6)}). 
Fix \(s_{0} \in \bbR \) and \(R \geq
1\); and suppose that
\[
\sup_{s\in [s_0-R-3, s_0+R+3]}\underline{\smm}(s) \leq \scrK.\]  
Then
\[
\begin{split}
  \Big|\frac{\partial A}{\partial  s}+ B_{A}\Big| &
\leq r \, (1 - |\alpha|^{2}) + \kappa;\\
 \Big|\frac{\partial A}{\partial  s}- B_{A}\Big| &
\leq r \, (1 - |\alpha|^{2}) + \kappa.
\end{split}\]

 at all points in \([s_{0} - R, s_{0} + R] \times M\). 
\end{lemma}

\subsubsection*{\textit{Proof of Lemma 5.11:}} 
Copy the proof of Lemma 3.9 in
\cite{T4} using Lemmas \ref{lem:5.2}, \ref{lem:5.7}, \ref{lem:5.8} and
\ref{lem:5.9} to replace their \cite{T4} analogs.

 Arguments from \cite{T4} refer to a connections on the pull-back of \(E\) over
 \(\bbR \times M\) that are denoted by \(\hat{A}\). Each version of \(\hat{A}\) differs
 from \(A\) where \(|\alpha|\) is near to 1. Any given version is defined from
 a specifed function \(\wp : [0, \infty ) \to [0, \infty )\) which is a
 non-decreasing function that obeys \(\wp(x) = x\) for \(x\) near zero and \(\wp
 (1) = 1\). With \(\wp\) in hand, set
\begin{equation}\label{(5.50)}
\hat A = A -\frac{1}{2}
 \wp (|\alpha|^{2})|\alpha|^{-2}
(\bar{\alpha}\nabla_{A}\alpha - \alpha \nabla_{A}\bar{\alpha}) .
\end{equation}

Here, and as previously, \(\nabla _{A}\) denotes the covariant on \(\bbR
\times M\) as defined by \(A\). The curvature, \(F_{\hat A}\), of this
connection is
\begin{equation}\label{(5.51)}
F_{\hat {A}} = (1 - \wp)F_{A} - \wp ' \nabla_{A}\bar{\alpha}
\wedge \nabla _{A}\alpha .
\end{equation}

If the assumptions of Lemma \ref{lem:5.10} hold, then
\begin{equation}\label{(5.52)}
|F_{\hat A}| \leq c_{0} \big(r^{-1} +
r\, e^{-\sqrt{r}\dist (\cdot, X_*)/\kappa}\big)
\end{equation}
on the whole of \(\bbR \times M\).
\smallbreak

\subsection{Behavior near \(\bbR  \times \gamma\)
when \(\gamma\) is elliptic}\label{sec:5(g)} 

The following is the analog here of Lemma 3.11 in \cite{T4}:

\begin{lemma}\label{lem:5.12}
Fix a positive real number \( L\). 
Then there exists a constant \(\kappa > 1\) depending on \( 
L\), 
with the following significance: For any \(r \geq \kappa\) and any
 solution \(\grd=(A, \psi)\) to the \(r\)'s version of instantion equation
(\ref{(4.6)}) satisfying (\ref{I-bdd}), and for any 
elliptic periodic orbit  \(\gamma\) with tubular neighborhood map as
described by (\ref{(4.1)}), one has:
 \[
|\beta|\leq \kappa r^{-1}, \quad |\nabla _{A}\beta|\leq \kappa
r^{-1/2}
\]
at all points in \(\bbR \times M\) with distance \(\kappa^{-1}\) or
less from \(\bbR \times \gamma \).
\end{lemma}

\subsubsection*{\textit{Proof of Lemma 5.12}:} Copy the proof of Lemma 3.11 in \cite{T4}.

\section{Generalizing to the non-monotone case}\label{sec:6}
\setcounter{equation}{0}

In the non-monotone case of either the Periodic Floer homology or
Seiberg-Witten Floer cohomology, the definition of the Floer
(co)homology in Section \ref{sec:1(a)} or \ref{sec:1(b)} no longer
works, because the moduli space \(\mathcal{M}_1(\Theta_+, \Theta_-)\)
or \(\scrM_1(\grc_-, \grc_+)\) used to define the differential fails
to be compact. To make sense of the infinite sum in (\ref{(1.5)}) or
(\ref{(1.9)}), the standard solution from the Floer theory literature
is to work with Novikov coefficients. This is done
in the Seiberg-Witten context in Chapter 30 in \cite{KM}.
We shall follow the formulation of Kronheimer and Mrowka, with
notation and terminology modified so as to better accomodate both
versions of Floer (co)homologies under discussion in one common
framework. By way of a
parenthetical remark, the formulation used in \cite{T6} for embedded contact homology and Seiberg-Witten Floer cohomology with group ring coefficients can be modified so as to apply in the present context with Novikov coefficients.

\subsection{Local coefficients}\label{sec:local-sys}

We shall use the term ``local system'' in the following more
general sense: In place of topological spaces, consider data sets of
the following type.
\begin{defn}\label{t-groupoid}
A {\em fundamental triple} consists of:
\begin{itemize}
\item a category  \(\grG\);
\item a groupoid \(\underline{\grG}\) with the same objects as
  \(\grG\), and such that for every pair \(a, b\in
  \mathpzc{Ob}(\grG)\), the set of morphisms \(\mathpzc{Mor}\, (\underline{\grG};
  a,b)\) in \(\underline{\grG}\)
  from \(a\) to \(b\) has countably many elements;
\item a functor  \(\pi_\grG\) from the former to the latter, which
  identifies the objects in both categories.
\end{itemize}
\end{defn}

{\it The model example.} Take \(\grG\) to be the category of points and paths in a
topological space of countable \(\pi_1\); let \(\underline{\grG}\) be
the fundamental groupoid of this space, and \(\pi_\grG\) defined by
taking the relative homotopy class. 

\begin{defn}   
Let \(R\) be a commutative ring , and let \(\op{Mod}_R\) denote the
category of \(R\)-modules. 
Given a fundamental triple \(\pi_\grG: \grG \to \underline{\grG}\), 
we call a functor \(\Lambda\) from \(\underline{\grG}\) to \(\op{Mod}_R\) a {\em local system of
  \(R\)-modules on \(\grG\)}. 
\end{defn}
\begin{defn}\label{or-module}
Let \(R\) be a ring such that \(-1\neq 1\) in \(R\).
Let \(\iota\) denote the generator of \(\bbZ/2\bbZ\). An {\em
  orientation functor} of a category \(\grC\) is a composition
\(\scrR\circ H\), where \(H\) is a functor from \(\grC\) to the category of
\(\bbZ/2\bbZ\)-torsors, and \(\scrR\) is the functor from the latter
category to \(\op{Mod}_{R[\bbZ/2\bbZ]}\),  which associated to
a \(\bbZ/2\bbZ\)-torsor \(Z\) 
the \(R[\bbZ/2\bbZ]\)-module 
\(R[Z]/(1+\iota)R[Z]\).

An {\em orientation local system} \(\mathpzc{O}\) of \(\grG\) 
is a orientation functor \(\scrR\circ H\) which is also a local system. 
When an element \(u_a\in H(a)\) is chosen for each object \(a\)
of \(\grG\),
\(\mathpzc{O}(h)\) may be identified with an element in \(\{1, -1\}\)
for each morphism \(h\) in \(\grG\).
The fundamental triple is said to be {\em orientable} if it is possible
to make such a system of choices 
\(\{u_a\}_{a\in\mathpzc{Ob}(\grG)}\) such that for any \(a,b\in  \mathpzc{Ob}(\grG)\),
\(H(h)\, u_a=u_b\) for all morphisms \(h \in  \mathpzc{Mor}
(\underline{\grG})\) from \(a\) to \(b\). Such a choice is said to be
an {\em orientation}.
\end{defn}

\subsubsection*{Floer thoery package.} Denote by \(G_\bullet\) the
groupoid of a single object \(\bullet\), with \(G\) as its set of
morphisms. 

A Floer theory in principle comes equipped with the following package:
\begin{itemize}
\item[{\bf [F1]}] a fundamental triple \(\{\pi_\grG\co \grG\to  \underline{\grG} \}\),
\item[{\bf [F2]}] a subset of \(\mathpzc{Ob}\, (\grG)\)  consisting of
  ``regular'' elements; denote the associated subcategories of \(\grG\) and
  \(\underline{\grG}\) by \(\grG^{reg}\) and \(\underline{\grG}^{reg}\) respectively;
\item[{\bf [F3]}] a functor \(\scrS \) from \(\grG^{reg}\) to 
  \(\bbZ_\bullet\), sending all objects to \(\bullet\).
This \(\scrS\) factors as \(\underline{\scrS}\circ
  \pi_{\grG}\). 
This is to define the notion of ``spectral flow'' or ``relative \(\bbZ\)-grading'', 
\item[{\bf [F4]}] a parameter space \(\scrP\) for the ``flow
  equation'' in the theory (such as of perturbations, metrics), 
so that corresponding to each \(\mu\in \scrP\), there is a functor \(\scrE_\mu\) from
\(\grG\) to \(\bbR_\bullet\), sending all objects to
\(\bullet\). 
Each \(\scrE_\mu\) factors as \(\underline{\scrE}_\mu\circ
\pi_\grG\). These are to be understood as notions of ``action''.
\item[{\bf [F5]}] a subcategory \(\grC_\mu\) of \(\grG^{reg}\) defined
  for elements \(\mu\) in a subset \(\scrP^{gen}\subset \scrP\)
  consisting of   ``generic'' elements, such that both \(\mathpzc{Ob}(\grC_\mu)\)
and the set \(\mathpzc{Mor}_1(\grC_\mu):=\{z\, |\, z\in
  \mathpzc{Mor}(\grC_\mu), \scrS(z)=1\}\) consist of countably many
  elements. 
Intuitively, \(\grC_\mu\) represents the category of 
 ``critical points''  and  ``broken trajectories'' of the gradient
 flow of the actional functional parametrized by \(\mu\). 
 \end{itemize}
Suppose \(R\) is not of characteristic 2, but \(R\) and
\(\mathpzc{O}\) are as in Definition \ref{or-module}
above. One needs additionally:
\begin{itemize}
\item[{\bf [F6]}] 
an orientation of \(\mathpzc{O}\) for the fundamental triple \(\pi_\grG\co
\grG\to  \underline{\grG}\), which induces an orientation functor
\(\mathpzc{O}_\mu \) of \(\grC_\mu\) 
for every \(\mu\in \scrP^{gen}\).
\end{itemize}
As previously explained, by choosing an element \(u_a\in H(a)\) for
each \(a\in \mathpzc{Ob}(\grC_\mu)\), 
one may identify \(\mathpzc{O}_\mu(h)\)
with a \(\op{sign}(z)\in \{-1, 1\}\), for each morphism \(z\) in \(\grC_\mu\).
This system of choice is called a {\em coherent orientation} for
\(\grC_\mu\) .
When \(\pi_\grG\co \grG\to \underline{\grG}\) comes from a (real) manifold
\(M\) and \(\grC_\mu\) the category of critical points and broken
trajectories of a Morse function on it, \(\mathpzc{O}\) is
defined from \(\det TM\), and \(\mathpzc{O}_\mu\) is defined from
\(\{\det T_a\mathcal{D}_a\}_{a\in \mathpzc{Ob}(\grC_\mu)}\), where
\(\mathcal{D}_a\) denotes the descending manifold from the critical
point \(a\).

Note that by definition, \(\underline{\grG}\) is
groupoid; hence each \(\mathpzc{Mor} (\underline{\grG}; a, b)\) may be
identified with a group \(\Pi\) after choosing a reference element in
it. With respect to this choice, \(\underline{\scrS}\) and
\(\underline{\scrE}_\mu\) identify respectively with the homomorphisms
\[\mathpzc{s} \co \Pi \to \bbZ \quad \text{and} \quad \mathpzc{e}_\mu\co \Pi\to \bbR.\]

\subsubsection*{Floer (co)homology with local coefficients. }
Given a Floer theory package described above
and a local system, \(\Lambda \), of \(R\)-modules on \(\grG\), 
a purported {\em Floer chain complex with local coefficient
  \(\Lambda\)} is defined from \(\grC_\mu\) by the 
following chain groups and differentials: 
\begin{equation}\label{floer-cpx}
\begin{split}
C_* (\grG, \mu ; \Lambda) &= \bigoplus_{a\in \mathpzc{Ob}\,
    (\grC_\mu)} \mathpzc{O}_\mu(a)\otimes \Lambda (a),  \\
\partial  = & 
\sum _{z\in \mathpzc{Mor}_1 (\grC_\mu), } \mathpzc{O}_\mu
  (z)\otimes \Lambda(\pi _\grG \, z) ,
\end{split}
\end{equation}

Since \(\scrS(h)-\scrS(h')= \mathpzc{s} (h-h')\), where 
\(h, h'\in \mathpzc{Mor}(\underline{\grG}; a, b)\) and hence \(h-h'\in
\Pi\), the functor \(\scrS\) defines a \(\bbZ /p\bbZ\)-grading on the Floer
complex by
\[
\op{gr}\, (a, b)=\scrS(h) \mod p \quad \text{for any \(h\in
  \mathpzc{Mor}(\underline{\grG}; a, b)\)}, 
\]
where \(p\) is the divisibility of \(\mathpzc{s}\).

The Floer cochain complex is obtained by replacing the category \(\grG\) above by
\(\grG^{op}\).

Let \(\mathpzc{Mor}_1 (\grC_\mu; a, b)\) denote the set consisting of
elements in \(\mathpzc{Mor}_1 (\grC_\mu)\) which go from the object
\(a\) to the object \(b\), and let \(\mathpzc{Mor}_{1, h} (\grC_\mu;
a, b)\) be its subset consisting of elements in the ``relative homotopy
class'' \(h\in \mathpzc{Mor}(\underline{\grG}; a, b)\). Let
\(\mathpzc{Mor}_1(\underline{\grG}; a, b)=\{h \in
\mathpzc{Mor}(\underline{\grG}; a, b) \, | \, \underline{\scrS}(h)=1\}\).
The chain groups and differentials
given can be written more concretely once a coherent orientation of
\(\grC_\mu\) is chosen:
\begin{equation}\label{twisted-cpx}
\begin{split}
C_* (\grG, \mu ; \Lambda) &= \bigoplus_{a\in \mathpzc{Ob}\,
    (\grC_\mu)} \Lambda (a),  \\
\partial  \, a= & \sum_{b\, \in \mathpzc{Ob}    (\grC_\mu)}\sum_{h\in
  \mathpzc{Mor}_1(\underline{\grG}; a, b)} \, \sigma_h(a, b)\, 
 \Lambda(h) \, b,\quad \text{where}\\
\sigma_h(a, b)= & \sum _{z\in \mathpzc{Mor}_{1,h} (\grC_\mu; a,b), } \op{sign}
  (z).
\end{split}
\end{equation}
 
In order for the map \(\partial \) to be well-defined, certain bounds
on the number of elements in \(\mathpzc{Mor}_{1} (\grC_\mu; a,b)\)
are required in order to make sense of the formula for
\(\partial\) above.  These are typically provided by results of the
following form:
\begin{equation}\label{general-cpt}
\parbox{0.8\hsize}{\textit{For a generic \(\mu\), the number of elements in
  \(\mathpzc{Mor}_1^L (\grC_\mu; a, b)\) is finite for any \(L\in \bbR\),} }
\end{equation}
where \(\mathpzc{Mor}_1^L (\grC_\mu; a, b)=\{ \, z\, |\, z\in
\mathpzc{Mor}_1(\grC_\mu; a, b), \, \scrE_\mu(z) \leq L\}\). 

Since \(\scrE_\mu\) factors through \(\pi_\grG\), this guarantees the
finiteness of \(\sigma_h (a,b)\) for each triple of \(a, b,
h\). However, \(\mathpzc{Mor}_1 (\underline{\grG}; a, b)\) often has
infinitely many elements:   with a choice a reference element, it
identifies with \(\ker \mathpzc{s}\).
Nevertheless, together with the following additional condition, (\ref{general-cpt}) guarantees that
\(\sigma_h(a, b)=0\) except for only finitely many elements \(h\) in  
\(\mathpzc{Mor}_1 (\underline{\grG}; a, b)\); hence  the formula for \(\partial\) is
(\ref{twisted-cpx}) is a finite sum when \(\mathpzc{Ob}(\grC_\mu)\)
has finitely many elements:
\(\mathpzc{e}_\mu|_{\ker\mathpzc{s}}=0\), or equivalently, 
\begin{equation}\label{mono-balanced}
\mathpzc{e}_\mu=-\kappa\mathpzc{s}  \quad \text{for some \(\kappa\in \bbR\)
  as elements in \(\op{Hom}\, (\Pi, \bbR)\).}
\end{equation}
When (\ref{mono-balanced}) fails, the infinite sum in the second line
of (\ref{twisted-cpx}) may be made sense of via particular choices of
the local coefficients \(\Lambda\), as follows.

Let \(L\in \bbR\) from this point on. 
Paraphrasing Definitions 30.2.1 and 30.2.2 in \cite{KM}, we say that:
\begin{defn}
subset \(S\subset \mathpzc{Mor}(\underline{\grG}; a, b)\) is
{\em \((\mathpzc{s},\mathpzc{e}_\mu )\)-finite} if:
\begin{itemize}
\item[{\rm (i)}] for all \(L\), the intersection \(S\cap  \mathpzc{Mor}^L
(\underline{\grG}; a, b)\) is finite.
\item[{\rm (ii)}] there exists \(d\) such that \(\scrS(h)\leq d\) for
    all \(h\in S\). 
\end{itemize} 
A local system of complete, filtered topological abelian groups \(\Lambda\) on
\(\grG\) is {\em \((\mathpzc{s},\mathpzc{e}_\mu)\)-complete} if, for any \(a, b\in
\mathpzc{Ob}(\grG)\) and any \((\mathpzc{s},\mathpzc{e}_\mu)\)-finite set \(S\subset
    \mathpzc{Mor}(\underline{\grG}; a, b)\), 
\begin{itemize}
\item[{\rm (i)}] the set \(\{\Lambda (h)\,
    |\, h\in S\}\subset \op{Hom} (\Lambda(a), \Lambda(b))\) is equicontinuous;
\item[{\rm (ii)}] the homomorphisms \(\Lambda (h)\) converges to 0 as
    \(h\) runs through \(S\).
\end{itemize} 
\end{defn}
In the above, \(\op{Hom} (\Lambda(a), \Lambda(b))\) is equipped with
the compact-open topology, and the notions of complete filtered groups
and equicontinuity are as explained in pp. 601-602 in \cite{KM}.
Note that the above definitions depend only on
\((\mathpzc{s},\mathpzc{e}_\mu)\), not on the particular choices of
\(\scrS\) and \(\scrE_\mu\), and a local system is
\((\mathpzc{s},\mathpzc{e}_\mu)\)-complete if and only if it is
\((\mathpzc{s},\lambda\mathpzc{e}_\mu)\)-complete, for any
\(\lambda\in \bbR^+\).
 
As was observed in \cite{KM}, the upshot of this definition is that for any
function \(\sigma\co \mathpzc{Mor}(\underline{\grG}; a, b)\to \bbZ\)
with \((\mathpzc{s},\mathpzc{e}_\mu)\)-finite support, 
a series of the form \(\sum_h\sigma(h)\Lambda(h)\) converges to a
continuous limit, and the terms of this series may be rearranged.
Apply this observation to the series \(\sum_{h\in
  \mathpzc{Mor}_1 (\underline{\grG}; a, b)} \sigma_h(a, b)\). 
The statement (\ref{general-cpt}) implies that the function \(h\mapsto
\sigma_h(a, b)\) is supported on a \((\mathpzc{s},\mathpzc{e}_\mu)\)-finite set, and hence
the convergence of this series. In sum, the formula
(\ref{twisted-cpx}) gives a valid definition of Floer homology of
local coefficient \(\Lambda\), once the following two conditions are
met:
\BTitem\label{Floer-condition}
\item There are finitely many elements in \(\grC_\mu\);
\item The statement (\ref{general-cpt}) holds. 
\ETitem

{\it Example.} In the case of (\ref{mono-balanced}), any (local)
coefficient is \((\mathpzc{s},\mathpzc{e}_\mu)\)-complete. In particular, Floer homology
with coefficient \(\bbZ\) is well-defined. The monotone cases
discussed in Sections \ref{sec:1(a)} and \ref{sec:1(b)}
belong to this case with
\(\kappa\neq 0\). (The \(\kappa=0\) case is the ``balanced'' case in the
terminology of \cite{KM}).

\subsection{The non-monotone case}

As explained in the previous subsection, the periodic Floer homology
and Seiberg-Witten cohomology can be defined in general with
appropriate choice of local coefficients. 

\subsubsection*{General Periodic Floer homology.}
Fix a volume form \(w_F\) on \(F\).
In the terminology of the previous subsection, this Floer theory
associates the following to each symplectic isotopy class of \(\ff\) and \(\Gamma\in H_1(M;\bbZ)\):
\begin{itemize}
\item one may take \(\grG\) to be the category with
integral 1-cycles in \(M\) in the class of \(\Gamma\) as objects, and with \(\mathpzc{Mor}(\grG;
\Theta, \Theta')\) consisting of integral 2-cycles
\(\Sigma\) in \(\bbR\times M\), which limit to \(\Theta\) or
\(\Theta'\) respectively as \(s\to
\infty\) or \(s\to -\infty\). Let \( \underline{\grG}\) be the
category with \(\mathpzc{Mor}\, (\underline{\grG};
\Theta, \Theta')=H_2 (M; \Theta, \Theta')\), and let 
\(\pi_\grG\) be the functor which sends an integral 1-cycle
\(\Theta\subset M\) to itself, and an integral 2-cycle \(\Sigma\subset \bbR\times M\)
described above to the relative homology class of its projection to \(M\).
\item The functor \(\scrS\) in this context is given by the relative grading
\(I(\Theta_+, \Theta_-; Z)\) described in Section \ref{sec:1(a)}. 
\item The parameter space \(\scrP\) may be taken to be the space of 
  pairs \((\ff, J)\), where \(\ff\) is a symplectomorphism in a
 a fixed Hamiltonian isotopy class, and \(J\in
  \mathcal{J}_f\). Using the same notation to denote the pullback of
  \(w_\ff\) under the projection \(\bbR\times M\to   M\), set
\[
\scrE_\mu(\Sigma)=\int_\Sigma  w_\ff, \quad \text{for \(\Sigma\in \mathpzc{Mor}(\grG)\).}
\]
\item We call a \(\mu=(\ff, J)\in \scrP\) ``generic'' in this context
  if \(\ff\) is nondegenerate, and \(J\in \mathcal{J}_{1\ff}\).
Let \(\grC_\mu\) be the category with
  \(\mathpzc{Ob}\, (\grC_\mu)=\mathcal{A}\) and with
  \(\mathpzc{Mor}(\grC_\mu)\) consisting of what is termed ``broken GFL''s in
  Definition 9.3 of \cite{H1}. By Lemma 9.3 in the same reference
and its generalization in Section 9.5, the set
\(\mathpzc{Mor}_1(\grC_\mu ; \Theta_+, \Theta_-)\) is the set
\(\mathcal{M}_1(\Theta_+, \Theta_-)/\bbR\) previously introduced in Section \ref{sec:1(a)}. 
\item It is explained in Section 9 of \cite{HT1} that a coherent
  orientation is provided by 
a choice of ordering of the positive hyperbolic periodic orbits in
  \(\Theta\) for each \(\Theta\in \mathcal{A}\).
\end{itemize}

In fact, the definitions of \(\grG\) and \(\scrP\) above may be chosen somewhat
differently, since they do not directly enter into the definition of
the Floer complex. Note also that in this Floer theory package,
\begin{equation}\label{pi-p}
\Pi=H_2(M;\bbZ), \quad \mathpzc{s}=c_\Gamma, \quad \mathpzc{e}_\mu=[w_\ff],
\end{equation}
regarding the latter two as homomorphisms from \(H_2(M;\bbZ)\) to
\(\bbZ\) and \(\bbR\) respectively. The condition
(\ref{mono-balanced}) becomes the defnition for monotone classes
in Section \ref{sec:1(a)}: The cases \(\kappa>0\) and \(\kappa
<0\) correspond to the positive monotone and negative monotone cases
respectively, and \(\kappa\neq 0\) in this context.

With \(\grC_\mu\) and \(\scrE_\mu\) thus specified, the validity of the
statement (\ref{general-cpt}) follows from the typical Gromov
compactness arguments as those which appeared in Lemma 9.8 of
\cite{H1}. It is also noted in Section \ref{sec:1(a)} that
\(\mathcal{A}\) consists of finitely many elements.
Thus, as both conditions in (\ref{Floer-condition}) are
met, given a \((c_\Gamma, [w_\ff])\)-complete  local coefficient
\(\Lambda\) on \(\grG\), the formula (\ref{floer-cpx}) gives a
well-defined chain complex. We call the associated homology the {\em
  periodic Floer homology with coefficient \(\Lambda\)} of the
parameter set \(\boldsymbol{\mu}_P =\{(F, w_F), \ff, \Gamma, J\}\), and denote it
  by  \[HP_* \big(\ff\co(F, w_F) \circlearrowleft   , \Gamma;
  \Lambda\big)_J,\] or simply \(HP_*(\ff, \Gamma; \Lambda)\) when
  there is no danger of ambiguity.  Under the monotonicity assumption,
  \(\Lambda\) can be taken to be \(\bbZ\); in this case the above
  definition of periodic Floer homology reduces to the version in
  Section \ref{sec:1(a)}.

\subsubsection*{General Seiberg-Witten cohomology.}
 In this Floer theory, choose
\begin{itemize}
\item \(\grG\) such that  its objects consist of gauge equivalence
  classes of configurations on \(M\), and such that \(\mathpzc{Mor}(\grG;
  \grc_-, \grc_+)=\scrB(\grc_-, \grc_+)\). Let \(\underline{\grG}\) be the category
  with the same objects, and with \(\mathpzc{Mor}(\underline{\grG};
  \grc_-, \grc_+)=\pi_0(\scrB(\grc_-, \grc_+))\). Let \(\pi_\grG\) be defined by
  taking the relative homotopy classes of elements in \(\scrB(\grc_-,
  \grc_+)\). 
\item Define \(\scrS\) using the relative grading \(\grI (\grc_-, \grc_+;
  \grh)\) introduced in Section \ref{sec:1(b)}.
\item Take the parameter space \(\scrP\) to be the set of triples
  \(\mu=(\varpi, g, (\grT, \grS))\).
Assume through out this article that \([\varpi]\neq 2\pi c_1(\grs)\) in
  \(H^2(M; \bbR)\).

Choose \(\scrE_\mu\) such that
\(\scrE_\mu(\gamma)=\frac{1}{2}\scrE^{top}_{\omega, \grq}(\gamma)\) in
the notation of \cite{KM}, where \(\omega=-i\varpi/4\), \(\grq=(\grT,
\grS)\).
 \item When \(\mu\) is generic in the sense explained in Section
  \ref{sec:1(b)}, let \(\grC_\mu\) be the category with
  \(\mathpzc{Ob}(\grC_\mu)=\scrC\), and with \(\mathpzc{Mor}(\grC_\mu; \grc_-,
    \grc_+)=\breve{M}^+(\grc_-, \grc_+)\) in the notation of \cite{KM}. In
    particular, \(\mathpzc{Mor}_1(\grC_\mu; \grc_-,
    \grc_+)\) is what is denoted by \(\scrM_1(\grc_-, \grc_+)/\bbR\)
    in Section \ref{sec:1(a)}.
\item The functor \(\mathpzc{O}_\mu\) assigns to each
  \(\grc\in \scrC^r\) the module in (22.1) of \cite{KM} (the two element
  set \(\Lambda\) in \cite{KM} is \(H\) in our notation), and to each
  \(\gamma\in \scrM_1(\grc_-, \grc_+)/\bbR\), the isomorphism
  \(\epsilon(\gamma)\) in (22.4) of \cite{KM}.
\end{itemize}
By the above  choice of \(\scrS\), \(\scrE_\mu\) and the discussion in
Section \ref{sec:1(a)}, in this Floer theory 
\begin{equation}\label{pi-s}
\Pi=H^1(M; \bbZ), \quad \mathpzc{s}=c_1(\grs), \quad
\mathpzc{e}_\mu=\mathpzc{p}_{\grs, \varpi}.
\end{equation}
The condition (\ref{mono-balanced}) now corresponds to
the positive monotone, negative monotone, and
balanced cases discussed in Chapter 29 of \cite{KM} when \(\kappa<0\),
\(\kappa>0\), and \(\kappa=0\) respectively.

It is shown in \cite{KM} that the conditions (\ref{irreducible}) and 
(\ref{Floer-condition}) hold
under our assumption that \(\mathpzc{p}_{\grs, \varpi}\neq 0\), and
thus for any \((c_1(\grs), \mathpzc{p}_{\grs, \varpi})\)-complete
local system \(\Lambda\), the construction in the previous subsection defines a Floer
homology with coefficient \(\Lambda\) associated to each parameter set 
\(\boldsymbol{\mu}_S=\{M, \grs, \varpi, g, \grq=(\grT,\grS)\}\). 
As previously observed, by replacing \(\grG\) with \(\grG^{op}\), the
same construction defines a Floer {\em co}homology. Thus,
for any \((c_1(\grs), -\mathpzc{p}_{\grs, \varpi})\)-complete local system
\(\Lambda\), there is a well-defined Floer cochain complex:
\[
\begin{split}
C^* (\grG, \mu ; \Lambda) &= \bigoplus_{a\in \mathpzc{Ob}\, (\grC_\mu)} \mathpzc{O}_\mu(a)\otimes \Lambda (a),  \\
\delta= & 
\sum _{z\in \mathpzc{Mor}_1 (\grC_\mu) } \mathpzc{O}_\mu
  (z^{-1})\otimes \Lambda(\pi _\grG \, z^{-1}) .
\end{split}
\] 
The resulting
cohomology is the {\em Seiberg-Witten Floer cohomology with coefficient
  \(\Lambda\)} for \(\boldsymbol{\mu}_S\). This  Floer cohomology only
depends on the data set (\ref{parameter-SW}).  We denote it by 
\(\textsc{Hm}^*(M, \grs  ,\mathpzc{p}_{\grs, \varpi}; \Lambda)\). To
compare with the notation in \cite{KM}, 
\[
HM^*(M, \grs  ,-\pi[\varpi]; \Lambda)=\textsc{Hm}^*(M, \grs  ,\mathpzc{p}_{\grs, \varpi}; \Lambda),
\]
and what we call a \((c_1(\grs), -\mathpzc{p}_{\grs, \varpi})\)-complete local
system is in the terminology of \cite{KM}, a \(c\)-dual complete local
system with \(c=\pi[\varpi]\).
When \(\grs\), \(\varpi\) are monotone with respect to each other,
the coefficient \(\Lambda\) will often be dropped from the notation
when it is set to be \(\bbZ\). 
\medbreak

\subsection{The general isomorphism theorem}
\label{sec:6(c)}
\par

We now generalize Theorem \ref{thm:1} to include the non-monotone
cases. Associate to each periodic Floer homology data set
\(\pmb{\mu}_P\) a Seiberg-Witten-Floer cohomology data set
\(\pmb{\mu}_S\), as described in Section \ref{sec:1(c)}.

Fix \(L \geq 1\). Given a pair, \(\Theta _{-}\) and \(\Theta _{+}\), from \(\mathcal{A}\), introduce
\(\mathcal{M} _{1L}(\Theta _{+}, \Theta _{-}) \subset \mathcal{M} _{1}(\Theta _{+},
\Theta _{-})\) to denote the elements \(\Sigma\) such that
\begin{equation}\label{(4.8)}
\int_C w_\ff \leq L,
\end{equation}
where \(C\) is the sole submanifold from \(\Sigma\) that is not an \(\bbR\)-invariant cylinder. The space \(\mathcal{M} _{1L}(\Theta _{+}, \Theta _{-})\)
has a finite set of components. This is asserted by Theorem 1.8 in
\cite{H1}. A value for \(L\) in \([1, \infty )\) is said to be
{\em generic} when there are no cases where (\ref{(4.8)}) is an
equality. The set of generic values is necessarily the complement of a
countable set.

As will be seen momentarily, the space \(\mathcal{M} _{1L}(\Theta _{+}, \Theta _{-})\)
has a precise analog on the Seiberg-Witten side of the ledger. To
define this analog, fix \(r \geq \)1 such that each solution to (\ref{(3.15)}) is
nondegenerate. Suppose that \(\grc _{-}\) and \(\grc _{+}\) are elements in
\(\scrC^r\). Given \(L \geq 1\), use \(\EuScript{M}_{1L}(\grc _{-}, \grc _{+})
\subset \EuScript{M}_{1}(\grc _{-}, \grc _{+})\) to denote the subspace of
instantons \(\grd = (A, \psi )\) that obey (\ref{I-bdd}), namely
\begin{equation}\label{I-bdd2}
\frac{i}{2\pi}\int_{\bbR\times M} \frac{\partial A}{\partial
  s} \, ds  \wedge w_\ff \leq L .
\end{equation}
Note in this regard that the above integral 
is absolutely
convergent given that \(\grc _{-}\) and \(\grc _{+}\) are both
non-degenerate.

Recall from Section 3 that Theorem \ref{thm:3.2} makes no
assumption on monotonicity. Thus, it constructs an isomorphism
\(\Phi^r\) from
\(\mathcal{A}\) to \(\scrC^r\). Keep in mind that both consist of finitely many
nondegenerate elements. 
Let \(\Theta_-\) and \(\Theta_+\) denote two elements in \(\mathcal{A}\), and let
\(\grc_-\) and \(\grc_+\) denote solutions to (\ref{(3.13)}) that give
\(\Phi^r(\Theta_-)\) and \(\Phi^r(\Theta_+)\). Theorem \ref{thm:4.1} as stated 
does require the assumption of monotonicity to guarantee that both
\(\mathcal{M}_1(\Theta_+, \Theta_-)\) and \(\scrM_1(\grc_-, \grc_+)\)
are finite. However, the theorem that follows is the non-monotone
analog.
\begin{thm}\label{thm:4.1a}  
Let \(\pmb{\mu}_P=\{(F, w_F), \ff, \Gamma, J\}\) and
\(\pmb{\mu}_S=\{M, \grs_\Gamma, \varpi_r, g, \grq\}\) be 
as in the statement of Theorem \ref{thm:1}. Fix a generic \(L\geq 1\).
Then there exists a constant \(\kappa _L> 1\) such that 
given \(N>\kappa_L\), \(\delta< \kappa^{-1}_L\) and a \((\delta,
N)\)-approximation \((\ff ', J')\) to \((\ff, J)\), there exists
\(\kappa'_L\geq \kappa_L\) so that for any \(r\geq\kappa_L'\), the
following hold for \(\mathcal{A}\) and \(\scrC^r\) defined using \((\ff ', J')\):
\begin{itemize}
\item[\rm (0)] Theorem \ref{thm:3.2} can be invoked so as to provide
the identification, \(\Phi ^{r}\), between \(\mathcal{A}\)  and
\(\scrC^r\). 
\item[{\rm (1)}]  Fix a pair \(\Theta _{-}\) and \(\Theta _{+} \) in \(\mathcal{A}\); and let \(\grc _{-}\) and \(\grc
_{+}\) denote solutions to (\ref{(3.15)}) whose respective gauge
equivalence classes are \(\Phi ^{r}(\Theta _{-})\) and \(\Phi
^{r}(\Theta _{+})\). Let \(\mathcal{M}_{1L}(\Theta_+, \Theta_-)\)
and \(\scrM_{1L}(\grc_-, \grc_+)\) be the versions defined using \((\ff', J')\).
\begin{itemize}
\item[{\rm (a)}]  The space \(\EuScript{M}_{1L}(\grc _{-}, \grc
_{+})\) has a finite set of components and each component
consists of non-degenerate solutions to (\ref{(4.6)}).
\item[{\rm (b)}]
 There is a smooth \(\bbR\)-equivariant diffeomorphism \[\Psi^{r}\co \mathcal{M} _{1L}(\Theta _{+}, \Theta _{-}) \to \EuScript{M}_{1L}(\grc _{-}, \grc _{+}).\] 
\end{itemize}
\item[{\rm (2)}]  There are coherent systems of orientations
for sheaves \newline
\(\{\Lambda (\Theta _{+}, \Theta_-)\}_{\Theta_-, \Theta_+\in \mathcal{A}}\) and 
\( \{\Lambda (\grc _{-}, \grc_{+})\}_{\grc_-, \grc_+\in \scrC^r}\)
such that any given \((\Theta _{-}, \Theta _{+})\) version of
\(\Psi^{r}\) maps the orientation \(\gro_{pFh}\) to the
orientation \(\gro_{SW}\).
\end{itemize}
\end{thm}

\subsubsection*{\bf\textit{Proof of Theorem \ref{thm:4.1a}:}} The arguments used to prove Theorem
\ref{thm:4.1} can be employed with only cosmetic changes to prove the
following: Fix a generic \(L\geq 1\). Then there exists an
\(L\)-dependent constant \(\kappa_L\) such that if \(r\geq \kappa _L\), the conclusions of assertions (1) and (2) hold with \(\mathcal{M}_{1}(\Theta_+, \Theta_-)\)
and \(\scrM_{1}(\grc_-, \grc_+)\) replaced respectively by \(\mathcal{M}_{1L}(\Theta_+, \Theta_-)\)
and \(\scrM_{1L}(\grc_-, \grc_+)\). Indeed, the only change to Part 1
is the change from \(\mathcal{M}_1(\cdot)\) to
\(\mathcal{M}_{1L}(\cdot)\). The gluing construction used in the proof
shows that \(\Psi^r\) maps \(\mathcal{M}_{1L}(\Theta_+, \Theta_-)\)
injectively into \(\scrM_{1L}(\grc_-, \grc_+)\) when \(r\) is greater than a purely \(L\)-dependent
constant. The proof that Theorem \ref{thm:4.1}'s version of
\(\Psi^r\) maps \(\mathcal{M}_{1L}(\Theta_+, \Theta_-)\) surjectively 
onto \(\scrM_{1L}(\grc_-, \grc_+)\) in the monotone case uses the monotonicity assumption only to obtain a
bound on (\ref{(4.9)}) from a bound on the index. The bound on (\ref{(4.9)}) is then used to derive the global results in Section 5. In the case at hand, the bound on (\ref{(4.9)}) is assumed because we are restricting attention to instantons in \(\scrM_{1L}(\grc_-, \grc_+)\). This understood, the arguments in Section 5 can be used verbatim to prove that \(\Psi^r\) maps \(\mathcal{M}_{1L}(\Theta_+, \Theta_-)\)
surjectively onto \(\scrM_{1L}(\grc_-, \grc_+)\). 
\bigbreak

Noting that the conditions (\ref{I-bdd}) and (\ref{I-bdd2}) correspond
under \(\Psi^r\), straight-forward modifications of the direct limit arguments 
in Sections 4d-4i of \cite{T1} then yield the following:
\begin{thm}\label{thm:2}
Let \(\boldsymbol{\mu}_P=\{(F,
w_F), \ff, \Gamma, J\}\) be a periodic Floer parameter set, where
\(\Gamma\in H_1(M; \bbZ)\) is an arbitrary class. Associate to it a
Seiberg-Witten parameter set \(\boldsymbol{\mu}_S=\{M, \grs_\Gamma, \varpi_r, g, \grq\}\) according
to the recipe described in Section \ref{sec:1(c)}. 
Let \(\Lambda_P\) be a \((c_\Gamma, [w_\ff])\)-complete local system for
the periodic Floer homology, and let \(\Lambda_S\) denote a
corresponding \((c_1(\grs_\Gamma), 2r\pi[w_\ff])\)-complete local
system for the Seiberg-Witten cohomology, in the sense specified below. 
Then 
there is an isomorphism between the two Floer (co)homologies
with coefficient \(\Lambda\)
\[
HP_*(\ff\co (F, w_F)\circlearrowleft, \Gamma ; \Lambda_P)_J\simeq
\textsc{Hm}^{-*}(M, \grs  _\Gamma, -2r\pi[w_\ff]; \Lambda_S),
\]
which reverses the relative \(\bbZ/p\bbZ\)-gradings. 
\end{thm}

Use the subscripts \(P\) and \(S\) respectively to label the
ingredients \(\pi_\grG: \grG\to \underline{\grG}\) etc. in the
periodic Floer theory and Seiberg-Witten-Floer theory, as described in
the previous subsection. 
To introduce \(\Lambda_S\), note that only its restriction to
\(\grC_{\mu_S}\) is relevant in the statement of the preceding
theorem. A \((c_\Gamma, [w_\ff])\)-complete local system \(\Lambda_P\)
on \(\grG_P\) specifies a local system on this restriction via the
maps \(\Phi^r\), \(\Psi^r\).
In particular, recall from Section \ref{sec:1(c)} 
that for a periodic Floer parameter set \(\boldsymbol{\mu}_P\)
 and its corresponding Seiberg-Witten parameter set
 \(\boldsymbol{\mu}_S\), one has
\(c_\Gamma=c_{\grs_\Gamma}\) and \(\mathpzc{p}_{\grs, \varpi_r}=-2r\pi[w_\ff]\). Moreover, comparing (\ref{pi-p}) and
 (\ref{pi-s}), and we see that the two Floer theories have identical
 \(\mathpzc{s}\), while \(\mathpzc{e}_\mu\) on the Seiberg-Witten side
 is \(-2r\pi\) times that on the periodic Floer homology side. Thus, 
 \(\Lambda_P\) is \((\mathpzc{s},
 \mathpzc{e}_\mu)\)-complete on the periodic Floer homology side if
 and only if \(\Lambda_S\) is \((\mathpzc{s},
 -\mathpzc{e}_\mu)\)-complete on the Seiberg-Witten-Floer theory side.
 In other words, the Seiberg-Witten Floer
 cohomology for \(\boldsymbol{\mu}_S\) with coefficient \(\Lambda_S\) is
 well-defined exactly when the periodic Floer homology for
 \(\boldsymbol{\mu}_P\) with coefficient \(\Lambda_P\) is.

\bigbreak

Corollary \ref{cor:PFH-inv} generalizes correspondingly: 

\begin{cor}\label{cor:PFH-inv2}
The general periodic Floer homology in the previous theorem depends
only on the Hamiltonian isotopy class of \(\ff\), the homology class
\(\Gamma\), and the local system \(\Lambda\). In fact, 
the periodic Floer homology  with ``maximally-twisted cofficients''
depends only on the symplectic isotopy class of \(\ff\) and
\(\Gamma\), in the following sense: 

Let \(\Lambda_{[w_\ff]}\supset \bbZ[H_1(M;\bbZ)]\)
denote the Novikov ring consisting of 
\(\bbZ\)-valued functions on \(H_1(M;\bbZ)\) that are finitely supported on
\([w_\ff]^{-1}(-\infty, \kappa]\) for any \(\kappa\in
\bbR\), regarding \([w_\ff]\) as an element in \(\op{Hom}(
H_2(M;\bbZ), \bbR)\). Then there is a \(\bbZ[H_1(M;\bbZ)]\)-module
with relative \(\bbZ/p\bbZ\)-grading, denoted \(HP^{[\ff]}_*(\Gamma)\),
such that 
\[
HP_* (\ff, \Gamma; \Lambda_{[w_\ff]})\simeq
HP^{[\ff]}_*(\Gamma)\otimes _{\bbZ[H_1(M;\bbZ)]}\Lambda_{[w_\ff]}, 
\] 
and \(HP_* ^{[\ff]}(\Gamma)\) depends only on the symplectic
isotopy class \([\ff]\) and \(\Gamma\).  
\end{cor}

\pf 
Because of the previous Theorem, it suffices to check the
corresponding statement for 
\(\textsc{Hm}^{*}(M, \grs  _\Gamma, -2r\pi[w_\ff];
\Lambda_{[w_\ff]})\). 
First, observe that \[
\overline{\textsc{Hm}}^{*}(M, \grs  _\Gamma, 0; \Lambda_{[w_\ff]})=\overline{HM}^{*}(M, \grs  _\Gamma, c_b;\Lambda_{[w_\ff]})=0,\]
due to the well-known fact that the Novikov homology of a
torus is trivial. Thus, we may again write \(\widehat{\textsc{Hm}}^{*}(M, \grs  _\Gamma, 0;
\Lambda_{[w_\ff]})=\widecheck{\textsc{Hm}}^{*}(M, \grs  _\Gamma, 0;
\Lambda_{[w_\ff]})=\textsc{Hm}^{*}(M, \grs  _\Gamma, 0;
\Lambda_{[w_\ff]})\). Now, the variant of Theorem 31.5.1 in \cite{KM}
with local coefficients states:
\[
\textsc{Hm}^{*}(M, \grs  _\Gamma, 0;
\Lambda_{[w_\ff]})\simeq \textsc{Hm}^{*}(M, \grs  _\Gamma, -2r\pi[w_\ff];
\Lambda_{[w_\ff]}).
\]
However, the coefficient ring \(\bbZ[H_1(M;\bbZ)]\) is
\((\grs_\Gamma,0)\)-complete, and \(\Lambda_{[w_\ff]}\) is flat over \(\bbZ[H_1(M;\bbZ)]\).
Thus 
\[
\textsc{Hm}^{*}(M, \grs  _\Gamma, 0;
\Lambda_{[w_\ff]})=\textsc{Hm}^{*}(M, \grs  _\Gamma, 0;\bbZ[H_1(M;\bbZ)])\otimes_{\bbZ[H_1(M;\bbZ)]}
\Lambda_{[w_\ff]}.
\]

For homologies classes \(\Gamma\) with degree 1, the
Periodic Floer homology reduces to the Floer homology of symplectic
fixed points. Thus, in this special case, Corollaries \ref{cor:PFH-inv},
\ref{cor:PFH-inv2} reduce to the invariance properties of the latter 
Floer homology in the existent literature.

\section{Appendices}

In Appendix A, we describe the appropriate parameter
space of the periodic Floer homology, and clarify its relation with
the parameter space of the Seiberg-Witten-Floer cohomology.
Note from the outset that in general, the coefficient system \(\Lambda\) will not be
considered to be among the parameters of a Floer theory
since it is chosen extraneously.
Moreover, it is clear from the discussion in Section \ref{sec:6} that 
local systems for two different Floer theories may be identified by
identifying their respective versions of \(\underline{\grG}\), which
depend on the intrinsic parameters of the theory. So do the functions
\(\mathpzc{s}\), \(\mathpzc{e}_\mu\) that determine the
completeness condition of the local systems. 

One application of the isomorphism theorem in this article is the
computation of Seiberg-Witten-Floer cohomology for degree 1
\(\Spin^c\) structures of all mapping tori, via the computation of
Floer homology of symplectic fixed points for oriented surfaces 
in the existent literature. We denote this latter Floer homology by \(HF\).  
In Appendix B, we describe the version of \(HF\)
which is equivalent to the periodic Floer homology. In \cite{C} and
its sequel, Cotton-Clay has a more refined definition of \(HF\),
which may be used to compute the version described below.

\subsection*{A. The dependence of \(HP\) on \(\ff\)} 

\setcounter{equation}{0}

Let \(h \in   \op{Symp}(F, w_F)\) be an element in the group of
symplectomorphisms of \((F, w_F)\), and let \(\op{Symp}_h(F,
w_F)\subset \op{Symp}(F, w_F)\) denote the path component containing
\(h\). 
We now define
family versions of \(M\), \(w_\ff\), and \(c_1(K^{-1})\) parametrized
by \(\op{Symp}_h(F, w_F)\) and its universal covering.

Let \(\underline{M}\to \op{Symp}_h(F, w_F)\) be defined as follows: 
\[
\underline{M}=\big(\op{Symp}_h(F, w_F)\times \bbR\times F\big)/\big((\ff, t,
x)\sim(f, t+2\pi, \ff(x))\big).
\] 
By projecting to the first factor of the product \(\op{Symp}_h(F,
w_F)\times \bbR\times F\), \(\underline{M}\) may be viewed as a bundle
over \(\op{Symp}_h(F,w_F)\) whose fiber over \(\ff\) is the mapping
torus \(M_\ff\) of \(\ff\). By projecting to the first two factors,
\(\underline{M}\) may be viewed as a fiber bundle over \(\op{Symp}_h(F,
w_F)\times S^1\); and the Euler class of the vertical tangent bundle
\(\underline{K}^{-1}\) over this fiber bundle, denoted by \(\underline{c}\), restricts to the
class \(c_1(K^{-1})\) on each \(M_\ff\subset \underline{M}\). Finally,
one may pull-back \(w_F\) from the projection to the last factor
\(F\). This induces a closed 2-form \(\underline{w}\), which restricts
to each \(M_\ff\) as \(w_\ff\).

Choose a connection \(\mathpzc{a}\) on the bundle \(\underline{M}\to \op{Symp}_h(F,
w_F)\) which preserves the structure of its fibers \(M_\ff\) as a
mapping torus of symplectomorphisms  of \((F, w_F)\). 
Let \(X\) denote a symplectic vector field on \((F,  
w_F)\) and thus an element in \(T_\ff\op{Symp}_h(F,w_F)\).  A connection as just  
defined assigns to \(X\) the following object:
\BTitem\label{symp-conn}
\item \(X_\ff\) is a section
 of the vertical tangent bundle \(\ker \pi_*\) over the
mapping torus \(\pi\co M_\ff\to S^1\),
\item \(X_\ff\) is a map sending each \(t\in S^1\) to  a symplectic
  vector field \(X_\ff
  (t)\) on \((F, w_F)\).
\item \(\pi_* X_\ff=X\).
\ETitem
The preceding construction defines various associated bundles of
\(\underline{M}\to \op{Symp}_h(F, w_F)\) and connections on them.
Note however that the monodromies of the various homology or cohomology bundles
\(H(M)\to \underline{H}(M)\to \op{Symp}_h(F, w_F)\) do not depend on
the choice of connection on \(\underline{M}\).

Let \(\widetilde{M}\to \widetilde{\op{Symp}}_h(F, w_F)\) be the pull-back
bundle over \(\underline{M}\) under the covering map of
\(\widetilde{\op{Symp}}_h(F, w_F)\), and let \(\widetilde{H}(M)\), \(\widetilde{c},
\)\(\widetilde{w}\) respectively denote the pull-backs of \(\underline{H}(M)\),
\(\underline{c}\), \(\underline{w}\) over \(\widetilde{M}\). The choice of
connection on \(\underline{M}\) above defines a trivialization of
\(\widetilde{M}\to \widetilde{\op{Symp}}_h(F, w_F)\), which induces
trivializations on its various associated bundles as well.

In order to compare with the Seiberg-Witten-Floer cohomology, we view
the 
periodic Floer homology \(HP\) as intrinsically an invariant of  a
fibered 3-manifold \(\pi\co M\to S^1\), a homology class \(\Gamma\in
H_1(M; \bbZ)\), a  closed 2-form \(w\in \Omega^2(M)\) which restricts to
\(w_F\) on the fibers, and an almost complex structure
\(J\) on \(\bbR\times M\) satisfying the conditions (\ref{def:J}), with
\(w_\ff\) therein replaced by \(w\). 
From this point of view, \(HP\) depends indirectly on 
the symplectomorphism \(\ff\) through the choice of the two-form
\(w=w_\ff\), and through the vector field \(\partial_t\). Thus, to
compare 
the periodic Floer homologies of different symplectomorphisms, 
first fix a lifting \(\ff_0\in  \widetilde{\op{Symp}}_h(F, w_F)\) of
\(h\). Identify \(M\) with the fiber \(M_{\ff_0}\subset \widetilde{M}\),
and set \(w=w_{\ff_0}\). Then for each triple consisting of \(\Gamma\in
H_1(M_{\ff_0};\bbZ)\), a suitable local system \(\Lambda\), and \(J\in \mathcal{J}_{\ff_0}\), the construction in
Sections \ref{sec:1} and \ref{sec:6} gives rise to a periodic Floer homology, which
we denote by \(HP\, (\ff_0, \Gamma; \Lambda)_J\). For an arbitrary \(\ff\in
\widetilde{\op{Symp}}_h(F, w_F)\), use the trivializations given
by \(\mathpzc{a}\) to define isomorphisms \(\iota_\ff\co M_\ff\to
M_{\ff_0}\), \(\iota^*_\ff\co \Omega^2(M_{\ff_0})\to
\Omega^2(M_\ff)\), \((\iota_\ff)_*\co H_1(M_\ff)\to H_1(M_{\ff_0})\),
and \((\iota_\ff)_* \co \mathcal{J}_\ff\to
\mathcal{J}_{\ff_0}\). 
 
Then for each triple consisting of a \(\Gamma_\ff\in H_1(M_\ff)\), a
local system \(\Lambda_\ff\) for \(\ff\),  and a \(J_\ff\in
\mathcal{J}_\ff\), we define \(HP \, (\ff, \Gamma_\ff ; \Lambda_\ff)_{J_\ff}\) to
be the periodic Floer homology associated to \((\iota_\ff)_*\Gamma\),
\((\iota_\ff^{-1})^*w_\ff\), \((\iota_\ff^{-1})^*\Lambda_\ff\), and \((\iota_\ff)_* J\).
The completeness condition for \(\Lambda_\ff\) is provided by
pull-backs of \(\widetilde{c}|_{M_\ff}+2e_{\Gamma_\ff}\) and \([w_\ff]\) via \((\iota_\ff^{-1})^*\).

In this manner, we have defined the periodic Floer homology over the
parameter space
\(\widetilde{H}_1(M;\bbZ)\times_{\widetilde{\op{Symp}}_h}\widetilde{\mathcal{J}}\). 
Here, \(\widetilde{\mathcal{J}}\) is the family version of
\(\mathcal{J}_\ff\) over
\(\widetilde{\op{Symp}}_h=\widetilde{\op{Symp}}_h(F, w_F)\).
On the other hand, since 
the Seiberg-Witten-Floer cohomology depends only on the \(\Spin^c\)
structure and the cohomology class of the perturbation 2-form,
one has \(HM\) defined over the parameter space
\(\Spin^c(M)\times H^2(M;\bbZ)\), where \(\Spin^c(M)\) denotes the set
of \(\Spin^c\) structures on \(M\). Fix a sufficiently large \(r\). 
The isomorphism \(\Phi\) from \(HP\)
  to \(HM\) constructed in the main theorems of this article
  factors through a map \(\mathpzc{q}\) between their parameter spaces,
  in the sense that
  \(\Phi(HP(\cdot))=HM(\mathpzc{q}(\cdot))\). To describe
  this map \(\mathpzc{q}\), first, note the aforementioned connection defines a trivialization
  \[\widetilde{H}_1(M;\bbZ)\times_{\widetilde{\op{Symp}}_h}\widetilde{\mathcal{J}}\to
  H_1(M;\bbZ)\times\widetilde{\op{Symp}}_h\times \mathcal{J}_{f_0}.\]  
The map \(\mathpzc{q}\) is the composition of this with the map which sends \((\Gamma, \ff, J)\in H_1(M;\bbZ)\times\widetilde{\op{Symp}}_h\times \mathcal{J}_{f_0}\) to 
\(\left( \grs_\Gamma,
  2r\big([w_{\ff_0}]+\mathpzc{b}(\theta_\ff)\big)+2\pi
  c_\Gamma\right)\in \Spin^c(M)\times H^2(M;\bbR)\).
In the above, \(\grs_\Gamma\) is defined as in Section \ref{sec:1(c)},
\(\theta_\ff\) denotes the flux of a symplectic isotopy from \(\ff_0\)
to \(\ff\), and the map \(\mathpzc{b}\) is part of the Mayer-Vietoris sequence:
\(\cdots H^1(F;\bbR)\stackrel{1-\ff^*}{\longrightarrow} H^1(F;\bbR) \stackrel{\mathpzc{b}}{\to}
H^2(M;\bbR)\cdots\). As a consequence, we have a family of periodic
Floer homologies parametrized by \(\Spin^c(M)\times H^2(M;\bbR)\), and 
this does not depend on the choice of connection on \(\underline{M}\).
When \(\ff\) vary through a path in \(\widetilde{\op{Symp}}_h\)
connecting two different lifts of the same symplectomorphism, the
corresponding \(\Spin^c\)-structure (determined by the vector field \(\partial_t\)) and
the class \([w_\ff]\) will change simultaneously. The periodic Floer
homology, and hence the corresponding Seiberg-Witten-Floer cohomology,
parametrized by these two elements in \(\Spin^c(M)\times H^2(M;\bbR)\)
will however be isomorphic.
 
\subsubsection*{\it Remark.} Cotton-Clay pointed out to us that there is no
such monodromy unless \(g=1\); therefore the above observation is
interesting only in very limited cases.

\subsection*{Appendix B: Relation with the
  Floer homology of fixed points}

When \(\Gamma\) has degree 1, a generator of \(CP\) may be
regarded as a section of \(\pi\co M\to S^1\), and by \cite{H1}, an
element in \(\mathcal{M}_1(\Theta_+, \Theta_-)\) used to define the
differential of the periodic Floer chain complex consists of a
holomorphic cylinder in \(\bbR\times M\). When \(J\) is admissible, 
such a cylinder may be
viewed as a holomorphic section \(\bbR\times (\bbR/(t\sim t+2\pi) \to
\bbR\times ((\bbR\times F)/\sim)=\bbR\times M\).  A section of
\(\pi\co M\to S^1\)
is by definition an element in the \(\ff\)-twisted loop space of
\(F\), which we denote by \(\scrL_\ff\). A section of \(\bbR\times
M\to \bbR\times S^1\) can be interpreted as a path in
\(\scrL_\ff\). Given an element \(\Gamma\in H_1(M; \bbZ)\), let
\(\scrL_{\ff ,\Gamma}\) denote the space of sections of \(\pi\co M\to
S^1\) in the homology class \(\Gamma\). The \(\ff\)-twisted loop space
decomposes as \(\scrL_\ff=\coprod _{\Gamma\in \scrN}\scrL_{\ff,
  \Gamma}\), where \(\scrN\) denotes the set of elements \(\Gamma\in
H_1(M;\bbZ)\) with \(\langle e_\Gamma\cdot, [F]\rangle=1\). This is an
affine space under \(H_1(F; \bbZ)/(1-f_*)\).
 
For lack of space, we shall only outline the construction of \(HF\), and refer
the reader to Sections 2 and 3 of \cite{L1} for more details. 

Fix the symplectic
structure \((F, w_F)\). The space of admissible almost complex
structures on \(\bbR\times M\) is identical with the space \(\mathcal{J}_K\)
in \cite{L1}, and an element from this space defines a metric on
\(\scrL_\ff\). Let \(\scrX_\ff\) be the space of elements satisfying the first two
conditions of (\ref{symp-conn}). 
Given an \(X_\ff\in \scrX_\ff\), Equation (13) of \cite{L1} defines an
action 1-form on \(\scrL_\ff\). 
One may use the formal Morse-Novikov theory of
this 1-form over \(\scrL_{\ff, \Gamma}\) and the aforementioned metric to
define a Floer homology for a suitable choice of coefficient \(\Lambda\). We 
denote this Floer homology by \(HF(\ff, \Gamma,  X_\ff;
\Lambda)_J\). The critical points of this action 1-form consists of
\(\ff\)-twisted periodic orbits of the symplectic vector fields
\(X_\ff\), and the flow lines satisfy the perturbed Cauchy-Riemann
equation. In particular, when \(X_\ff=0\), we have an isomorphism
between the \(HF\) chain complexes and the \(HP\) chain
complexes for the same triple \(\ff, \Gamma, J\). Furthermore, in the
language of Section \ref{sec:local-sys}, there is a
natural functor from the fundamental groupoid of \(\scrL_{\ff, \Gamma}\)
to the version of \(\underline{\grG}\) in the periodic Floer homology,
whose morphism sets are \(H_2(M; \Theta_+, \Theta_-)\). (This is an
variant of the map \(\underline{\op{im}}\) defined in Section 3.1.3 of
\cite{L1}, and hence we denote it by the same notation). The functors \(\scrS\), \(\scrE_\mu\) in \(HF\) (defined from
the Conley-Zehnder index and the aforementioned action 1-form) both
factor through \(\underline{\op{im}}\). Thus, by restricting the choice of local
systems to those which factor through this map, we may use the same
groupoid as the \(\underline{\grG}\) in \(HF\) as well. To
continue the comparison, recall the basic fact (\ref{pi-s}) on \(HP\).
When \(X_\ff=0\), the
\(HP\)-version of \(\mathpzc{e}_\mu\), i.e. \([w_\ff]\),
pulls back by
\(\underline{\op{im}}\) to its counterpart in \(HF\). However, 
the \(HF\)-version of \(\mathpzc{s}\), which pulls back to \(c_1(K^{-1})\), is 
seemingly different from its counterpart in \(HP\), namely 
\(c_\Gamma=c_1(K^{-1})+2e_\Gamma\). However, this is inconsequential. The fact
that in the degree 1 case, \(\mathcal{M}_1(\Theta_+, \Theta_-)\)
consists of holomorphic cylinders implies that in both \(HF\) and
\(HP\) of degree 1, the relative homology classes of the flow
lines used to define their differentials fall in the image of
\(\underline{\op{im}}\). Let \(\grH\subset H_2(M;\bbZ)\) denote the image
of \(\underline{\op{im}}\) as in \cite{L1}. Then \(e_\Gamma|\grH=0\).
 
This has the following implications: First, the relative grading of
the degree 1 periodic Floer homology can be refined to take value in
\(\bbZ/p_1\bbZ\), where \(p_1\geq p\) is the divisibility of \(c_1(K^{-1})\)
on \(\grH/\op{Tors}\).  Second, the \((\mathpzc{s},
\mathpzc{e}_\mu)\)-completeness condition on the coefficient systems
is weaker in the degree 1 case, since \(\mathpzc{s},
\mathpzc{e}_\mu\) can now be viewed as a linear function from the
smaller group \(\grH\subset H_2(M;\bbZ)\). In particular, the monotonicity
condition (\ref{mono-balanced}) requires only that \([w_\ff]|_\grH\)
is proportional to \(c_1(K^{-1})|\grH\). 
For example, if \(\ff \) is monotone in the sense that \([w_\ff]\) and
\(c_1(K^{-1})\) are proportional, then the periodic Floer homology with
\(\bbZ\)-coefficients is well-defined for all \(\Gamma\in \scrN\) at
this \(\ff\), and thus there is a well-defined notion of
\(HP_1(\ff; \bbZ)=\bigoplus_{\Gamma\in \scrN} HP\, (\ff, \Gamma;\bbZ)\),
agreeing with the usual notion of \(HF\, (\ff; \bbZ)\). (Cp. Theorem 4.1.3 of \cite{L1}).

Finally, we compare the parameter spaces of these two Floer theories. 
In contrast to the view point we adopted for \(HP\) in Appendix A, we view the
symplectomorphism \(\ff\) as an intrinsic parameter of \(HF\). 
The space
\(\scrP_{HF}=\underline{\scrN}\times_{\op{Symp}_h}\underline{\scrX}\times_{\op{Symp}_h}\underline{\mathcal{J}}\)
parametrizes the above construction of \(HF\). Here,
\(\underline{\scrN}\), \(\underline{\scrX}\), \(\underline{\mathcal{J}}\)
are respectively the family versions of \(\scrN\),
\(\scrX\), \(\mathcal{J}_K\) parametrized by \(\op{Symp}_h(F, w_F)\).
The connection \(\mathpzc{a}\) induces connections on each of these
bundles, and hence a connection on the fibered product \(\scrP_{HF}\)
as well. 

We now describe the map 
from the \(HP\)-parameter space 
to the \(HF\)-parameter space over which the 
isomorphism between the two Floer homologies is defined. 
When restricted to
degree 1 \(H_1(M;\bbZ)\) classes, the \(HP\) parameter space is 
trivialized by the connection \(\mathpzc{a}\) as
\(\scrP_{HP_1}=\scrN\times\widetilde{\op{Symp}}_h\times \mathcal{J}\). As explained
above, the same connection gives a trivialization of the
pull-back of the \(HF\)-parameter space over
\(\widetilde{\op{Symp}}_h\):
\(
\widetilde{\scrP}_{HF}\simeq  \scrN\times\widetilde{\op{Symp}_h}\times \scrX_{\ff_0}\times\mathcal{J}.
\) (Cf. Section 3.1.5 in \cite{L1}, where it is also shown that parallel
transport along the connection on \(\scrP_{HF}\) yields
isomorphisms among  the \(HF\) Floer chain complexes so parametrized).

Let \(\widetilde{\mathpzc{m}}=\op{id}\times \widetilde{s}_0\times\op{id}\), where \(\widetilde{s}_0\) denotes the
zero section of \(\widetilde{\op{Symp}_h}\times \scrX_{\ff_0}\to
\widetilde{\op{Symp}_h}\).  (Note that the zero section
\(\widetilde{s_0}\) is nonconstant relative to the trivilization induced
by \(\mathpzc{a}\), and depends on the choice of the lift \(\ff_0\)).
This is the map underlying the
isomorphism between the parametrized versions of \(HP_1\) and
\(HF\). Notice that by the invariance of \(HF\) under
the parallel transport determined by \( \mathpzc{a}\), one may identify \(HF\) on the zero
section with a fiber \(\scrN\times \{\ff_0\}\times \scrX_{\ff_0}\times
\mathcal{J}\) of \(\widetilde{\scrP}_{HF}\). As a consequence, one may take
\(\mathpzc{m}\co
\scrN\times\widetilde{\op{Symp}}_h\times \mathcal{J}\to \scrN\times \scrX_{\ff_0}\times
\mathcal{J}\) to be the map underlying the isomorphism
between the parametrized versions of \(HF\) and \(HP_1\), where  
\(\mathpzc{m}\) sends \((n, \ff, J)\) to \((n, \theta_\ff, J)\), where
\(\ff\mapsto \theta_\ff\) is as in the definition of
\(\mathpzc{q}\) in Appendix A.

\end{document}